\documentstyle[12pt,amstex]{article}
\newtheorem{theorem}{Theorem}[subsection]
\newtheorem{corollary}[theorem]{Corollary}
\newtheorem{lemma}[theorem]{Lemma}
\newtheorem{proposition}[theorem]{Proposition}
\newtheorem{remark}[theorem]{Remark}
\newtheorem{example}[theorem]{Example}
\DeclareSymbolFont{AMSa}{U}{msa}{m}{n}
\DeclareMathSymbol{\boxtimes}{\mathbin}{AMSa}{"02}
\def\fin{\nolinebreak \raise-0.1em\hbox{$\Box$}}
\def\stackunder#1#2{\mathrel{\mathop{#2}\limits_{#1}}}
\def\dual{ ^{\rm v}}
\def\bidual{ ^{\rm vv}}
\def\primedual{ ^{\prime \rm v}}
\def\ideal{{\cal J}}
\def\LL{{\cal L}}

\def\Ext{{\cal E}xt }
\def\Tor{{\cal T}or }

\def\prefix#1{\mbox{{\footnotesize $( #1 )$}-}}
\newcommand{\complex}{\setlength{\unitlength}{0.01em}
\begin{picture}(70,70)
\put(0,0){C}
\put(27,3){\line(0,1){62}}
\end{picture}
}
\newcommand{\proj}{{\rm I} \! {\rm P}}
\newcommand{\nat}{{\rm I} \! {\rm N}}
\newcommand{\OO}{{\cal O}}
\newcommand{\Op}{{{\cal O}_{\proj_3}}}
\def\kercoker#1{\!\!\!\stackunder{
    \begin{array}{c}
        \phantom{.}\searrow \phantom{M} \nearrow \phantom{.}\\
         #1 \\
        \phantom{.}\nearrow \hfill \searrow\phantom{.} \\
       ^0\hfill ^0
    \end{array}}{\longrightarrow }\!\!\!}
\makeatletter

\def\diagram{\leftwidth=\z@ \rightwidth=\z@ \topheight=\z@
\botheight=\z@ \setbox\@picbox\hbox\bgroup}

\def\enddiagram{\egroup\wd\@picbox\rightwidth\unitlength
\ht\@picbox\topheight\unitlength \dp\@picbox\botheight\unitlength
\hskip\leftwidth\unitlength\box\@picbox}

\def\bfig{\begin{diagram}}
\def\efig{\end{diagram}}
\newcount\wideness \newcount\leftwidth \newcount\rightwidth
\newcount\highness \newcount\topheight \newcount\botheight

\def\ratchet#1#2{\ifnum#1<#2 \global #1=#2 \fi}

\def\putbox(#1,#2)#3{%
\horsize{\wideness}{#3} \divide\wideness by 2
{\advance\wideness by #1 \ratchet{\rightwidth}{\wideness}}
{\advance\wideness by -#1 \ratchet{\leftwidth}{\wideness}}
\vertsize{\highness}{#3} \divide\highness by 2
{\advance\highness by #2 \ratchet{\topheight}{\highness}}
{\advance\highness by -#2 \ratchet{\botheight}{\highness}}
\put(#1,#2){\makebox(0,0){$#3$}}}

\def\putlbox(#1,#2)#3{%
\horsize{\wideness}{#3}
{\advance\wideness by #1 \ratchet{\rightwidth}{\wideness}}
{\ratchet{\leftwidth}{-#1}}
\vertsize{\highness}{#3} \divide\highness by 2
{\advance\highness by #2 \ratchet{\topheight}{\highness}}
{\advance\highness by -#2 \ratchet{\botheight}{\highness}}
\put(#1,#2){\makebox(0,0)[l]{$#3$}}}

\def\putrbox(#1,#2)#3{%
\horsize{\wideness}{#3}
{\ratchet{\rightwidth}{#1}}
{\advance\wideness by -#1 \ratchet{\leftwidth}{\wideness}}
\vertsize{\highness}{#3} \divide\highness by 2
{\advance\highness by #2 \ratchet{\topheight}{\highness}}
{\advance\highness by -#2 \ratchet{\botheight}{\highness}}
\put(#1,#2){\makebox(0,0)[r]{$#3$}}}

\def\adjust[#1]{} 

\newcount \coefa
\newcount \coefb
\newcount \coefc
\newcount\tempcounta
\newcount\tempcountb
\newcount\tempcountc
\newcount\tempcountd
\newcount\xext
\newcount\yext
\newcount\xoff
\newcount\yoff
\newcount\gap%
\newcount\arrowtypea
\newcount\arrowtypeb
\newcount\arrowtypec
\newcount\arrowtyped
\newcount\arrowtypee
\newcount\height
\newcount\width
\newcount\xpos
\newcount\ypos
\newcount\run
\newcount\rise
\newcount\arrowlength
\newcount\halflength
\newcount\arrowtype
\newdimen\tempdimen
\newdimen\xlen
\newdimen\ylen
\newsavebox{\tempboxa}%
\newsavebox{\tempboxb}%
\newsavebox{\tempboxc}%

\newdimen\w@dth

\def\setw@dth#1#2{\setbox\z@\hbox{$#1$}\w@dth=\wd\z@
\setbox\@ne\hbox{$#2$}\ifnum\w@dth<\wd\@ne \w@dth=\wd\@ne \fi
\advance\w@dth by 1.2em}

\def\t@^#1_#2{\def\n@one{#1}\def\n@two{#2}\mathrel{\setw@dth{#1}{#2}
\mathop{\hbox to \w@dth{\rightarrowfill}}\limits
\ifx\n@one\empty\else ^{\box\z@}\fi
\ifx\n@two\empty\else _{\box\@ne}\fi}}
\def\t@@^#1{\@ifnextchar_ {\t@^{#1}}{\t@^{#1}_{}}}
\def\to{\@ifnextchar^ {\t@@}{\t@@^{}}}

\def\t@left^#1_#2{\def\n@one{#1}\def\n@two{#2}\mathrel{\setw@dth{#1}{#2}
\mathop{\hbox to \w@dth{\leftarrowfill}}\limits
\ifx\n@one\empty\else ^{\box\z@}\fi
\ifx\n@two\empty\else _{\box\@ne}\fi}}
\def\t@@left^#1{\@ifnextchar_ {\t@left^{#1}}{\t@left^{#1}_{}}}
\def\toleft{\@ifnextchar^ {\t@@left}{\t@@left^{}}}

\def\two@^#1_#2{\def\n@one{#1}\def\n@two{#2}\mathrel{\setw@dth{#1}{#2}
\mathop{\vcenter{\hbox to \w@dth{\rightarrowfill}\kern-1.7ex
                 \hbox to \w@dth{\rightarrowfill}}%
       }\limits
\ifx\n@one\empty\else ^{\box\z@}\fi
\ifx\n@two\empty\else _{\box\@ne}\fi}}
\def\tw@@^#1{\@ifnextchar_ {\two@^{#1}}{\two@^{#1}_{}}}
\def\two{\@ifnextchar^ {\tw@@}{\tw@@^{}}}

\def\tofr@^#1_#2{\def\n@one{#1}\def\n@two{#2}\mathrel{\setw@dth{#1}{#2}
\mathop{\vcenter{\hbox to \w@dth{\rightarrowfill}\kern-1.7ex
                 \hbox to \w@dth{\leftarrowfill}}%
       }\limits
\ifx\n@one\empty\else ^{\box\z@}\fi
\ifx\n@two\empty\else _{\box\@ne}\fi}}
\def\t@fr@^#1{\@ifnextchar_ {\tofr@^{#1}}{\tofr@^{#1}_{}}}
\def\tofro{\@ifnextchar^ {\t@fr@}{\t@fr@^{}}}

\def\mon{\mathop{\m@th\hbox to
      14.6\P@{\lasyb\char'51\hskip-2.1\P@$\arrext$\hss
$\mathord\rightarrow$}}\limits} 
\def\leftmono{\mathrel{\m@th\hbox to
14.6\P@{$\mathord\leftarrow$\hss$\arrext$\hskip-2.1\P@\lasyb\char'50%
}}\limits} 
\mathchardef\arrext="0200       

\def\settypes(#1,#2,#3){\arrowtypea#1 \arrowtypeb#2 \arrowtypec#3}
\def\settoheight#1#2{\setbox\@tempboxa\hbox{#2}#1\ht\@tempboxa\relax}%
\def\settodepth#1#2{\setbox\@tempboxa\hbox{#2}#1\dp\@tempboxa\relax}%
\def\settokens[#1`#2`#3`#4]{%
     \def\tokena{#1}\def\tokenb{#2}\def\tokenc{#3}\def\tokend{#4}}
\def\setsqparms[#1`#2`#3`#4;#5`#6]{%
\arrowtypea #1
\arrowtypeb #2
\arrowtypec #3
\arrowtyped #4
\width #5
\height #6
}
\def\setpos(#1,#2){\xpos=#1 \ypos#2}

\def\settriparms[#1`#2`#3;#4]{\settripairparms[#1`#2`#3`1`1;#4]}%

\def\settripairparms[#1`#2`#3`#4`#5;#6]{%
\arrowtypea #1
\arrowtypeb #2
\arrowtypec #3
\arrowtyped #4
\arrowtypee #5
\width #6
\height #6
}

\def\resetparms{\settripairparms[1`1`1`1`1;500]\width 500}

\resetparms

\def\mvector(#1,#2)#3{
\put(0,0){\vector(#1,#2){#3}}%
\put(0,0){\vector(#1,#2){26}}%
}
\def\evector(#1,#2)#3{{
\arrowlength #3
\put(0,0){\vector(#1,#2){\arrowlength}}%
\advance \arrowlength by-30
\put(0,0){\vector(#1,#2){\arrowlength}}%
}}

\def\horsize#1#2{%
\settowidth{\tempdimen}{$#2$}%
#1=\tempdimen
\divide #1 by\unitlength
}

\def\vertsize#1#2{%
\settoheight{\tempdimen}{$#2$}%
#1=\tempdimen
\settodepth{\tempdimen}{$#2$}%
\advance #1 by\tempdimen
\divide #1 by\unitlength
}

\def\putvector(#1,#2)(#3,#4)#5#6{{%
\ifnum3<\arrowtype
\putdashvector(#1,#2)(#3,#4)#5\arrowtype
\else
\ifnum\arrowtype<-3
\putdashvector(#1,#2)(#3,#4)#5\arrowtype
\else
\xpos=#1
\ypos=#2
\run=#3
\rise=#4
\arrowlength=#5
\ifnum \arrowtype<0
    \ifnum \run=0
        \advance \ypos by-\arrowlength
    \else
        \tempcounta \arrowlength
        \multiply \tempcounta by\rise
        \divide \tempcounta by\run
        \ifnum\run>0
            \advance \xpos by\arrowlength
            \advance \ypos by\tempcounta
        \else
            \advance \xpos by-\arrowlength
            \advance \ypos by-\tempcounta
        \fi
    \fi
    \multiply \arrowtype by-1
    \multiply \rise by-1
    \multiply \run by-1
\fi
\ifcase \arrowtype
\or \put(\xpos,\ypos){\vector(\run,\rise){\arrowlength}}%
\or \put(\xpos,\ypos){\mvector(\run,\rise)\arrowlength}%
\or \put(\xpos,\ypos){\evector(\run,\rise){\arrowlength}}%
\fi\fi\fi
}}

\def\putsplitvector(#1,#2)#3#4{
\xpos #1
\ypos #2
\arrowtype #4
\halflength #3
\arrowlength #3
\gap 140
\advance \halflength by-\gap
\divide \halflength by2
\ifnum\arrowtype>0
   \ifcase \arrowtype
   \or \put(\xpos,\ypos){\line(0,-1){\halflength}}%
       \advance\ypos by-\halflength
       \advance\ypos by-\gap
       \put(\xpos,\ypos){\vector(0,-1){\halflength}}%
   \or \put(\xpos,\ypos){\line(0,-1)\halflength}%
       \put(\xpos,\ypos){\vector(0,-1)3}%
       \advance\ypos by-\halflength
       \advance\ypos by-\gap
       \put(\xpos,\ypos){\vector(0,-1){\halflength}}%
   \or \put(\xpos,\ypos){\line(0,-1)\halflength}%
       \advance\ypos by-\halflength
       \advance\ypos by-\gap
       \put(\xpos,\ypos){\evector(0,-1){\halflength}}%
   \fi
\else \arrowtype=-\arrowtype
   \ifcase\arrowtype
   \or \advance \ypos by-\arrowlength
       \put(\xpos,\ypos){\line(0,1){\halflength}}%
       \advance\ypos by\halflength
       \advance\ypos by\gap
       \put(\xpos,\ypos){\vector(0,1){\halflength}}%
   \or \advance \ypos by-\arrowlength
       \put(\xpos,\ypos){\line(0,1)\halflength}%
       \put(\xpos,\ypos){\vector(0,1)3}%
       \advance\ypos by\halflength
       \advance\ypos by\gap
       \put(\xpos,\ypos){\vector(0,1){\halflength}}%
   \or \advance \ypos by-\arrowlength
       \put(\xpos,\ypos){\line(0,1)\halflength}%
       \advance\ypos by\halflength
       \advance\ypos by\gap
       \put(\xpos,\ypos){\evector(0,1){\halflength}}%
   \fi
\fi
}

\def\putmorphism(#1)(#2,#3)[#4`#5`#6]#7#8#9{{%
\run #2
\rise #3
\ifnum\rise=0
  \puthmorphism(#1)[#4`#5`#6]{#7}{#8}#9%
\else\ifnum\run=0
  \putvmorphism(#1)[#4`#5`#6]{#7}{#8}#9%
\else
\setpos(#1)%
\arrowlength #7
\arrowtype #8
\ifnum\run=0
\else\ifnum\rise=0
\else
\ifnum\run>0
    \coefa=1
\else
   \coefa=-1
\fi
\ifnum\arrowtype>0
   \coefb=0
   \coefc=-1
\else
   \coefb=\coefa
   \coefc=1
   \arrowtype=-\arrowtype
\fi
\width=2
\multiply \width by\run
\divide \width by\rise
\ifnum \width<0  \width=-\width\fi
\advance\width by60
\if l#9 \width=-\width\fi
\putbox(\xpos,\ypos){#4}
{\multiply \coefa by\arrowlength
\advance\xpos by\coefa
\multiply \coefa by\rise
\divide \coefa by\run
\advance \ypos by\coefa
\putbox(\xpos,\ypos){#5} }%
{\multiply \coefa by\arrowlength
\divide \coefa by2
\advance \xpos by\coefa
\advance \xpos by\width
\multiply \coefa by\rise
\divide \coefa by\run
\advance \ypos by\coefa
\if l#9%
   \putrbox(\xpos,\ypos){#6}%
\else\if r#9%
   \putlbox(\xpos,\ypos){#6}%
\fi\fi }%
{\multiply \rise by-\coefc
\multiply \run by-\coefc
\multiply \coefb by\arrowlength
\advance \xpos by\coefb
\multiply \coefb by\rise
\divide \coefb by\run
\advance \ypos by\coefb
\multiply \coefc by70
\advance \ypos by\coefc
\multiply \coefc by\run
\divide \coefc by\rise
\advance \xpos by\coefc
\multiply \coefa by140
\multiply \coefa by\run
\divide \coefa by\rise
\advance \arrowlength by\coefa
\ifcase\arrowtype
\or \put(\xpos,\ypos){\vector(\run,\rise){\arrowlength}}%
\or \put(\xpos,\ypos){\mvector(\run,\rise){\arrowlength}}%
\or \put(\xpos,\ypos){\evector(\run,\rise){\arrowlength}}%
\fi}\fi\fi\fi\fi}}

\newcount\numbdashes \newcount\lengthdash \newcount\increment

\def\howmanydashes{
\numbdashes=\arrowlength \lengthdash=40
\divide\numbdashes by \lengthdash
\lengthdash=\arrowlength
\divide\lengthdash by \numbdashes
\increment=\lengthdash
\multiply\lengthdash by 3
\divide\lengthdash by 5
}

\def\putdashvector(#1)(#2,#3)#4#5{%
\ifnum#3=0 \putdashhvector(#1){#4}#5
\else
\ifnum#2=0
\putdashvvector(#1){#4}#5\fi\fi}

\def\putdashhvector(#1,#2)#3#4{{%
\arrowlength=#3 \howmanydashes
\multiput(#1,#2)(\increment,0){\numbdashes}%
{\vrule height .4pt width \lengthdash\unitlength}
\arrowtype=#4 \xpos=#1
\ifnum\arrowtype<0 \advance\arrowtype by 7 \fi
\ifcase\arrowtype
\or \advance\xpos by 10
    \put(\xpos,#2){\vector(-1,0){\lengthdash}}
    \advance\xpos by 40
    \put(\xpos,#2){\vector(-1,0){\lengthdash}}
\or \advance \xpos by 10
    \put(\xpos,#2){\vector(-1,0){\lengthdash}}
    \advance\xpos by  \arrowlength
    \advance\xpos by  -50
    \put(\xpos,#2){\vector(-1,0){\lengthdash}}
\or \advance\xpos by 10
    \put(\xpos,#2){\vector(-1,0){\lengthdash}}
\or \advance\xpos by \arrowlength
    \advance\xpos by -\lengthdash
    \put(\xpos,#2){\vector(1,0){\lengthdash}}
\or {\advance\xpos by 10
    \put(\xpos,#2){\vector(1,0){\lengthdash}}}
    \advance\xpos by \arrowlength
    \advance\xpos by -\lengthdash
    \put(\xpos,#2){\vector(1,0){\lengthdash}}
\or \advance\xpos by \arrowlength
    \advance\xpos by -\lengthdash
    \put(\xpos,#2){\vector(1,0){\lengthdash}}
    \advance\xpos by -40
    \put(\xpos,#2){\vector(1,0){\lengthdash}}
   \fi
}}

\def\putdashvvector(#1,#2)#3#4{{%
\arrowlength=#3 \howmanydashes
\ypos=#2 \advance\ypos by -\arrowlength
\multiput(#1,#2)(0,\increment){\numbdashes}%
    {\vrule width .4pt height \lengthdash\unitlength}
\arrowtype=#4 \ypos=#2
\ifnum\arrowtype<0 \advance\arrowtype by 7 \fi
\ifcase\arrowtype
\or \advance\ypos by \arrowlength \advance\ypos by -40
    \put(#1,\ypos){\vector(0,1){\lengthdash}}
    \advance\ypos by -40
    \put(#1,\ypos){\vector(0,1){\lengthdash}}
\or \advance\ypos by 10
    \put(#1,\ypos){\vector(0,1){\lengthdash}}
    \advance\ypos by \arrowlength \advance\ypos by -40
    \put(#1,\ypos){\vector(0,1){\lengthdash}}
\or \advance\ypos by \arrowlength \advance\ypos by -40
    \put(#1,\ypos){\vector(0,1){\lengthdash}}
\or \advance\ypos by 10
    \put(#1,\ypos){\vector(0,-1){\lengthdash}}
\or \advance\ypos by 10
    \put(#1,\ypos){\vector(0,-1){\lengthdash}}
    \advance\ypos by \arrowlength \advance\ypos by -40
    \put(#1,\ypos){\vector(0,-1){\lengthdash}}
\or \advance\ypos by 10
    \put(#1,\ypos){\vector(0,-1){\lengthdash}}
    \advance\ypos by 40
    \put(#1,\ypos){\vector(0,-1){\lengthdash}}
\fi
}}

\def\puthmorphism(#1,#2)[#3`#4`#5]#6#7#8{{%
\xpos #1
\ypos #2
\width #6
\arrowlength #6
\arrowtype=#7
\putbox(\xpos,\ypos){#3\vphantom{#4}}%
{\advance \xpos by\arrowlength
\putbox(\xpos,\ypos){\vphantom{#3}#4}}%
\horsize{\tempcounta}{#3}%
\horsize{\tempcountb}{#4}%
\divide \tempcounta by2
\divide \tempcountb by2
\advance \tempcounta by30
\advance \tempcountb by30
\advance \xpos by\tempcounta
\advance \arrowlength by-\tempcounta
\advance \arrowlength by-\tempcountb
\putvector(\xpos,\ypos)(1,0)\arrowlength\arrowtype
\divide \arrowlength by2
\advance \xpos by\arrowlength
\vertsize{\tempcounta}{#5}%
\divide\tempcounta by2
\advance \tempcounta by20
\if a#8 %
   \advance \ypos by\tempcounta
   \putbox(\xpos,\ypos){#5}%
\else
   \advance \ypos by-\tempcounta
   \putbox(\xpos,\ypos){#5}%
\fi}}

\def\putvmorphism(#1,#2)[#3`#4`#5]#6#7#8{{%
\xpos #1
\ypos #2
\arrowlength #6
\arrowtype #7
\settowidth{\xlen}{$#5$}%
\putbox(\xpos,\ypos){#3}%
{\advance \ypos by-\arrowlength
\putbox(\xpos,\ypos){#4}}%
{\advance\arrowlength by-140
\advance \ypos by-70
\ifdim\xlen>0pt
   \if m#8%
      \putsplitvector(\xpos,\ypos)\arrowlength\arrowtype
   \else
   \putvector(\xpos,\ypos)(0,-1)\arrowlength\arrowtype
   \fi
\else
   \putvector(\xpos,\ypos)(0,-1)\arrowlength\arrowtype
\fi}%
\ifdim\xlen>0pt
   \divide \arrowlength by2
   \advance\ypos by-\arrowlength
   \if l#8%
      \advance \xpos by-40
      \putrbox(\xpos,\ypos){#5}%
   \else\if r#8%
      \advance \xpos by40
      \putlbox(\xpos,\ypos){#5}%
   \else
      \putbox(\xpos,\ypos){#5}%
   \fi\fi
\fi
}}

\def\putsquarep<#1>(#2)[#3;#4`#5`#6`#7]{{%
\setsqparms[#1]%
\setpos(#2)%
\settokens[#3]%
\puthmorphism(\xpos,\ypos)[\tokenc`\tokend`{#7}]{\width}{\arrowtyped}b%
\advance\ypos by \height
\puthmorphism(\xpos,\ypos)[\tokena`\tokenb`{#4}]{\width}{\arrowtypea}a%
\putvmorphism(\xpos,\ypos)[``{#5}]{\height}{\arrowtypeb}l%
\advance\xpos by \width
\putvmorphism(\xpos,\ypos)[``{#6}]{\height}{\arrowtypec}r%
}}

\def\putsquare{\@ifnextchar <{\putsquarep}{\putsquarep%
   <\arrowtypea`\arrowtypeb`\arrowtypec`\arrowtyped;\width`\height>}}
\def\square{\@ifnextchar< {\squarep}{\squarep
   <\arrowtypea`\arrowtypeb`\arrowtypec`\arrowtyped;\width`\height>}}
\def\squarep<#1>[#2`#3`#4`#5;#6`#7`#8`#9]{{
\setsqparms[#1]
\diagram
\putsquarep<\arrowtypea`\arrowtypeb`\arrowtypec`
\arrowtyped;\width`\height>
(0,0)[#2`#3`#4`{#5};#6`#7`#8`{#9}]
\enddiagram
}}                                                 
\def\putptrianglep<#1>(#2,#3)[#4`#5`#6;#7`#8`#9]{{%
\settriparms[#1]%
\xpos=#2 \ypos=#3
\advance\ypos by \height
\puthmorphism(\xpos,\ypos)[#4`#5`{#7}]{\height}{\arrowtypea}a%
\putvmorphism(\xpos,\ypos)[`#6`{#8}]{\height}{\arrowtypeb}l%
\advance\xpos by\height
\putmorphism(\xpos,\ypos)(-1,-1)[``{#9}]{\height}{\arrowtypec}r%
}}

\def\putptriangle{\@ifnextchar <{\putptrianglep}{\putptrianglep
   <\arrowtypea`\arrowtypeb`\arrowtypec;\height>}}
\def\ptriangle{\@ifnextchar <{\ptrianglep}{\ptrianglep
   <\arrowtypea`\arrowtypeb`\arrowtypec;\height>}}
\def\ptrianglep<#1>[#2`#3`#4;#5`#6`#7]{{
\settriparms[#1]
\diagram
\putptrianglep<\arrowtypea`\arrowtypeb`
\arrowtypec;\height>
(0,0)[#2`#3`#4;#5`#6`{#7}]
\enddiagram
}}                                            

\def\putqtrianglep<#1>(#2,#3)[#4`#5`#6;#7`#8`#9]{{%
\settriparms[#1]%
\xpos=#2 \ypos=#3
\advance\ypos by\height
\puthmorphism(\xpos,\ypos)[#4`#5`{#7}]{\height}{\arrowtypea}a%
\putmorphism(\xpos,\ypos)(1,-1)[``{#8}]{\height}{\arrowtypeb}l%
\advance\xpos by\height
\putvmorphism(\xpos,\ypos)[`#6`{#9}]{\height}{\arrowtypec}r%
}}

\def\putqtriangle{\@ifnextchar <{\putqtrianglep}{\putqtrianglep
   <\arrowtypea`\arrowtypeb`\arrowtypec;\height>}}
\def\qtriangle{\@ifnextchar <{\qtrianglep}{\qtrianglep
   <\arrowtypea`\arrowtypeb`\arrowtypec;\height>}}
\def\qtrianglep<#1>[#2`#3`#4;#5`#6`#7]{{
\settriparms[#1]
\width=\height                                
\diagram
\putqtrianglep<\arrowtypea`\arrowtypeb`
\arrowtypec;\height>
(0,0)[#2`#3`#4;#5`#6`{#7}]
\enddiagram
}}

\def\putdtrianglep<#1>(#2,#3)[#4`#5`#6;#7`#8`#9]{{%
\settriparms[#1]%
\xpos=#2 \ypos=#3
\puthmorphism(\xpos,\ypos)[#5`#6`{#9}]{\height}{\arrowtypec}b%
\advance\xpos by \height \advance\ypos by\height
\putmorphism(\xpos,\ypos)(-1,-1)[``{#7}]{\height}{\arrowtypea}l%
\putvmorphism(\xpos,\ypos)[#4``{#8}]{\height}{\arrowtypeb}r%
}}

\def\putdtriangle{\@ifnextchar <{\putdtrianglep}{\putdtrianglep
   <\arrowtypea`\arrowtypeb`\arrowtypec;\height>}}
\def\dtriangle{\@ifnextchar <{\dtrianglep}{\dtrianglep
   <\arrowtypea`\arrowtypeb`\arrowtypec;\height>}}
\def\dtrianglep<#1>[#2`#3`#4;#5`#6`#7]{{
\settriparms[#1]
\width=\height                                
\diagram
\putdtrianglep<\arrowtypea`\arrowtypeb`
\arrowtypec;\height>
(0,0)[#2`#3`#4;#5`#6`{#7}]
\enddiagram
}}

\def\putbtrianglep<#1>(#2,#3)[#4`#5`#6;#7`#8`#9]{{%
\settriparms[#1]%
\xpos=#2 \ypos=#3
\puthmorphism(\xpos,\ypos)[#5`#6`{#9}]{\height}{\arrowtypec}b%
\advance\ypos by\height
\putmorphism(\xpos,\ypos)(1,-1)[``{#8}]{\height}{\arrowtypeb}r%
\putvmorphism(\xpos,\ypos)[#4``{#7}]{\height}{\arrowtypea}l%
}}

\def\putbtriangle{\@ifnextchar <{\putbtrianglep}{\putbtrianglep
   <\arrowtypea`\arrowtypeb`\arrowtypec;\height>}}
\def\btriangle{\@ifnextchar <{\btrianglep}{\btrianglep
   <\arrowtypea`\arrowtypeb`\arrowtypec;\height>}}
\def\btrianglep<#1>[#2`#3`#4;#5`#6`#7]{{
\settriparms[#1]
\width=\height                               
\diagram
\putbtrianglep<\arrowtypea`\arrowtypeb`
\arrowtypec;\height>
(0,0)[#2`#3`#4;#5`#6`{#7}]
\enddiagram
}}

\def\putAtrianglep<#1>(#2,#3)[#4`#5`#6;#7`#8`#9]{{%
\settriparms[#1]%
\xpos=#2 \ypos=#3
{\multiply \height by2
\puthmorphism(\xpos,\ypos)[#5`#6`{#9}]{\height}{\arrowtypec}b}%
\advance\xpos by\height \advance\ypos by\height
\putmorphism(\xpos,\ypos)(-1,-1)[#4``{#7}]{\height}{\arrowtypea}l%
\putmorphism(\xpos,\ypos)(1,-1)[``{#8}]{\height}{\arrowtypeb}r%
}}

\def\putAtriangle{\@ifnextchar <{\putAtrianglep}{\putAtrianglep
   <\arrowtypea`\arrowtypeb`\arrowtypec;\height>}}
\def\Atriangle{\@ifnextchar <{\Atrianglep}{\Atrianglep
   <\arrowtypea`\arrowtypeb`\arrowtypec;\height>}}
\def\Atrianglep<#1>[#2`#3`#4;#5`#6`#7]{{
\settriparms[#1]
\width=\height                                     
\diagram
\putAtrianglep<\arrowtypea`\arrowtypeb`
\arrowtypec;\height>
(0,0)[#2`#3`#4;#5`#6`{#7}]
\enddiagram
}}

\def\putAtrianglepairp<#1>(#2)[#3;#4`#5`#6`#7`#8]{{%
\settripairparms[#1]%
\setpos(#2)%
\settokens[#3]%
\puthmorphism(\xpos,\ypos)[\tokenb`\tokenc`{#7}]{\height}{\arrowtyped}b%
\advance\xpos by\height
\puthmorphism(\xpos,\ypos)[\phantom{\tokenc}`\tokend`{#8}]%
{\height}{\arrowtypee}b%
\advance\ypos by\height
\putmorphism(\xpos,\ypos)(-1,-1)[\tokena``{#4}]{\height}{\arrowtypea}l%
\putvmorphism(\xpos,\ypos)[``{#5}]{\height}{\arrowtypeb}m%
\putmorphism(\xpos,\ypos)(1,-1)[``{#6}]{\height}{\arrowtypec}r%
}}

\def\putAtrianglepair{\@ifnextchar <{\putAtrianglepairp}{\putAtrianglepairp%
   <\arrowtypea`\arrowtypeb`\arrowtypec`\arrowtyped`\arrowtypee;\height>}}
\def\Atrianglepair{\@ifnextchar <{\Atrianglepairp}{\Atrianglepairp%
   <\arrowtypea`\arrowtypeb`\arrowtypec`\arrowtyped`\arrowtypee;\height>}}

\def\Atrianglepairp<#1>[#2;#3`#4`#5`#6`#7]{{
\settripairparms[#1]
\settokens[#2]
\width=\height                                
\diagram
\putAtrianglepairp                            
<\arrowtypea`\arrowtypeb`\arrowtypec`
\arrowtyped`\arrowtypee;\height>
(0,0)[{#2};#3`#4`#5`#6`{#7}]
\enddiagram
}}

\def\putVtrianglep<#1>(#2,#3)[#4`#5`#6;#7`#8`#9]{{%
\settriparms[#1]%
\xpos=#2 \ypos=#3
\advance\ypos by\height
{\multiply\height by2
\puthmorphism(\xpos,\ypos)[#4`#5`{#7}]{\height}{\arrowtypea}a}%
\putmorphism(\xpos,\ypos)(1,-1)[`#6`{#8}]{\height}{\arrowtypeb}l%
\advance\xpos by\height
\advance\xpos by\height
\putmorphism(\xpos,\ypos)(-1,-1)[``{#9}]{\height}{\arrowtypec}r%
}}

\def\putVtriangle{\@ifnextchar <{\putVtrianglep}{\putVtrianglep
   <\arrowtypea`\arrowtypeb`\arrowtypec;\height>}}
\def\Vtriangle{\@ifnextchar <{\Vtrianglep}{\Vtrianglep
   <\arrowtypea`\arrowtypeb`\arrowtypec;\height>}}
\def\Vtrianglep<#1>[#2`#3`#4;#5`#6`#7]{{
\settriparms[#1]
\width=\height                                 
\diagram
\putVtrianglep<\arrowtypea`\arrowtypeb`
\arrowtypec;\height>
(0,0)[#2`#3`#4;#5`#6`{#7}]
\enddiagram
}}

\def\putVtrianglepairp<#1>(#2)[#3;#4`#5`#6`#7`#8]{{
\settripairparms[#1]%
\setpos(#2)%
\settokens[#3]%
\advance\ypos by\height
\putmorphism(\xpos,\ypos)(1,-1)[`\tokend`{#6}]{\height}{\arrowtypec}l%
\puthmorphism(\xpos,\ypos)[\tokena`\tokenb`{#4}]{\height}{\arrowtypea}a%
\advance\xpos by\height
\puthmorphism(\xpos,\ypos)[\phantom{\tokenb}`\tokenc`{#5}]%
{\height}{\arrowtypeb}a%
\putvmorphism(\xpos,\ypos)[``{#7}]{\height}{\arrowtyped}m%
\advance\xpos by\height
\putmorphism(\xpos,\ypos)(-1,-1)[``{#8}]{\height}{\arrowtypee}r%
}}

\def\putVtrianglepair{\@ifnextchar <{\putVtrianglepairp}{\putVtrianglepairp%
    <\arrowtypea`\arrowtypeb`\arrowtypec`\arrowtyped`\arrowtypee;\height>}}
\def\Vtrianglepair{\@ifnextchar <{\Vtrianglepairp}{\Vtrianglepairp%
    <\arrowtypea`\arrowtypeb`\arrowtypec`\arrowtyped`\arrowtypee;\height>}}
\def\Vtrianglepairp<#1>[#2;#3`#4`#5`#6`#7]{{
\settripairparms[#1]
\settokens[#2]
\diagram
\putVtrianglepairp                             
<\arrowtypea`\arrowtypeb`\arrowtypec`
\arrowtyped`\arrowtypee;\height>
(0,0)[{#2};#3`#4`#5`#6`{#7}]
\enddiagram
}}

\def\putCtrianglep<#1>(#2,#3)[#4`#5`#6;#7`#8`#9]{{%
\settriparms[#1]%
\xpos=#2 \ypos=#3
\advance\ypos by\height
\putmorphism(\xpos,\ypos)(1,-1)[``{#9}]{\height}{\arrowtypec}l%
\advance\xpos by\height
\advance\ypos by\height
\putmorphism(\xpos,\ypos)(-1,-1)[#4`#5`{#7}]{\height}{\arrowtypea}l%
{\multiply\height by 2
\putvmorphism(\xpos,\ypos)[`#6`{#8}]{\height}{\arrowtypeb}r}%
}}

\def\putCtriangle{\@ifnextchar <{\putCtrianglep}{\putCtrianglep
    <\arrowtypea`\arrowtypeb`\arrowtypec;\height>}}
\def\Ctriangle{\@ifnextchar <{\Ctrianglep}{\Ctrianglep
    <\arrowtypea`\arrowtypeb`\arrowtypec;\height>}}
\def\Ctrianglep<#1>[#2`#3`#4;#5`#6`#7]{{
\settriparms[#1]
\width=\height                               
\diagram
\putCtrianglep<\arrowtypea`\arrowtypeb`
\arrowtypec;\height>
(0,0)[#2`#3`#4;#5`#6`{#7}]
\enddiagram
}}                                           
\def\putDtrianglep<#1>(#2,#3)[#4`#5`#6;#7`#8`#9]{{%
\settriparms[#1]%
\xpos=#2 \ypos=#3
\advance\xpos by\height \advance\ypos by\height
\putmorphism(\xpos,\ypos)(-1,-1)[``{#9}]{\height}{\arrowtypec}r%
\advance\xpos by-\height \advance\ypos by\height
\putmorphism(\xpos,\ypos)(1,-1)[`#5`{#8}]{\height}{\arrowtypeb}r%
{\multiply\height by 2
\putvmorphism(\xpos,\ypos)[#4`#6`{#7}]{\height}{\arrowtypea}l}%
}}

\def\putDtriangle{\@ifnextchar <{\putDtrianglep}{\putDtrianglep
    <\arrowtypea`\arrowtypeb`\arrowtypec;\height>}}
\def\Dtriangle{\@ifnextchar <{\Dtrianglep}{\Dtrianglep
   <\arrowtypea`\arrowtypeb`\arrowtypec;\height>}}
\def\Dtrianglep<#1>[#2`#3`#4;#5`#6`#7]{{
\settriparms[#1]
\width=\height                              
\diagram
\putDtrianglep<\arrowtypea`\arrowtypeb`
\arrowtypec;\height>
(0,0)[#2`#3`#4;#5`#6`{#7}]
\enddiagram
}}                                          
\def\setrecparms[#1`#2]{\width=#1 \height=#2}%

\def\recursep<#1`#2>[#3;#4`#5`#6`#7`#8]{{%
\width=#1 \height=#2
\settokens[#3]
\settowidth{\tempdimen}{$\tokena$}
\ifdim\tempdimen=0pt
  \savebox{\tempboxa}{\hbox{$\tokenb$}}%
  \savebox{\tempboxb}{\hbox{$\tokend$}}%
  \savebox{\tempboxc}{\hbox{$#6$}}%
\else
  \savebox{\tempboxa}{\hbox{$\hbox{$\tokena$}\times\hbox{$\tokenb$}$}}%
  \savebox{\tempboxb}{\hbox{$\hbox{$\tokena$}\times\hbox{$\tokend$}$}}%
  \savebox{\tempboxc}{\hbox{$\hbox{$\tokena$}\times\hbox{$#6$}$}}%
\fi
\ypos=\height
\divide\ypos by 2
\xpos=\ypos
\advance\xpos by \width
\bfig
\putCtrianglep<-1`1`1;\ypos>(0,0)[`\tokenc`;#5`#6`{#7}]%
\puthmorphism(\ypos,0)[\tokend`\usebox{\tempboxb}`{#8}]{\width}{-1}b%
\puthmorphism(\ypos,\height)[\tokenb`\usebox{\tempboxa}`{#4}]{\width}{-1}a%
\advance\ypos by \width
\putvmorphism(\ypos,\height)[``\usebox{\tempboxc}]{\height}1r%
\efig
}}

\def\recurse{\@ifnextchar <{\recursep}{\recursep<\width`\height>}}

\def\puttwohmorphisms(#1,#2)[#3`#4;#5`#6]#7#8#9{{%
\puthmorphism(#1,#2)[#3`#4`]{#7}0a
\ypos=#2
\advance\ypos by 20
\puthmorphism(#1,\ypos)[\phantom{#3}`\phantom{#4}`#5]{#7}{#8}a
\advance\ypos by -40
\puthmorphism(#1,\ypos)[\phantom{#3}`\phantom{#4}`#6]{#7}{#9}b
}}

\def\puttwovmorphisms(#1,#2)[#3`#4;#5`#6]#7#8#9{{%
\putvmorphism(#1,#2)[#3`#4`]{#7}0a
\xpos=#1
\advance\xpos by -20
\putvmorphism(\xpos,#2)[\phantom{#3}`\phantom{#4}`#5]{#7}{#8}l
\advance\xpos by 40
\putvmorphism(\xpos,#2)[\phantom{#3}`\phantom{#4}`#6]{#7}{#9}r
}}

\def\puthcoequalizer(#1)[#2`#3`#4;#5`#6`#7]#8#9{{%
\setpos(#1)%
\puttwohmorphisms(\xpos,\ypos)[#2`#3;#5`#6]{#8}11%
\advance\xpos by #8
\puthmorphism(\xpos,\ypos)[\phantom{#3}`#4`#7]{#8}1{#9}
}}

\def\putvcoequalizer(#1)[#2`#3`#4;#5`#6`#7]#8#9{{%
\setpos(#1)%
\puttwovmorphisms(\xpos,\ypos)[#2`#3;#5`#6]{#8}11%
\advance\ypos by -#8
\putvmorphism(\xpos,\ypos)[\phantom{#3}`#4`#7]{#8}1{#9}
}}

\def\putthreehmorphisms(#1)[#2`#3;#4`#5`#6]#7(#8)#9{{%
\setpos(#1) \settypes(#8)
\if a#9 %
     \vertsize{\tempcounta}{#5}%
     \vertsize{\tempcountb}{#6}%
     \ifnum \tempcounta<\tempcountb \tempcounta=\tempcountb \fi
\else
     \vertsize{\tempcounta}{#4}%
     \vertsize{\tempcountb}{#5}%
     \ifnum \tempcounta<\tempcountb \tempcounta=\tempcountb \fi
\fi
\advance \tempcounta by 60
\puthmorphism(\xpos,\ypos)[#2`#3`#5]{#7}{\arrowtypeb}{#9}
\advance\ypos by \tempcounta
\puthmorphism(\xpos,\ypos)[\phantom{#2}`\phantom{#3}`#4]{#7}{\arrowtypea}{#9}
\advance\ypos by -\tempcounta \advance\ypos by -\tempcounta
\puthmorphism(\xpos,\ypos)[\phantom{#2}`\phantom{#3}`#6]{#7}{\arrowtypec}{#9}
}}

\def\setarrowtoks[#1`#2`#3`#4`#5`#6]{%
\def\toka{#1}
\def\tokb{#2}
\def\tokc{#3}
\def\tokd{#4}
\def\toke{#5}
\def\tokf{#6}
}
\def\hex{\@ifnextchar <{\hexp}{\hexp<1000`400>}}
\def\hexp<#1`#2>[#3`#4`#5`#6`#7`#8;#9]{%
\setarrowtoks[#9]
\yext=#2 \advance \yext by #2
\xext=#1 \advance\xext by \yext
\bfig
\putCtriangle<-1`0`1;#2>(0,0)[`#5`;\tokb``\tokd]
\xext=#1 \yext=#2 \advance \yext by #2
\putsquare<1`0`0`1;\xext`\yext>(#2,0)[#3`#4`#7`#8;\toka```\tokf]
\advance \xext by #2
\putDtriangle<0`1`-1;#2>(\xext,0)[`#6`;`\tokc`\toke]
\efig
}

\begin{document}
\title{\mbox{\bf On moduli spaces of 4- or 5-instanton bundles}}
\author{}
\date{}
\maketitle
{\vskip -1.5cm {\noindent \bf F. Han}

\bigskip

\begin{abstract}
  We study the scheme of multi-jumping lines of an $n$-instanton
  bundle mainly for $n\leq 5$. We apply it to prove the irreducibility
  and smoothness of the moduli space of 5-instanton. Some particular
  situations with higher $c_2$ are also studied.
\end{abstract}
\tableofcontents
{\noindent
{\large{\bf Introduction}}
\bigskip}

An $n$-instanton is an algebraic rank 2 vector bundle $E$ over
$\proj_3(\complex)$ with Chern classes $c_1=0$, $c_2=n$, which is
stable (i.e: $h^0E=0$), and satisfies the natural cohomology property
$h^1E(-2)=0$. Denote by $I_n$ the moduli space of those bundles.

In general, little is known about $I_n$. A classical description of its
tangent space gives that every irreducible component of $I_n$ is at
least $8n-3$ dimensional. In fact it has been shown in the early
eighties that for $n\leq 4$ the moduli space
$I_n$ is irreducible (Cf [H], [E-S], [Ba2]) and smooth (Cf [LeP]).

Define a line $L$ of $\proj_3$ to be  a jumping line of order $k$ for $E$
(denoted $k$-jumping line) when the
restriction of $E$ to $L$ is $\OO_L(k) \oplus \OO_L(-k)$. Such a line
is called a multi-jumping line if $k\geq 2$. It seems 
difficult to study the moduli space without some results on the
schemes of jumping lines. Despite the fact that those schemes have a
determinantal structure involving strong properties
when $E$ is in some openset of $I_n$,  results satisfied by
any $n$-instanton are much weaker. Indeed, it seems hopeless to expect more than the
fact that any $n$-instanton has a
scheme of multi-jumping lines with good codimension in $G$ (i.e is a
curve). Surprisingly this statement which was already
conjectured in [E-S] is still open for $n>3$.

In the first 3 parts, we will study this scheme of multi-jumping lines
for any $n$-instanton with $n\leq 5$; we will prove the following:
\begin{description}
\item[Proposition \ref{c2=5courbe}] {\it The scheme of multi-jumping
lines of any $4$- or $5$-instanton is a curve in the Grassmann manifold $G$.} 
\end{description}
The techniques developed here will also apply to 
some particular families of instanton with higher $c_2$ (Cf
\ref{t'hooft} and \ref{multisautfamille}) .

As a consequence of this result, we will recover the results of Barth and
of Le Potier about $I_4$ and obtain the following one
which is also announced by G.Trautmann and A.Tikhomirov using
independent methods. 
\begin{description}
\item[Theorem \ref{irred5} and \ref{lissite}]  {\it The moduli space
 of $5$-instanton is smooth and irreducible of dimension 37.} 
\end{description}
The proofs of the irreducibility of $I_n$ for $n\leq 4$ were
all different. Our method will be to extend the one
of Ellingsrud and Stromme to \mbox{$n=4$} and $n=5$. Using the jumping
lines containing a fixed point $P$ of $\proj_3$, they constructed a
morphism from some open subset $U_P$ of $I_n$ to the moduli space $\Theta$ of
semi-table $\theta$-characteristic over plane curves of degree $n$.

The conjecture on the codimension of
multi-jumping lines  means that the $(U_P)_{(P\in \proj_3)}$
cover $I_n$. Ellingsrud and Stromme described  2 other difficulties to
extend their method:

\begin{itemize}
\item[-] Understand the images of the morphisms $U_P
\rightarrow \Theta$. This however seems only to be a real problem for
$n\geq 12$, as when $n \leq 11$ those morphism have a dense image.
\item[-] Obtain a description of the fiber of such a morphism.
\end{itemize}
In fact, the data describing this fiber are satisfying a
condition which was empty for $n=3$, and which need to be understood
for higher $n$. An explicit description of the fiber will be given in
the \S\ref{applications} for any $n$, but the key result needed to have
a grip on the moduli spaces is to understand when this fiber is
singular. This will be understood for $n=4$, and $5$, and then we will
have to check that the bundles in the ramification of all the morphism
are nevertheless smooth points of $I_n$ to obtain the above theorem.
We will need there the study of the family \ref{lisse}.

\bigskip
To prove the proposition \ref{c2=5courbe}, we will assume that there
is some $4$- or $5$-instanton $E$ such that its scheme of multi-jumping
lines contains some surface. All along the proof, we will try to
translate this excess of multi-jumping lines to some properties of
$E$ with respect to its restrictions to planes: indeed the problem of
the dimension in $\proj_3\dual$ of the non-stable planes has first been
studied by Barth (Cf [Ba1]) in 1977, and recently extended by Coanda
(Cf [Co]). Those results will be fundamental for us, because they will
be used at every end of proof, and furthermore, many proofs are
inspired by those of Coanda.

\bigskip
In the \S 1 we will study the restrictions of $E$ to planes, in other words
it gives properties of stable and semi stable rank 2 vector bundles
over the plane with $c_2\leq 5$ because the $h^1E(-2)=0$ condition
implies that $E$ has no unstable planes. For instance, we will give
informations on the intersection of the scheme of multi-jumping lines
with some $\beta$-plane, or on the bidegree of the hypothetical
congruence $S$ of multi-jumping lines.

In the \S 2 we will take care of the \prefix{k\geq 3}jumping lines, to
obtain that a general hyperplane section of $S$ can avoid those
points. This is the major difficulty for $c_2=5$, and it will enable
us to prove in the \S 3 that a general hyperplane
section of $S$ has no trisecant. The remaining cases will then be
strongly bounded, and we will first get rid of the situation where $S$ is
ruled which is easy for $c_2=5$, but not for $c_2=4$. The case $\deg
S=4$ will also be particular, and the final step will need the
residual calculus of \ref{appendiceresiduel}.

\bigskip
{\noindent
{\large{\bf Acknowledgments:}}
\bigskip}

I'd like to thanks here my thesis advisor Laurent Gruson for all his
help and his constant interest.

\bigskip
{\noindent
{\large{\bf Notations:}}
\bigskip}

In the following, we will denote by $p$ and $q$ the projections from
the points/lines incidence variety to $\proj_3$ and to the Grassmann
manifold of lines of $\proj_3$ denoted by $G$. (or to
$\proj_2$ and $\proj_2\dual$ in the plane situation). Let's define by
$\tau= p^*\OO_{\proj_3}(1)$ and $\sigma = q^*\OO_G(1)$. When a blow up
$\widetilde{\proj_3^{\prime}}$ of $\proj_3$ along a line will be used,
we will denote by $p^{\prime}$ and $q^{\prime}$ the projections to
$\proj_3$ and $\proj_1$.

Let's call by $\alpha$-plane and $\beta$-plane the 2 dimensional
planes included in $G$ which are made of the lines of $\proj_3$
containing a same point, or included in a same plane, and finally
denote by $\proj(F)=Proj(Sym F\dual)$.

Denote by $M$ the scheme of multi-jumping lines of the $n$-instanton
$E$ with $n\leq 5$. We will say that a plane $H$ is a jumping plane if
the restriction $E_H$ is semi-stable but not stable.
\begin{description}
\item[Remark]  {\it The scheme of multi-jumping lines $M$ of a 4- or
    5-instanton $E$ is at most 2 dimensional.}
\end{description}
Indeed, if this was not true, any plane would contain infinitely many
multi-jumping lines, and we will show in the \ref{infinitebis} that it
is not possible in a stable plane. Then every plane would be a jumping
plane and according to [Ba1] it is not possible when $n>1$. \fin

\bigskip
So we will finally denote by $S$ the hypothetical purely 2-dimensional
subscheme of $M$ defined by the following sequence:
\begin{center}
  $0 \longrightarrow \ideal \longrightarrow \OO_M \longrightarrow
  \OO_S \longrightarrow 0$
\end{center}
where $\ideal$ is the biggest subsheaf of $\OO_M$ which is at most
1-dimensional.

\section{Some properties of $S$ for $c_{\bf{2}}\bf{\leq 5}$}

\subsection{A link between trisecant lines to $S$ and jumping plane}

We'd like to study here the lines $d$ of $\proj_5$ which meet $M$ in
length 3 or more. Those lines have to be included in the Grassmannian
G, and each of them can be identified to a plane pencil of lines of
$\proj_3$. Let $h$ be the $\beta $-plane containing $d$, and $H$ the
associated plane of $\proj_3$. As $R^2q_{*}$ is zero, we have an  isomorphism
between the restriction of $R^1q_{*}p^{*}E$ to $h$, and $R^1q_{*}p^{*}E_H$, where $%
E_H$ is the restriction of $E$ to $H$, and also between $R^1q_{*}p^{*}E_H$
restricted to a line of $h$ and the sheaf constructed analogously
when blowing up the associated point of $H$. (So the projections over
$H$ and $\proj_1$ will be still denoted by $p$ and $q$).
\begin{proposition}
\label{pinceaux}For $c_2\leq 5$, let $d$ be a line of
$\proj_5$ meeting $M$ in a 0 dimensional scheme of length 3 or more
not containing jumping lines of order 3 or more, then the plane $H$
associated to $d$ is not a stable plane for $E$. 
\end{proposition}

Let $d$ be such a line, and assume that $H$ is stable (ie that $E_H$
is stable), and blow up $H$ at the 
point $D$ included in all the lines represented by $d$. The resolution
of this blow up is:
\begin{center}
$0\longrightarrow \OO_{d\times H}(-\tau -\sigma )\longrightarrow 
\OO_{d\times H}\longrightarrow \OO_{\widetilde{H}}\longrightarrow 0$
\end{center}
Twisting it by $p^*E_H$ gives using the functor $q_*$:
\begin{center}
$0\rightarrow q_*p^*E_H\rightarrow H^1E_H(-1)\otimes \OO_d (-1)
\longrightarrow H^1E_H\otimes \OO_d \rightarrow R^1q_*p^*E_H \rightarrow 0$
\end{center}
The length induced by the scheme of multi-jumping lines on $d$ is
$c_1(R^1q_*p^*E_H)$ and is at least 3 by hypothesis. The sheaf
$q_*p^*E_H$ is locally free of rank $-h^1E_H+h^1E_H(-1)=2$ over
$d=\proj_1$. So it has to split in $\OO_d(-a) \oplus \OO_d (-b)$ with
$a,b>0$ and $a+b=h^1E_H(-1)-c_1(R^1q_*p^*E_H)$. But $h^1E_H(-1)$ is
$c_2$, so it gives already a contradiction for $c_2\leq 4$.

For $c_2=5$, one has necessarily $a=b=1$, so $h^0(q_*p^*E_H(\sigma))$
which is also $h^0(\ideal_D \otimes E_H(1))$ is 2. Those 2 sections
$s_1$ and $s_2$ are thus proportional on the conic $Z_{s_1\wedge
  s_2}$ which has to be singular in $D$ because this point is in the
first Fitting ideal of $2\OO_H \rightarrow E_H(1)$.

If $Z_{s_1\wedge s_2}$ is made of 2 distinct lines $d_1$ and $d_2$, we
have the following exact sequence because both $s_1$ and $s_2$ vanish
at $D$:
\begin{center}
$0\longrightarrow 2\OO_H\stackrel{(s_1,s_2)}{\longrightarrow }%
E_H(1)\longrightarrow \OO_{d_1}(\alpha)\oplus \OO_{d_2}(\beta)\longrightarrow 0$
\end{center}
Computing Euler-Poincar\'e characteristics gives $\alpha+\beta=-3$, so
one of them is less or equal to $-2$, thus one of the $d_i$ is a
\prefix{k\geq 3}jumping line containing $D$ which contradicts the hypothesis.

It is the same when $Z_{s_1\wedge s_2}$ is a double line $d_1$. We
have the exact sequence:
\begin{center}
$0\longrightarrow 2\OO_H\stackrel{(s_1,s_2)}{\longrightarrow }
E_H(1)\longrightarrow \LL \longrightarrow 0$
\end{center}
where $\LL$ is a sheave of rank 1 on the plane double structure of
$d_1$ (in fact $\LL$ is not locally free at $D$). We have the linking sequence:
\begin{center}
$0\longrightarrow \OO_{d_1}(\alpha)\longrightarrow \LL \longrightarrow
\OO_{d_1}(\beta)\longrightarrow 0$ 
\end{center}
once again $\alpha+\beta=-3$, so either
$h^0\OO_{d_1}(\beta)<h^1\OO_{d_1}(\alpha)$, or
$h^1\OO_{d_1}(\beta) \neq 0$, but in both cases $R^1q_*p^*\LL_{(d_1)}$
is non zero, so $d_1$ would be a \prefix{k\geq 3}jumping line containing $D$ as
previously.

\bigskip

We can't hope such a result in a jumping plane (ie semi-stable but not
stable), but we will show in \ref{pinceauxsauteurfini} that the
problems arise only at finitely many points of the plane.

\begin{lemma}
  \label{secantes}
  Let $F$ be a bundle over $\proj_n$ with $c_1F=0$. If $F(k)$ has a
  section of vanishing locus $Z$, then for any $i\geq k-1$, the scheme
  of jumping lines of order at least $i+2$ is isomorphic to the scheme
  of lines at least \prefix{i+k+2}secant to $Z$.
\end{lemma}

We deduce from the following sequence using the standard construction:
\begin{center}
$0\longrightarrow \OO_{\proj_n} \longrightarrow F(k)\longrightarrow
\ideal_Z(2k) \longrightarrow 0$
\end{center}
an isomorphism $R^1q_*p^*F(i) \simeq R^1q_*p^*\ideal_Z(k+i)$ when
$i\geq k-1$. And the supports of those sheaves are by definition (Cf
[G-P]) the above schemes. 
\begin{lemma}\label{pinceauxsauteurfini}
For $c_2\leq 5$, any $\beta -$plane associated to a jumping plane
contain only a finite 
number of lines of $\proj_5$ meeting $M$ in a scheme of length 3 or
more made only of 2-jumping lines.
\end{lemma}
If $H$ is a jumping plane, we can take $k=0$ in lemma \ref{secantes}
and study the bisecant 
lines to $Z$. To prove the lemma it is enough to show that for any
point $N$ not in $Z$ with no \prefix{k\geq 3}jumping lines through $N$, the
line $N\dual$ is not trisecant to $M$ where $N\dual$ is the lines of
$H$ containing $N$. So let $N$ be such a point, and study the bisecant
lines to $Z$ in $N\dual$ by blowing up $H$ at $N$. One has from the
standard construction the exact sequence because $N\not\in Z$:
\begin{center}
$0\rightarrow q_*p^*\ideal_Z \rightarrow H^1\ideal_Z(-1) \otimes 
\OO_{N\dual}(-1) \longrightarrow H^1\ideal_Z \otimes \OO_{N\dual}
\rightarrow  R^1q_*p^*\ideal_Z\rightarrow 0$
\end{center}
where $q_*p^*\ideal_Z$ is isomorphic to $\OO_{\proj_1}(-a)$ for some
$a>0$ with $a=\deg Z-c_1(R^1q_*p^*\ideal_Z)$ because there are only a
finite number of bisecant to $Z$ through $N$ ($N\not\in Z$). So if $Z$
has at least 3 
bisecant (with multiplicity) through $N$, then $a\leq 2$ for $c_2\leq
5$. If $a=1$, then $Z$ is included in a line containing $N$ which must
be a $c_2$-jumping line, and it conflicts with the hypothesis. So only
the case $a=2$, and $c_2=5$ is remaining. But then, the section of
$q_*p^*\ideal_Z(2\sigma)$
would give a conic $C$ singular in $N$ and containing $Z$. We have the linking
sequence where $d_1$ may be equal to $d_2$:
\begin{center}
  $0 \longrightarrow  \OO_{d_1}(-1) \longrightarrow \OO_C
  \longrightarrow \OO_{d_2} \longrightarrow 0$ 
\end{center}
As previously $Z$ can't be included in a line so it cuts $d_2$ and
$h^0(\ideal_Z \otimes\OO_{d_2})=0$. Twisting the previous sequence by
$\ideal_Z$ gives when computing Euler-Poincar\'e's characteristics
$h^1(\ideal_Z \otimes\OO_{d_1}(-1))+h^1(\ideal_Z 
\otimes\OO_{d_2})=4$ because $Z$ is included in $C$. So one of the
$d_i$ would be a trisecant line to $Z$, thus it would be a 3-jumping
line through $N$ which contradicts the hypothesis and gives the lemma.
\begin{proposition}
 \label{infinitebis}Any plane $H$ containing infinitely many
 multi-jumping lines is a jumping plane. In that situation, the support
 of $R^1q_*p^*E$ induce on the $\beta $-plane
associated to $H$ a  reduce line with some possible
 embedded points.\\When $c_2=5$, there is a 3-jumping line in this
 pencil.\\When $c_2=4$ there are 4 embedded points.
\end{proposition}

Let's first notice that $H$ can't be a stable plane. So assume that it
is stable, then, for $c_2 \leq 5$, the bundle $E_H(1)$ has a section
vanishing on some scheme $Z$ of length at most 6. We showed in the
proof of \ref{pinceauxsauteurfini} the link between multi-jumping
lines in $H$ and trisecant lines to $Z$, but if $Z$ had infinitely many
trisecant lines (necessarily passing through a same point of $Z$),
then on the blow up of $H$ at this point, the sheaf $p^*E(\tau-3x)$
would have a section (where $x=\tau-\sigma$ is the exceptional
divisor) which is not possible because it has a negative $c_2$.

So $H$ is a jumping plane and the section of $E_H$ has infinitely many
bisecant lines (necessarily containing a same point). Let's blow up
$H$ at this point, then $p^*E_H$ has a section
vanishing 2 times on the exceptional divisor $x$. The residual scheme
of $2x$ is empty when $c_2=4$ and a single point when $c_2=5$. The
line of $H$ containing this point and the point blown up is thus a
trisecant line to the vanishing of the section of $E_H$, so it is a
3-jumping line.

Although it is not required by the following, it is interesting to
understand more this situation. For example, when $c_2=4$ we can
understand the scheme structure of multi-jumping lines in this plane.
Let $s$ be the section of $E_H$, $Z$ its vanishing locus and $\ideal_Z$
its ideal. In fact $Z$ is the complete intersection of 2 conics
singular at the same point. Let $I$ be the point/line incidence
variety in $H\dual \times H$, and pull back on $I$ the resolution of $\ideal_Z$: 
\begin{center}
$R^1q_*\OO_I(-4\tau )\longrightarrow 2R^1q_*\OO_I(-2\tau )
\longrightarrow R^1q_*p^* \ideal_Z \longrightarrow 0$
\end{center}
to obtain by relative duality:
\begin{center}
$q_*\OO_I(2\tau +\sigma )\dual \longrightarrow 2q_* \OO_I(-\sigma )
\longrightarrow R^1q_*p^* \ideal_Z \longrightarrow 0$
\end{center}
As $I= \proj(\Omega _{H\dual}(2)\dual)$ with $0\rightarrow \Omega
_{H\dual}(2) \rightarrow 3 \OO_{H\dual}(1) \rightarrow \OO_{H\dual}(2)
\rightarrow 0$, one has: \hbox{$q_* \OO_I(k\tau )=Sym_k(\Omega
  _{H\dual}(2))$}, and we are reduced to study the degeneracy locus of
a map: $2\OO_{H\dual} \longrightarrow Sym_2(\Omega _{H\dual}(2))$.
This locus contain a reduced line because through a general point $P$
of $H$ there is no conic singular in $P$ containing $Z$ so $a\geq 3$
in the proof of \ref{pinceauxsauteurfini}. So the residual scheme of
this line in the degeneracy locus has the good dimension, so we can
compute its class with the method of the appendix
\ref{appendiceresiduel}, which is simpler here because the excess is a
Cartier divisor. So this locus is given by:
\begin{center}
$c_2(Sym_2(\Omega _{H\dual}(2)))-h.c_1(Sym_2(\Omega _{H\dual}(2)))+h^2=4h^2$
\end{center}
where $h$ is the hyperplane class of $H\dual$. The geometric
interpretation of this residual locus should be that it is made of
twice the lines which occur as doubled lines in the pencil of
singular conics containing $Z$. So the ideal of the multi-jumping
lines of $E_H$ should look like $(y^2,yx^2(x-1)^2)$ in $H\dual$.

\subsection{Some preliminary results}

In the last section we studied the link between trisecant lines to $S$
and jumping planes. We can notice that some hypothetical bad
singularities of $S$ give rise to vector bundles with a line in
$\proj_3\dual$ of jumping planes. Unfortunately, such bundles exist,
for example when the instanton has a $c_2$-jumping line (Cf [S]). A
well understanding of those bundle would give many shortcuts,
particularly when $c_2=4$, but the following example shows that some
may still be unknown. 
\begin{example}
There are 4-instanton with a 2-jumping line such that every plane
containing this line is a jumping plane.
\end{example}
Let's construct this bundle from an elliptic curve $C$ and two degree 2
invertible sheaves $\LL$ and $\LL^{\prime}$ such that $(\LL\dual
\otimes \LL^{\prime})^{\otimes 4}=\OO_C$, where 4 is the smallest
integer with this property.

Let $\Sigma\stackrel{\pi }{\rightarrow }C$ be the ruled surface
$\proj(\LL \oplus \LL^{\prime})$ and $\overline{\Sigma}$ its quartic
image in $\proj_3=\proj (H^0(\LL) \oplus H^0(\LL^{\prime}))$. Let
$D=\pi^*(\LL\dual)^{\otimes 4} \otimes \OO_{\Sigma}(4)$. The linear
system $|D|$ is \mbox{$h^0((\LL\dual)^{\otimes 4}\otimes
  Sym_4(\LL\oplus \LL^{\prime}))-1$} dimensional which is 1 due to the relation
between $\LL$ and $\LL^{\prime}$. So define the bundle $E$ by the
following exact sequence:
\begin{center}
$0\longrightarrow E(-2)\longrightarrow 2\OO_{\proj_3} \longrightarrow
\phi _{*} (\pi ^{*}(\LL\dual)^{\otimes 4}\otimes {\OO}_{\Sigma}(4))\longrightarrow 0$
\end{center}
where $\phi$ is the morphism from $\Sigma$ to $\proj_3$, and where the
dimension of $|D|$ proves that $E$ is an instanton.

Let $s$ and $s^{\prime}$ be the only sections of  $\pi ^{*}\LL\dual
\otimes \OO_{\Sigma}(1)$ and of $\pi ^{*}\LL\primedual\otimes
\OO_{\Sigma}(1)$. The linear system $|D|$ is made of curves of
equation $\lambda s^4+ \mu s^{\prime4}$ with generic element a degree
8 smooth elliptic curve section of $E(2)$. Furthermore, when $\lambda $
or $\mu$ vanishes, this element of $|D|$ provides to the double line
of $\overline{\Sigma}$ associated to $s$ or $s^{\prime}$ a multiple
structure. Every plane containing one of those lines (notes $d$ and
$d^{\prime}$) must be a jumping plane, and every ruling of
$\overline{\Sigma}$ is a 4-secant to those section of $E(2)$, so they
are 2-jumping lines. In any plane containing $d$ or $d^{\prime}$, the
jumping section vanishes on a scheme of ideal $(y^2,x(x-1))$ where $y$
is the equation of $d$ or $d^{\prime}$. So the lines $d$ or
$d^{\prime}$ must be 2-jumping lines, and those sections of $E(2)$
induce on the exceptional divisor of $\proj_3$ blown up in $d$ or
$d^{\prime}$ a curve of bidegree $(2,2)$.
\begin{remark}
\label{dtemultiple}For every $c_2$, if there is a line $\delta $,at
least 1-jumping for some instanton $E$ such that $\ideal _\delta
^k\otimes E(k)$ has a section, then $\delta $ is a $k$-jumping line,
and for $k \geq 2$ this section is irreducible.
\end{remark}

Let's first notice that if $\ideal _\delta^k\otimes E(k)$ has a
section when $\delta$ is a jumping line, then this section is
set-theoretically $\delta$. Indeed, consider the blowing up
$\widetilde{\proj_3^{\prime}}$ of $\proj_3$ along $\delta$, and denote
by $p^{\prime}$ et $q^{\prime}$ the projections of $\widetilde{\proj_3^{\prime}}$ on
$\proj_3$ et $\proj_1$. The restriction of the section of $p^{\prime
  *}E(k\sigma )$ to the exceptional divisor is a section of
$\OO_{\proj_1 \times \proj_1}(k,a) \oplus \OO_{\proj_1 \times
  \proj_1}(k,-a)$ where $\delta$ is a $a$-jumping line. As $a>0$, this
section is in fact a section of $\OO_{\proj_1 \times \proj_1}(k,a)$,
so it has no embedded nor isolated points. So the section of $p^{\prime
  *}E(k\sigma )$ don't have irreducible components which meet the
exceptional divisor, but the section of $\ideal^k_{\delta} \otimes
E(k)$ has to be connected in $\proj_3$ for $k\geq 2$ because
$h^1E(-2)=0$. So this section must be set-theoretically $\delta$ when
$k\geq 2$.

On another hand, the class in the Chow ring of
$\widetilde{\proj_3^{\prime}}$ of a divisor of the exceptional $\proj_1 \times
\proj_1$ of bidegree $(k,a)$ is: $k\tau ^2+(a-k)\tau\sigma $, but this
curve is included in a section of $p^{\prime *}E(k\sigma )$ of class
$c_2\tau^2$ which gives $a=k$.

\bigskip
The main problem encountered to classify those bundle comes from the
fact that the section of $p^{\prime *}E(k\sigma )$ may have a multiple
structure not included in the exceptional divisor as in the
previous example. So we will not try to classify all those bundle, but
we will need the following:
\begin{lemma}
\label{dteplansaut}When $c_2=4$ and $h^0E(1)=0$, if there is a
\prefix{k\leq 2}jumping line $d$, such that every plane containing $d$
is  a jumping one, then
$\ideal_d ^{\otimes 2}\otimes E(2)$ has a section
\end{lemma}

Still blow up $\proj_3$ along $d$, the class of the exceptional
divisor $x=\proj_1 \times d$ is $\tau - \sigma$, and we have the exact
sequence:
\begin{center}
$0\longrightarrow \OO_{\proj_3\times\proj_1}(-\tau -\sigma
)\longrightarrow \OO_{\proj_3\times \proj_1}\longrightarrow 
\OO_{\widetilde{\proj_3^{\prime}}}\longrightarrow 0$
\end{center}
which gives when twisted by $p^{\prime *}E(\tau )$ using the functor
$q_*^{\prime }$ the following exact sequence because $h^0E(1)=0$.
\begin{center}
$0\rightarrow \!q_{*}^{\prime }p^{\prime *}E(\tau )\rightarrow \!H^1E\otimes 
\OO_{\proj\!_1}(-1) \kercoker{K} \!\!H^1E(1)\otimes \OO_{\proj\!_1}\rightarrow
\!R^1q_{*}^{\prime }p^{\prime *}E(\tau )\rightarrow \!0$
\end{center}
The sheaf $R^1q_{*}^{\prime }p^{\prime *}E(\tau )$ is locally free of
rank 1 over $\proj_1$ because every plane containing $d$ is a jumping
one and $c_2=4$. On another hand, the resolution of $x$ in
$\widetilde{\proj_3^{\prime}}$ gives the exact sequence:
\begin{center}
$0\rightarrow p^{\prime *}E(\sigma )\rightarrow p^{\prime *}E(\tau
)\rightarrow p^{\prime *}E_x(\tau )\rightarrow 0$
\end{center}
But $p^{\prime *}E_x(\tau )\simeq \OO_{\proj\!_1\times
  d}(0,-a+1)\oplus \OO_{\proj\!_1\times d}(0,a+1)$ with $a\leq 2$, and
  the map \linebreak\mbox{$R^1q_{*}^{\prime }p^{\prime *}E(\sigma
)\rightarrow R^1q_{*}^{\prime }p^{\prime *}E(\tau )$} is a surjection,
  so $R^1q_{*}^{\prime }p^{\prime *}E(\tau )\simeq \OO_{\proj_1}(b)$
  with \mbox{$b\geq 1$} due to the surjection of $H^1E\otimes \OO_{\proj_1}$
  onto $R^1q_{*}^{\prime }p^{\prime *}E$. But $h^1(K)=0$ so
  \mbox{$h^1(q_{*}^{\prime }p^{\prime *}E(\tau ))\leq 2$} and as
  $q_{*}^{\prime }p^{\prime *}E(\tau )$ is locally free of rank 3, we
  can deduce that $h^0(q_{*}^{\prime }p^{\prime *}E(\tau +\sigma
  ))\neq 0$. So we have a section of $\ideal_d\otimes E(2)$ whose
  restriction to every plane containing $d$ is proportional to the
  jumping section because the jumping section don't have a vanishing
  locus included in a line ( $d$ is not $c_2$-jumping). So in every
  plane containing $d$, the section of $\ideal_d E(2)$ must vanish on
  a conic which has to be twice $d$, so it gives a section of
  $\ideal_d^{\otimes 2}E(2)$.

\begin{example}
 \label{alphaun}There are bundles  $E$ with $h^1E(-2)\mbox{mod }2=1$,
$c_2=4$ having a line congruence of multi-jumping lines.   
\end{example}
Take 2 skew lines and put on each of them a quadruple structure made
by the complete intersection of 2 quadrics singular along those lines.
So the disjoint union of these two elliptic quartics gives a section
of $E(2)$ where $E$ has $c_2=4$ and $h^1E(-2)\mbox{mod }2=1$ according
to [C]. Every line meeting those 2 lines is 4-secant to this section
which gives a $(1,1)$ congruence of 2-jumping lines.

\begin{lemma}
\label{t'hooft}For every $c_2$, the t'Hooft bundles (ie $h^0E(1)\neq
0$) have a 1-dimensional scheme of multi-jumping lines.
\end{lemma}

We can remark that using deformation theory around such a bundle Brun
and Hirshowitz (Cf [B-H]) proved that any instanton in some open
subset of the irreducible component of $I_n$ containing the t'Hooft
bundles has a smooth curve as scheme of multi-jumping lines. This is
much stronger but unfortunately the t'Hooft bundles are not in this
open set because they always have a 3-jumping line when $n \geq 3$. ( Take a
quadrisecant to the section of $E(1)$ which is made of $n+1$ disjoint
lines Cf [H])

\smallskip
So, let $s$ be a section of $E(1)$ with vanishing locus $Z$. For any
$k\in \nat$, according to the \ref{secantes}, the scheme of
\prefix{k+2}jumping lines is isomorphic to 
the scheme of \prefix{k+3}secant lines to $Z$. The lemma is thus
immediate when the lines making $Z$ are reduced. If it is not the
case, a 2-parameter family of trisecant to $Z$ may arise from 2 kinds
of situations:

a) There is a line $d$ such that the scheme induce on $d$ by $Z$ has a
congruence of bisecant lines which also meet another line $d^{\prime}$
of $Z$. This congruence must then be set-theoretically the lines
meeting $d$ and $d^{\prime}$, so for any plane $H$ containing
$d^{\prime}$, the lines passing through $d\cap H$ are bisecant to $Z$
at this point. The line $d$ would then be equipped by $Z$ of a multiple
structure doubled in every plane containing $d$, which is impossible
because those planes would be unstable.

b) The scheme induce by $Z$ on some line $d$ (noted $Z_d$) has a 2
-dimensional family of trisecant lines which can't lye in a same
plane. Every plane containing $d$ is a jumping one, because the
restriction of the section of $E(1)$ to this plane is vanishing on
$d$, and has to contain infinitely many multi-jumping lines by
hypothesis. As those lines are bisecant to the jumping section which
is 0 dimensional, those lines must form a pencil through a point of
$d$ by hypothesis. So $d$ is a multi-jumping line which contradicts the
\ref{dtemultiple} because $s$ is a section of $\ideal_d \otimes E(1)$.

\bigskip
We can now reformulate using the lemma \ref{t'hooft} the
theorem of Coanda (Cf [Co]) in a form which will be used
at many times in the following: 
\begin{theorem}
  \label{coanda}
  {\rm \bf (Coanda)} For any $c_2$, if an instanton has not a 1-dimensional family of
  multi-jumping lines, then it has at most a 1-dimensional family of
  jumping planes.
\end{theorem}
Indeed, Coanda showed that any bundle with a 2-dimensional family of
jumping planes is either a special t'Hooft bundle which has its
multi-jumping lines in good dimension according to the \ref{t'hooft}, or
another kind of bundles which are not an instanton. In fact, the last
one are in the ''big family'' of Barth-Hulek (Cf [Ba-Hu]).

\subsection{A bound on the bidegree of $S$}

Let $(\alpha, \beta)$ be the bidegree of $S$, where $\alpha$ (resp
$\beta$) is the number of lines of the congruence $S$ passing through
a general point of $\proj_3$ (ie in a $\alpha$-plane), (resp in a
general plane of $\proj_3$ (ie in a $\beta$-plane)).
\begin{proposition}
\label{basdeg}For any $c_2$, if an instanton has at most a
2-dimensional scheme of multi-jumping lines , and if  $h^0(\ideal
_P^{\otimes c_2-4} \otimes E(c_2-4))$ is zero for a general point of
$\proj_3$, then 
$\alpha \leq 2c_2-6$.\\ Similarly, any stable plane $H$ without
\prefix{c_2} and \prefix{c_2-1}jumping lines and such that $h^0(\ideal
_P^{\otimes k}\otimes E_H(k))$ is zero when $P$ is general in $H$ for
$0\leq k\leq c_2-4$, then the $\beta$-plane $h$ associated to $H$ cuts the
scheme of multi-jumping lines in length at most $2c_2-6$. So if those
hypothesis are satified for the general plane then $\beta \leq 2c_2-6$. 
\end{proposition}
{\bf Remark:} The above hypothesis about $h^0(\ideal
_P^{\otimes k}\otimes E_H(k))$ are always satisfied when $c_2=4$ or $5$
because one has in any stable plane $H$ without \prefix{c_2} and
\prefix{c_2-1}jumping lines: $h^0(E_H(1))\leq 2$. Furthermore, we will
show for those $c_2$, under the assumption of the existence of $S$, that
the general member of any 2-dimensional 
family of planes is stable from the theorem of Coanda stated in the
\ref{coanda}, and don't contain \prefix{c_2} or \prefix{c_2-1}jumping
lines according to the \ref{maxsauteuses}, \ref{fini4}, and \ref{courbetris}.

\bigskip
Start first with the bound of $\alpha$, so take an $\alpha$-plane $p$
cutting $S$ in length $\alpha$ and denote by $P$ its associated point
of $\proj_3$. The standard construction associated to the blow up of
$\proj_3$ at $P$ gives the following exact sequence above the
exceptional divisor $p$:
\begin{center}
$0\longrightarrow q_*p^*E \longrightarrow H^1E(-1)\otimes (\OO_p
\oplus \OO_p(-1)) \longrightarrow H^1E\otimes \OO_p \longrightarrow
R^1q_*p^*E \longrightarrow 0$
\end{center}
The sheaf $q_*p^*E$ is thus a second local syzygy, so it has a
3-codimensional singular locus hence it is a vector bundle denoted by
$F$ in the following. For every instanton, one has $h^1E(-1)=c_2$ and
$h^1E=2c_2-2$. So $\alpha$ is given by
$\chi(R^1q_*p^*E)=c_2-2+\chi(F)$.

The stability of $E$ implies that $F$ has no sections, but it is also
locally free of first Chern class $-c_2$,
so $\chi(F)=-h^1F+h^0F(c_2-3)$. Take a $\proj_1$ in $p$ without any
multi-jumping lines. The restriction of this sequence to this $\proj_1$
gives an injection of $F_{\proj_1}$ in $c_2\OO_{\proj_1} \oplus
c_2\OO_{\proj_1}(-1)$, so the bundle $F_{\proj_1}$ has to split into
$\OO_{\proj_1}(-a) \oplus \OO_{\proj_1}(-b)$ with $a+b=c_2$ because
this line avoids the multi-jumping lines by hypothesis. We can bound
$h^0F_{\proj_1}(c_2-3)$ by $c_2-4$ except when $a=0$ or $1$ (which may
happen, for example if $E$ had infinitely many jumping planes). But in
the cases $a=0$ or $1$, one has $h^0F_{\proj_1}(c_2-3)=c_2-2-a$, and
as we can assume that $P$ is not on the vanishing of a section of
$E(1)$ because for $c_2\geq 2$ $h^0E(1)\leq 2$, we have $h^0F(1)=0$,
so there is an injection of $H^0F_{\proj_1}(1)$ into $H^1F$. Then in
the cases $a=0$ or $1$, one has $h^1F \geq 2-a$.
But $h^0F(c_2-4)$ is zero by hypothesis so one has \linebreak$h^0F(c_2-3)
\leq h^0F_{\proj_1}(c_2-3)$ thus $\chi(F)$ is also bounded by $c_2-4$
in the cases $a=0$ or $1$, which gives $\alpha \leq 2c_2-6$

\bigskip
Let's now take care of the bound of $\beta$. Take a stable plane $H$
without \prefix{c_2} and \prefix{c_2-1}jumping lines, and denote by $\beta^{\prime}$
the length of the intersection of the $\beta$-plane $h$ with the
scheme of multi-jumping lines $M$. Consider now the incidence 
variety $I\subset H\dual \times 
H$. The resolution of $I$ twisted by $p^*E_H$ gives the following
exact sequence:
\begin{center}
$0\longrightarrow q_*p^*E_H \longrightarrow H^1E_H(-1)\otimes
\OO_{H\dual}(-1) \longrightarrow H^1E_H\otimes
\OO_{H\dual}\longrightarrow R^1q_*p^*E_H \longrightarrow 0$  
\end{center}
One has $h^1E_H(-1)=c_2$ and $h^1E_H=c_2-2$, so
\mbox{$\beta^{\prime}=\chi(R^1q_*p^*E_H) = c_2-2+\chi (F)$}, where this time $F$ is
$q_*p^*E_H$ which is still locally free of rank 2 and $c_1F=-c_2$, so
we have $\chi (F)=-h^1F+h^0F(c_2-3)$. As $H$ is stable without
\prefix{c_2} and \prefix{c_2-1}jumping lines the vanishing locus of a
section of $E_H(1)$ is at most on one conic from the \ref{secantes},
so we have $h^0E_H(1) \leq 2$, and we can take a point $P$ of $H$
which is neither in a multi-jumping line nor in the vanishing locus of
some section of $E_H(1)$. The line $p\subset H\dual$ associated to $P$
don't contain any multi-jumping lines, so we have an injection: $0\rightarrow
F_p \longrightarrow c_2\OO_p(-1)$ and $F_p=\OO_{\proj_1}(-a)\oplus
\OO_{\proj_1}(-b)$ with $a+b=c_2$. But $P$ is not in the vanishing
locus of a section of $E_H(1)$ so $a$ and $b$ are at least 2 then
$h^0F_p(c_2-3) \leq c_2-4$. By hypothesis, $h^0F_p(k)=0$ for $0 \leq k
\leq c_2-4$, so for those $k$ we have by induction that all the
$h^0F(k)$ are zero using the sequence of restriction of $F(k)$ to $p$.
In particular we have $h^0F(c_2-4)=0$ so $h^0F(c_2-3)$ is bounded by
$h^0F_p(c_2-3) \leq c_2-4$, which gives $\beta^{\prime} \leq 2c_2-6$.
Furthermore, the general plane is stable according to [Ba1], so if it
satisfies the conditions on $h^0(\ideal_P^{\otimes k}\otimes E_H(k))$
then we have $\beta=\beta^{\prime}\leq 2c_2-6$. \fin

\section{``No'' k-jumping lines on $S$ when $k\geq 3$}\label{fini3sautet+}

The aim of this section is to prove that $S$ is made of 2-jumping
lines except at a finite number of points. Let's first take care of
some extremal cases:

\subsection{Some too particular jumping phenomena}
We want to prove here that under the hypothesis of the existence of
$S$ with $c_2\leq 5$, there is at most a \prefix{c_2-1-k}dimensional family of k-jumping
lines for $3 \leq k \leq 5$.
\begin{lemma}
 \label{trisautplan}For $c_2\leq 5$ and $i=0$ or $1$, the support of $R^1q_*p^*E(1+i)$
 meets any $\beta $-plane $h$ which doesn't have \prefix{k\geq 4+i}jumping
 lines in a scheme of length at most $c_2-3-i$.
\end{lemma}
Let's first get rid of the case of a stable plane $H$. For $c_2\leq
5$, we can pick a section of $E_H(1)$, and note $Z$ its vanishing
locus. One has $h^1E_H(1)=h^1\ideal_Z(2)$ which is at most 1 because
$\deg Z\leq 6$ and $Z$ has no \prefix{5+i}secant because there is no
\prefix{4+i}jumping lines in $H$. So $R^1q_*p^*E_H(1)$ is the cokernel of some map
$(c_2-2)\OO_h(-1) \rightarrow \OO_h$, and its support is 0-dimensional
according to the \ref{infinitebis}, so it has length 0 or 1.
Similarly, if there is a 4-jumping line in $H$, then
$R^1q_*p^*E_H(2)$ is the cokernel of some map $2\OO_h(-1) \rightarrow
\OO_h$ , so the support of this sheaf has length 1.

If $H$ is a jumping plane, then denote by $Z$ the vanishing locus of
the section of $E_H$. So we have to study the \prefix{3+i}secant lines
to $Z$ with $\deg Z =c_2$, so it is the same problem as the study of the
\prefix{2+i}jumping lines of a stable bundle over $H$ with $c_2\leq
4$. So the result for $i=0$ is a deduced from \ref{basdeg}, and to
obtain $i=1$, remark that we have shown above that there is at most one
3-jumping line in a 
stable plane when $c_2 \leq 4$, so there is at most one 4-jumping line
in any plane when $c_2 \leq 5$. \fin

\bigskip
So we can obtain the following:
\begin{lemma}
The jumping lines of order at least 3 of a 4- or 5-instanton make an
at most 1-dimensional scheme.
\end{lemma}

Indeed, if $E$ has a line congruence of \prefix{k\geq 3}jumping lines
of bidegree $(\alpha,\beta)$, then $\alpha \leq 1$ and $\beta \leq 1$
because we have shown above that when $c_2\leq 5$ there was no two
3-jumping lines in a stable plane, and there is at most a one
parameter family of jumping planes from Coanda's theorem stated in
\ref{coanda}. But the congruence $(0,1)$, $(1,0)$, $(1,1)$ contain lines (Cf [R]), 
and for those $c_2$ it is impossible to have a plane pencil of 3-jumping lines
according to the \ref{infinitebis}.
\begin{remark}
\label{maxsauteuses}When $c_2\leq 5$, if $E$ has a $c_2$-jumping line
$d$,then its multi-jumping lines are 3-codimensional in $G$.
\end{remark}
Let's first recall that the instanton property imply directly that any
plane is semi-stable. Furthermore, any semi-stable plane containing a
$c_2$-jumping line doesn't contain another multi-jumping line because
the vanishing locus $Z$ of a section of $E_H(1)$ or of $E_H$ is
included in $d$ because $d$ is \prefix{\deg Z}secant to $Z$, and the
multi-jumping lines are at least 2-secant to $Z$.

But the lines of
$\proj_3$ meeting $d$ is a hypersurface of $G$, so it would cut $S$ in
infinitely many points and we could find a multi-jumping line in the
same plane as $d$ which is impossible.

\bigskip
So we will assume in the following that $E$ has no \prefix{c_2}jumping lines.
\begin{proposition}
\label{fini4}The assumption of the existence of $S$ implies that $E$
has at most a finite number of 4-jumping lines when $c_2=5$.
\end{proposition}

Assume here that there is a 4-jumping line $q$. Consider the
hypersurface of $S$ made of the lines of the congruence meeting $q$
(i.e: $T_qG\cap S$). But, when
$c_2=5$, any plane $H$ containing $q$ and another multi-jumping line can't
be stable from the \ref{secantes}, because  if $H$ was stable, the
vanishing locus of a section of $E_H(1)$ would 
be a scheme of degree 6 with $q$ as 5-secant line, so it 
couldn't have another trisecant. We had also showed in the
\ref{infinitebis} that any plane containing $q$ has necessarily a
finite number of multi-jumping lines, so the curve  $T_qG\cap S$
can't be a plane curve, thus any plane containing $q$ must be a
jumping plane, and if there was infinitely many such lines $q$, there
would exist a 2-dimensional family of jumping planes, and we could
conclude with Coanda's theorem as in the \ref{coanda}.

\subsection{Curves of 3-jumping lines}
One wants here to prove that there is no curve of \prefix{k\geq
  3}jumping lines lying on $S$. When $c_2=4$ this is exactly the
  method of Coanda, but there are many more problems when $c_2=5$
  because the ruled surface used may have a singular locus. Let's
  first consider the case of a curve of 3-jumping lines.
\begin{proposition}
\label{courbetris}No 4-instanton with $h^0E(1)=0$ has a curve
  of jumping lines of order exactly 3.\\When $c_2=5$, if $E$ has a surface $S$ of
  multi-jumping lines, then there is no curve of 3-jumping lines on $S$.
\end{proposition}
Assume that there is such a curve of 3-jumping lines (ie without
4-jumping lines), and take
$\Gamma$ be an integral curve made only of 3-jumping lines. Denote by
$\Sigma$ the ruled surface associated to $\Gamma$, and
$\widetilde{\Sigma } \stackrel{\pi }{\longrightarrow }\Sigma $ its
smooth model.

\begin{itemize}
\item When $c_2=4$
\end{itemize}
We showed in the \ref{trisautplan} that two 3-jumping
lines can't cut one another, so according to Coanda's result (Cf
[Co]), the curve $\Gamma$ is either a regulus of some quadric $Q$ of
$\proj_3$, or the tangent lines to a skew cubic curve, so in both cases,
one has  $\widetilde{\Sigma }\simeq \proj_1 \times \proj_1$. Let's now
apply the method used by Coanda in [Co]. As there is no 4-jumping
lines in $\Gamma$, we have the exact sequence for some divisor $A$ of
degree $a$ in Num($\Sigma $):
\begin{center}
$0 \longrightarrow \OO_{\widetilde{\Sigma }}(-A,3)\longrightarrow \pi
^*E_\Sigma \longrightarrow \OO_{\widetilde{\Sigma }}(A,-3) \longrightarrow 0$
\end{center}
which gives $c_2(\pi ^*E_\Sigma )=d.c_2=2a3$, where $d=\deg \Sigma$ is
2 or 4. So it gives a contradiction because here $c_2=4$.

\begin{itemize}
\item When $c_2=5$
\end{itemize}
We also have such an exact sequence, but this time
$\widetilde{\Sigma }= \proj_{\widetilde{\Gamma }}(F\dual)$ where
$\widetilde{\Gamma }$ is the normalization of $\Gamma$, and $F$ is a
rank 2 vector bundle having a section (denote by $\LL$ its cokernel),
and such that $h^0F(l)=0$ for any invertible sheaf $l$of negative
degree. The surjection from $F$ to $\LL$ gives a section $C_0$ of
$\widetilde{\Gamma }$ into  $\widetilde{\Sigma }$. Let $e=-\deg F$,
the intersection form in \mbox{Num$\widetilde{\Sigma }$} is given by
 $\left( 
\begin{array}{cc}
0 & 1 \\ 
1 & -e
\end{array}
\right) $ in the basis $(f,C_0)$ where $f$ is the class of a fiber, so
we have the relation $5d=6a+9e$ which proves that $d$ is a multiple of
3. On another hand, the dual ruled surface $\Sigma\dual$ has no triple
points because there is at most two 3-jumping lines in a same plane
from the \ref{trisautplan}, so we have the relation from [K]:
\begin{equation}
  \label{nopointtriple}
  [(d-2)(d-3)-6g](d-4)=0
\end{equation}
because $\widetilde{\Sigma }$ and $\widetilde{\Sigma\dual }$ have same
degree and genus. So when $d=3$ we can work as in [Co] because it
gives $g=0$, and $\widetilde{\Sigma }=\proj(\OO_{\proj_1}\oplus
\OO_{\proj_1}(1))$, so $e=1$ and $a=1$, then we have:\\
\centerline{$\pi ^* \OO_\Sigma
(1) \stackunder{Num}{\sim } \OO_{\widetilde{\Sigma }}(2,1)$, and
$h^2(\pi^* E_\Sigma
(-2))=h^0(\OO_{\proj_1}(-a+4+2g-2-e)\otimes Sym_3F)$}\\ which
imply $h^2E_{\Sigma}(-2)\geq 1$ when taking sections in the following
sequence where $C^{\prime}$ is at most 1-dimensional. 
\begin{center}
$0 \longrightarrow E_\Sigma (-2)\longrightarrow \pi _*\pi ^*E_\Sigma
(-2)\longrightarrow E_\Sigma (-2)\otimes \omega _{C^{\prime}} \longrightarrow 0$
\end{center}
The resolution of $\Sigma$ in $\proj_3$ twisted by $E(-2)$ gives then
$h^3E(-5)\geq 1$, so $h^0E(1)\geq 1$ which contradicts the hypothesis.
Thus $d$ is a multiple of 3 and $d\geq 6$. We will now need the
following 3 lemmas before going on proving the proposition \ref{courbetris}

\begin{lemma}
  A generic ruling of $\Sigma$ can't meet points of multiplicity 3 or
  more of $\Sigma$
\end{lemma}

If the opposite was true, there would be a point $P$ of $\Sigma$
included in 3 distinct ruling of $\Sigma$. Those 3 lines would give 3
triples points of the pentic curve representing the jumping lines
through $P$. This pentic must then be reducible, and those triple
points have to lie in a same line which contradicts the fact that a
plane has at most two 3-jumping lines (Cf \ref{trisautplan}).
\begin{lemma}
\label{lieudoublereduit}Any generic ruling of $\Sigma $ meets only
reduced components of the double locus of $\Sigma $
\end{lemma}

Assume that the double locus of $\Sigma$ has a reducible component
meeting all the rulings. So any point of this component is a point $P$
of the double locus where the tangent cone $C_P\Sigma$ is a double
plane, so this plane must contain the 2 rulings passing through $P$.
On another hand, there are pinch points on $\Sigma$ because the number
$2d+4(d-1)$ (Cf [K]) is not zero when $d$ is multiple of 3 and solution of
\ref{nopointtriple}. So let $P^{\prime}$ be a pinch point of $\Sigma$,
then there is through $P^{\prime}$ a ruling $d$ of $\Sigma$ such that
the tangent plane to $\Sigma$ at any point of $d$ is a same plane $H$.
This means that the tangent space to $\Gamma$ at $d$ is included in
the $\beta$-plane $h$ associated to $H$. So $\Gamma\cap h_{\{d\}}$ has
length at least 2. But under the above hypothesis we can find another
ruling of $\Sigma$ in $H$. Indeed, by assumption $d$ meets a non
reduced component of the double locus in some point $P$, so there is
another ruling $d^{\prime}$ of $\Sigma$ through $P$ and we showed
that $C_P\Sigma$ was set-theoretically the plane containing $d$ and
$d^{\prime}$, but as the tangent space is constant along $d$, $H$ has
to be included in $C_P\Sigma$ so $d^{\prime}$ is also in $H$ and
$h\cap \Gamma$ has length at least 3 which contradicts the \ref{trisautplan}. 
\begin{lemma}
\label{nplanssaut}A generic ruling of $\Sigma $ meets at least 4 other
distinct rulings of $\Sigma$.
\end{lemma}

a) One have first to proves that a generic ruling of $\Sigma$ meets the
double locus in  distinct points. According to the
\ref{lieudoublereduit}, the only case where this could be not
satisfied is when the generic ruling of $\Sigma$ is
tangent to some component of the double locus. That is
the case where $\Sigma$ is a developable surface. The cone situation
is immediate for the lemma \ref{nplanssaut}, so we will first consider
the case where $\Sigma$ is a the developable surface of tangent lines
to some curve $C$.

For any point $P$ of $C$, the osculating plane $H_P$ to $C$ must
coincide with the tangent plane to $\Sigma$ at any point of the ruling
$d_P$ which is tangent to $C$ at $P$. Furthermore, the tangent space
to $\Gamma$ at $d_P$ is made of the pencil of lines of $H_P$
containing $P$, so it is included in the $\beta$-plane associated to
$H_P$.

We can also remark that for any point $P^{\prime}$ of $C$ distinct
from $P$, the planes $H_P$ and $H_{P^{\prime}}$ are distinct,
otherwise we the $\beta$-plane associated would cut $\Gamma$ in length
4 which contradicts the \ref{trisautplan}. So we can translate this by
the fact that the dual surface $\Sigma\dual$ is a developable surface
of tangent lines to a smooth curve $C\dual$.

Consider the projection of $C\dual$ to $\proj_2$ from a point of
$\proj_3\dual$ not in $\Sigma\dual$. The degree of $\Gamma$ is also
the degree of the dual curve of this projection, which is according to
Pl\"ucker's formula $d=2(\deg C\dual -1)+2g$ because this projection
don't have cusp. So when $\Sigma$ is developable, the case $d=6, g=2$
is not possible , then  $d\geq 9$. Let $t$ be a general point of
$\Gamma$, then the contact between $T_tG$ and $\Gamma$ can't be of
order 6 or more, otherwise the osculating plane to $\Gamma$ at $t$
would be included in $G$, so there would exist a $\beta$-plane meeting
$G$ in length 3 or more which contradicts the \ref{trisautplan}. But
the contact between $T_tG$ and $\Gamma$ must be even (Cf [Co] lemma
6), so $T_tG\cap \Gamma-\{t\}$ is made of $d-4$ distinct points
because in general $t$ is not bitangent to $C$ because $\Gamma$ is
reduced. So, when $\Sigma$ is developable a generic ruling of $\Sigma$
meets at least 5 other distinct ones.

b) So we can assume that $\Sigma$ is not developable, and using the
\ref{lieudoublereduit}, that a generic ruling of $\Sigma$ meets the
double locus in distinct points. One have thus to compute the number
of such points. But as $\Sigma$ is not developable, in general a point
$t$ of $\Gamma$ is such that the intersection of $T_tG$ with
$\Gamma-\{t\}$ has length $d-2$, which gives the lemma because $d\geq
6$. \fin

\bigskip
We are now ready to continue the proof of \ref{courbetris}. We
obtained in the \ref{nplanssaut} that through any ruling $t$ of
$\Sigma$ there are at least 4 distinct planes containing another
ruling of $\Sigma$. So according to the \ref{trisautplan}, we have at
least 4 jumping planes through $t$. We can assume that there are only
a finite number of jumping planes through $t$ when $t$ is general
because the bundle $E$ don't have a 2-dimensional family of jumping
planes (\ref{coanda}). Let's blow up $\proj_3$ along $t$ and denote by
$p^{\prime}$ and $q^{\prime}$ the projections over $\proj_3$ and
$\proj_1$. We have from the standard construction the following sequence:
\begin{center}
$0\rightarrow \!q_{*}^{\prime }p^{\prime *}E(\tau )\rightarrow \!H^1E\otimes 
\OO_{\proj\!_1}(-1) \kercoker{K} H^1E(1)\otimes \OO_{\proj\!_1}\rightarrow
\!R^1q_{*}^{\prime }p^{\prime *}E(\tau )\rightarrow \!0$
\end{center}
where $h^1E=8,h^1E(1)=7$, and where $R^1q_*^{\prime }p^{\prime
  *}E(\tau)$ is not locally free at the jumping planes but has rank 1
  because in a stable plane $H$ containing the 3-jumping line $t$, the 6
  points of the vanishing locus of a section of $E_H(1)$ must lie on a
  conic because $t$ is 4-secant to this locus. The sheaf
  $q_*^{\prime }p^{\prime *}E(\tau)$ is thus locally free of rank 2,
  and we have: $h^0(q_*^{\prime }p^{\prime *}E(\tau))=0$, $h^1(K)=0$,
  $h^1(q_*^{\prime }p^{\prime *}E(\tau ))=h^0(K)$, et
  $h^0(R^1q_*^{\prime } p^{\prime *} E(\tau ))=7-h^1(q_*^{\prime }
  p^{\prime *}E(\tau ))$.

  If $h^1(q_*^{\prime }p^{\prime *}E(\tau ))$ is 0 or 1 then
  $q_*^{\prime } p^{\prime *}E(\tau )= \OO_{\proj_1}(-1) \oplus
  \OO_{\proj_1} (-1\ or\ -2)$, so $[q_*^{\prime }p^{\prime *}E(\tau
  )](1)$ would have a section and $\ideal_t\otimes E(2)$ also. But
  this phenomena must occurred for a generic $t$, so we have $h^0E(2)
  \geq 2$, but the vanishing locus of those sections of $E(2)$ would
  then be in a quartic surface, which must contain every 3-jumping
  line because those lines are 5-secant to those vanishing locus. But
  it contradicts $\deg \Sigma \geq 6$.

  So we have $h^0(R^1q_*^{\prime }p^{\prime *}E(\tau )\leq 5$. On
  another hand $R^2q_*^{\prime }p^{\prime *}E(\tau )$ is zero from
  Grauert's theorem, so we have by base change for all $y$ in $\proj_1$:\\
$R^1q_*^{\prime }p^{\prime *}E(\tau )_{\{y\}}\simeq H^1E_Y(1)$ where
  $Y$ is the plane associated to $y$. So we have the exact sequence:
\begin{center}
$0\longrightarrow \stackunder{k\ times}{\oplus }\complex
\longrightarrow R^1q_*^{\prime } p^{\prime *}E(\tau )\longrightarrow
R^1q_*^{\prime }p^{\prime *}E(\tau )\bidual \longrightarrow 0$
\end{center}
where $k$ is the number of jumping planes through $t$, so $k\geq 4$
from the \ref{nplanssaut}.  But $R^1q_*^{\prime }p^{\prime *}E(\tau
)\bidual \simeq \OO_{\proj_1}(l)$ with $l\geq 0$, so necessarily
$l=0,k=4$. Then $d=6$ so $g=2$ from the formula \ref{nopointtriple}.
We also have $h^0(R^1q_*^{\prime }p^{\prime *}E(1))=5$ and
$q_*^{\prime} p^{\prime *}E(1) \simeq \OO_{\proj_1}(-2) \oplus
\OO_{\proj_1}(-2)$, so $h^0\ideal_t^2\otimes E(3) \geq 2$.

Take $s$ and $s^{\prime}$ 2 sections of $\ideal_t^2\otimes E(3)$.
Those sections are proportional on a sextic surface containing the 4
jumping planes $H_i$ passing through $t$ because the restriction of
$s$ and $s^{\prime}$ to $H_i$ are multiple of $t^2$ and $E_{H_i}(1)$
can't have a section not proportional to the section of $E_{H_i}$
because there are no 5-jumping lines according to the
\ref{maxsauteuses}. So denote by $Q$ the quadric surface such that
$Z_{s\wedge s^{\prime}}=\bigcup\limits_{i=1}^k H_i \cup Q$.

First notice that for a general plane $H$ containing $t$, the sections
$s$ and $s^{\prime}$ gives, after the division by $t^2$, 2 sections of
$E_H(1)$ which are proportional on a conic which must be $Q\cap H$.
Furthermore, this conic must contain $t$ because $t$ is 4-secant to the
vanishing locus of any section of $E_H(1)$. Denote by $d_H$ the line
such that $t \cup d_H=Q \cap H$.

Next, assume that $\Gamma$ is drawn on $S$ as in the hypothesis of
\ref{courbetris}. The curve $T_tG\cap S$ is made of multi-jumping lines
meeting $t$, and is an hyperplane section of $S$. We want here to
understand this curve by cutting it with planes containing $t$. First
remark that there is no irreducible component of 
$T_tG\cap S$ included in some $\beta$-plane containing $t$, otherwise
its associated plane would be a jumping one according to the
\ref{infinitebis}, so it would be one of the $H_i$, but the $H_i$
contain two 3-jumping lines, so they contain only a finite number of
multi-jumping lines using again the \ref{infinitebis}. So the points of
$T_tG\cap S$ which represent lines lying in a stable plane $H$
containing $t$ make  a dense subset of $T_tG\cap S$. Now we can
notice that in a general plane $H$ containing $t$, the only possible
multi-jumping lines are trisecant to any vanishing locus of a section
of $E_H(1)$, so they must be $t$ or $d_H$. But $d_H$ is always in $Q$,
so the curve $T_tG\cap S$ must be the rulings of $Q$, and it would
imply $\deg S=2$, then the curve $\Gamma$ of degree 6 and genus 2
could not lie on $S$ because $h^0\OO_{\Gamma}(1)=6$, and it gives the
proposition \ref{courbetris}.

\bigskip
We'd like now to enlarge the result \ref{courbetris} to the case where
the previous curve $\Gamma$ contains 4-jumping lines when $c_2=5$.
\begin{proposition}
\label{sauteuses4}When $c_2=5$, the surface $S$ don't contain a curve
of jumping lines of order 3 in general with some 4-jumping lines
\end{proposition}

So we assume that there is such a curve, and let $\Gamma$ be an
irreducible curve of jumping lines of order 3 in general with at least
a 4-jumping line $q$. Denote by $\Sigma$ the ruled surface associated
to $\Gamma$.

First notice that $\chi(E(2))=0$, and the existence of $q$ implies
$h^1E(2)\neq 0$, so $E(2)$ has a section $s$ vanishing on some degree
9 curve $Z$. Furthermore, any 3-jumping lines is 5-secant to $Z$. The
part of the proof of the \ref{courbetris} bounding the number of
jumping planes containing a general ruling $t$ of $\Sigma$ is still
valid in this situation, so we still have $\deg \Gamma \leq 6$.

a) Now notice that $\Gamma$ can't be a conic. In fact this is the worst
case because the bundles of the family \ref{famille} have such a conic
of 3-jumping lines with four 4-jumping lines. So the required
contradiction need to use the existence of $S$. In that case
$\Sigma$ is a quadric so we have:
\begin{center}
$0\longrightarrow \OO_{\proj_1 \times \proj_1}(-a,3) \longrightarrow
E_{\Sigma} \longrightarrow \OO_{\proj_1 \times \proj_1}(a,-3)
\longrightarrow \complex ^n\longrightarrow 0$
\end{center}
which gives $10=6a+n$ where $n$ is the number of 4-jumping lines on
$\Gamma$. So $a$ would be 0 or 1 then $h^0E_{\Sigma}(2)\geq 12$. The
resolution of the quadric twisted by $E(2)$ gives the sequence:
\begin{center}
$0 \longrightarrow H^0E(2)\longrightarrow H^0E_{\Sigma}(2)\longrightarrow H^1E$
\end{center}
So \mbox{$h^0E(2)\geq 4$}, and we can solve this case with the \ref{multisautfamille}.

\bigskip
For the other possible degrees of $\Gamma$, the situation is easier
because none of the ruling of $\Sigma$ can meet $q$ because in any
plane $H$ containing $q$, this line is according to the \ref{secantes}
\prefix{\deg Z-1}secant to the
0-dimensional vanishing locus $Z$ of a section of $E_H$ or of
$E_H(1)$, and any ruling of $\Sigma$ is at least 3-secant to $Z$.

So if $\Sigma$ has a double locus, it must meet $q$ in a pinch point, so $q$
is a torsal line. In other words, there is a line $L$ included in $G$
meeting $\Gamma$ in length at least 2 around $q$. So, the support of
the sheaf $R^1q_{*}p^{*}E(\tau )_L$ must contain twice the point $q$.
So this sheaf can't be equal to its restriction at the reduced point
$q$, but it has rank 2 at $q$, so $R^1q_{*}p^{*}E(\tau )_L\geq 3$
which is not possible around a 4-jumping line, because the plane
associated to $L$ can't be stable, hence according to the
\ref{secantes}, it is for $c_2=5$ the same problem as the study of the
\ref{pinceaux} for $c_2=4$ around a 3-jumping line, and the proof of
\ref{pinceaux} was still valid for $c_2=4$ in this situation.

If $\Sigma$ has no singular locus, then as it can't be a quadric, it
is the tangential surface of a skew cubic, hence all its lines are
torsals, so we can conclude as previously. \fin

We also have the following:
\begin{corollary}
 \label{c2=5deg}For $c_2=5$, any integral curve in $G$ made only of 2-jumping
 lines has a degree multiple of 4. The surface $S_{red}$ has a degree
 multiple of 4.
\end{corollary}

Indeed, any integral curve of degree $d$ made only of 2-jumping lines
gives using the previous notations, a ruled surface such that
$5d=4a+4e$. But the results \ref{maxsauteuses}, \ref{fini4},
\ref{courbetris}, and \ref{sauteuses4} imply the general hyperplane
section of $S_{red}$ don't contain any \prefix{k\geq3}jumping lines. 
\section{Trisecant lines to $S$ and applications}

In the first section we obtained an interpretation of some trisecant lines to
$S$. The key application of this interpretation is the folowing proposition:
\subsection{Trisecant lines to the general hyperplane section of $S$}
\begin{proposition}
\label{c2=5tris} When $c_2\leq 5$, the general hyperplane section of $S_{red}$ has no
trisecant lines. 
\end{proposition}

Indeed, on one hand the generic hyperplane section of $S_{red}$ contains no
\linebreak\prefix{k\geq 3}jumping lines according 
to the section \ref{fini3sautet+}, and on the other hand, the lines of
$\proj_5$ which are in a hyperplane make a 2-codimensional subscheme
of the Grassmann manifold $G(1,5)$. So if the general hyperplane
section of $S_{red}$ had a trisecant
line then $S_{red}$ would have a 2-parameter family of
trisecant lines meeting $S$ in 2-jumping lines. But any of those trisecant lines
would give a jumping plane according to the \ref{pinceaux}, and as
there is only a finite number of such trisecant lines in a same
jumping $\beta$-plane according to the \ref{pinceauxsauteurfini}, we would have
a 2 parameter family of jumping planes which would contradict the
\ref{coanda}. \fin

\subsection{The cases where $S$ is ruled}\label{reglee}

\begin{itemize}
\item Let's first consider the case where $S$ is a ruled surface when
$c_2=5$.
\end{itemize}
Every ruling of $S$ must 
contain a 3-jumping line according to the  \ref{infinitebis}, but
there is at most a finite number of 3-jumping lines on $S$ from the
\ref{courbetris} and the \ref{sauteuses4}. So $S$ must be a cone with
vertex a 3-jumping line $t$ , and this cone is in the $\proj_4$
constructed by projectivisation of the tangent space $T_tG$. According
to the \ref{infinitebis}, any plane contain at most 
one pencil of multi-jumping lines, so $S_{red}$ must have bidegree
$(1,\beta)$, and the bound on $\beta$ of the \ref{basdeg} and the
\ref{c2=5deg} imply that $\beta=3$. But the projection from $t$ of
this cone gives a curve of bidegree $(1,3)$ in the quadric obtained by
projection of $T_tG\cap G$. So there are many planes containing 3
distinct rulings of $S_{red}$, then the generic hyperplane section of
$S_{red}$ would have trisecant lines which would contradict the
\ref{c2=5tris}.

\begin{itemize}
\item So assume here that $S$ is a ruled surface and that $c_2=4$.
\end{itemize}

Denote by $C$ the reduced curve described by the center of the pencils
of lines associated to the rulings of $S$.
\begin{itemize}
\item[-] If $\deg C\geq 3$ 
\end{itemize}
A general hyperplane $H$ meeting $C$ in distinct points contain a line of
the congruence $S$ passing through each point of $C\cap H$. As $\beta
\leq 2$ and $\deg C\geq 3$, those lines must be bisecant to $C$, and
the other bisecant to $C$ are not in the congruence $S$. Hence
$S$ is a join between 2 components of $C_1$ and $C_2$ of $C$ with
$\deg C_1=1$ because $S$ must be ruled. So every plane $H$ containing
$C_1$ contain a pencil of multi-jumping lines as in the
\ref{infinitebis} which is centered at a point of $C_2\cap H$.  So
those plane can't be stable and 
according to the \ref{dteplansaut} $\ideal^2_{C_1}\otimes E(2)$ has a section.
Furthermore, this section vanishes at least with multiplicity 4 on $C_1$ and
also on $C_2$, but it has degree 8 and must be connected, then $S$
would be a $\beta$-plane which is not possible. 
\begin{itemize}
\item[-] If $\deg C=1$
\end{itemize}
Then $S$ is a cone with vertex the point representing $C$, and once
again every plane containing $C$ would contain infinitely many
multi-jumping lines, hence from the \ref{dteplansaut}
$\ideal^2_{C}\otimes E(2)$ has a section $s_0$ vanishing at least with
multiplicity 4 
on $C$ and on the curve $C^{\prime}$, where $C^{\prime}$ is the curve
obtained by the center of the pencils of multi-jumping lines
which are in planes containing $C$. So we must have $\deg C^{\prime}
=\alpha =1$ and the connexity of the vanishing locus imply that $C=C^{\prime}$.
So, in any plane $H$ containing $C$, the bundle $E_H$
has a section $s_H$ of type \ref{infinitebis} such that the connected
component of $Z_{s_H}$ is on $C$, and any section $\sigma$ of $E_H(1)$
is proportional  
to $s_H$, otherwise its vanishing locus $Z_{\sigma}$ would be in the line
$Z_{\sigma \wedge s_H}$ which would be a 4-jumping line, but it would contradict
the \ref{infinitebis}. Therefore, every section of $\ideal_C
\otimes E(2)$ must be a section of $\ideal_C^2 \otimes E(2)$.

Furthermore, if we blow up 
$\proj_3$ along $C$, the bundle $q^{\prime}_*p^{\prime *}E$ has rank
1, so $h^0\ideal_C^2 \otimes E(2)=1$. This proves that the sections of
$E(2)$ have no base curves, and the number of base points is the number
of residual points in the intersection of the 3 quartics of
$H^0\ideal_{Z_{s_0}}(4)$ which is 0 from [Fu] p155. So from Bertini's
theorem, there is a smooth section $s$ of $E(2)$, such that $Z_s$ is also connected
from $h^1E(-2)=0$, and $Z_s \cap C=\emptyset$. The quartic surface
$Z_{s_0 \wedge s}$ must cut any 
plane $H$ containing $C$ in twice $C$, and a conic which must be singular in some
$P_H\in C$
because it contains $Z_{s_H}$. Furthermore, this conic doesn't contain $C$
because $Z_s\cap C= \emptyset$ and the lines through $P_H$ are $0$- or
4-secant to $Z_s$. So this quartic is a ruled
surface with $C$ as directrix, but  $C$ is in its singular locus so
there are rulings of $Z_{s_0 \wedge s}$ through $P_H$ which are not in
$H$. So $C$ must be a triple curve of $Z_{s_0 \wedge s}$, which
imply that this quartic has a rational basis. But it 
contradicts the formula of Segre (Cf [G-P2]) applied to $Z_s$ because it is a
degree $8$ elliptic curve 4-secant to the rulings of $Z_{s_0 \wedge s}$. 

\begin{itemize}
\item[-]  When $\deg C=2$
\end{itemize}
The plane $H$ containing the conic $C$ must then contain infinitely many
multi-jumping lines which have to contain a same point $O$ according to
the \ref{infinitebis}. But the lines of the congruence are secant to
$C$, so only one of them contain a fixed point $P$ of $H-\{C,O\}$, then
$\alpha =1$. We had already eliminated the cases of the congruences
$(1,1)$ ($\deg C=1$), and the congruences $(1,2)$ are the joins
between $C$ and a line $d$ cutting $C$ in 1 point (Cf [R]). As
previously, every plane containing $d$ would be of type
\ref{infinitebis}, and the  section of $\ideal_d^2 \otimes E(2)$ would
vanish on $d\cup C$ at least with multiplicity 4, which is not
possible because it has degree 8. \fin 

\subsection{``No'' trisecant cases}
The aim of this section is to get rid of the remaining cases.
\begin{lemma}\label{degS=4}
  For $c_2=4$ and 5, the congruence $S_{red}$ can't have degree 4.
\end{lemma}
The degree 4 surfaces are classified in [S-D], and as we already
eliminate the cases of ruled surfaces in \S\ref{reglee}, only the
complete intersection 
of 2 quadrics in $\proj_4$ and the Veronese are remaining.

\begin{itemize}
  \item  For $c_2=4$
\end{itemize}
The Veronese has bidegree $(3,1)$ or $(1,3)$ so it would contradict the
\ref{basdeg}. If $S$ is the complete intersection of $G$ with another
quadric and some $\proj_4$ then $S$ has bidegree $(2,2)$ in G (Cf
[A-S]). Furthermore, in that case $S$ is locally complete
intersection, so we can compute from the \ref{exresiduel} that the
scheme of multi-jumping lines $M$ is made of $S$ with a non empty residual
curve $C$. There is a 2 dimensional family of $\beta$-planes which
meet $C$, so from the \ref{coanda} its general element must be a
stable plane, then it must cut $S$ in good dimension according to the
\ref{infinitebis}. But those $\beta$-planes will contain $\beta +1=3$
multi-jumping lines (with multiplicity) which is not possible when
$c_2=4$ (Cf \ref{basdeg}). We can remark that it is still true when
$C$ is drawn on $S$. Indeed, one has the exact sequence:
\begin{center}
  $0 \longrightarrow  \LL \longrightarrow \OO_M  \longrightarrow \OO_S
  \longrightarrow 0$ 
\end{center}
where $C$ is the support of $\LL$. Restrict it to a stable
$\beta$-plane $h$ which meet $C$ to have:
\begin{center}
 $\Tor_1(\OO_S,\OO_h)  \longrightarrow \LL \otimes \OO_h
 \longrightarrow \OO_{M\cap h} \longrightarrow  \OO_{S\cap h} \longrightarrow 0$  
\end{center}
As $h$ is stable, it cuts $S$ in dimension 0, so
$Tor_1(\OO_S,\OO_h)=\ideal_S \cap \ideal_h /\ideal_S \ideal_h=0$,
which still gives length$\OO_{M\cap h}\geq 3$.
\begin{itemize}
  \item When $c_2=5$
\end{itemize}
The complete intersection of 2 quadrics in $\proj_4$ is a Del-Pezzo
surface isomorphic to $\proj_2$ blown up in 5 points embedded by
$3L-E_1-...-E_5$. According to the \S\ref{fini3sautet+}, the surface
$S$ contains at most a finite number of \prefix{k\geq 3}jumping lines,
so we can find on $S$ some cubic curves made only of 2-jumping lines,
which contradicts the \ref{c2=5deg}. Similarly, in the Veronese case,
there would have many conics made only of 2-jumping lines which also
contradicts the \ref{c2=5deg}. \fin

\begin{proposition}
  \label{c2=5courbe}
  For $c_2=4$ or $5$, the scheme of multi-jumping lines of a
  mathematical instanton is a curve in $G$.
\end{proposition}

In fact the case $c_2=4$ is already done because $S$ can't have degree
4 (\ref{degS=4}), and the bound of its bidegree of the \ref{basdeg}
imply that $\alpha =1$ or $\beta =1$. But those congruences are ruled
(Cf [R]) and we can conclude with the \S\ref{reglee}.

\bigskip
For $c_2=5$, let's first prove the following:
\begin{lemma}
  The set of points $p$ of $S$ such that $\dim T_pS_{red}=4$ is finite.
\end{lemma}
Assume that there are infinitely many such $p$, then the  general
hyperplane section contains one of them, so in general, the surface $S_{red}$
has no trisecant through $p$ (Cf \ref{c2=5tris}), so the curve $T_tG \cap
S_{red}$ must be an union of lines containing $p$, and as $S$ is not
ruled (Cf \S\ref{reglee}), it can happen only at finitely many $p$
which proves the lemma. \fin

\bigskip
So a general hyperplane section $\Gamma$ of $S_{red}$ is locally complete
intersection with at most finitely many $p$ such that $\dim
T_p\Gamma=2$. According to the \ref{c2=5tris}, the curve $\Gamma $ has
no trisecant, and we can adapt LeBarz's formula of [LeBa2]. In fact,
when projecting $\Gamma $ to $\proj_2$, the singularities of $\Gamma $
won't gives embedded points, so they have to be removed from the
contribution of the apparent double points. Hence the number of
trisecant to $\Gamma$ is $(d-4)((d-2)(d-3)-6\pi)$ where $\pi$ is the
arithmetic genus of $\Gamma $.

So we have from the \ref{c2=5tris} and \ref{degS=4} that
$(d-2)(d-3)-6\pi=0$. Thus $d-1$ is not a multiple of 3, and then
$\Gamma $ is locally complete intersection with maximal genus in
$\proj_4$, because if $\epsilon$ is the remainder of the division of
$d-1$ by 3, and $m=\left[ \frac{d-1}3\right] $, then
$\frac{(d-2)(d-3)}6 =\frac{\epsilon ^2-3\epsilon +2}6+
\frac{3m(m-1)}2+ m\epsilon $ because $\epsilon\neq 0 $.

According to the \ref{degS=4} and the \ref{c2=5deg}, we have $\deg
\Gamma =8$, and as it has maximal genus in $\proj_4$, it must be a
canonical curve, and it can't be trigonal (ie have a $g^1_3$) because it
has no trisecant. So according to Riemann-Roch $\Gamma$ hence $S$ must
be the complete intersection of 3 quadrics. So $S$ is the Kummer K-3
surface and has bidegree $(4,4)$ (Cf [A-S]). We can compute from the
\ref{exresiduel} that the scheme of multi-jumping lines has a residual
curve because $S$ is locally complete intersection, and we construct
as in the \ref{degS=4} when $c_2=4$ a 2 parameter family of
$\beta$-planes containing 5 multi-jumping lines (with multiplicity),
which contradicts the \ref{coanda} or the \ref{basdeg}. \fin
\section{Applications to moduli spaces}\label{applications}

\subsection{The construction of Ellingsrud and Stromme}

Let us recall here the construction made in [E-S], and show how the
previous results could be used to study I$_n$.

Let $N$ be a point such that there exists an instanton $E$ which
has only a one dimensional family of jumping lines through $N$ which
are all of order $1$. Denote by 
 $U_N$ the subscheme of I$_n$ made of 
instantons which don't have a multi-jumping line through $N$ and
which have a non jumping line through $N$. The result
\ref{c2=5courbe} and Grauert-Mullich's theorem implies that the
$(U_N)_{N\in \proj_3}$ cover $I_n$ for $n=4$ or $5$.

To describe $U_N$, first blow up $\proj_3$ at $N$, denote by $Y$ the
  plane parameterizing the lines containing $N$, $p$ and $q$ the
  projections on $\proj _3$ and $Y$, and $\tau=p^*\Op(1)$,
  $\sigma=q^*\OO_Y(1)$, and by $V=H^0(\OO_Y(1))$. Let $H$ be
a $\complex $-vector space of dimension $n$. ($H$ will stand for
$H^2E(-3)$). A point $E$ of  $U_N$ is characterized  (after the choice of
an isomorphism:  $H \simeq H^2E(-3)$) by the following data where $F$
  is naturally $q_*p^*E$ and $\theta (1)$ is $R^1q_*p^*E(-1)$:

\begin{center}
$\left\{ 
\begin{array}{cllc}
1) & & m\in V\otimes S_2H\dual & \\ 
2) & \text{a surjection: }& 2\OO_{\proj(V\dual)}\rightarrow
\theta (2)\rightarrow 0 & \text{(denote by } F \text{ its kernel)} \\ 
3) & \text{and a surjection: }& q^{*}F\dual\rightarrow q^{*}\theta
(\sigma +\tau )\rightarrow 0 & 
\end{array}
\right. $
\end{center}

In fact, the symmetric map $m$ is just the restriction to the $\alpha $-plane
associated to $N$ of the $2n+2$ rank element of $\stackrel{2}{\Lambda }H^0(\OO%
_{\proj\!_3}(1))\otimes S_2H\dual$ of Tjurin, (Cf [T2]
or [LeP]) and it can be obtained from the following sequence coming from
Beilinson's spectral sequence.
\begin{equation}
  \label{reseau}
  0\longrightarrow H\otimes \OO_Y\stackrel{m}{
  \longrightarrow }H\dual\otimes \OO_Y(1)\longrightarrow
\theta (2)\longrightarrow 0
\end{equation}
We have also to recall from [E-S] that the data 3) exists only under the
condition that the restriction of $F$ to $C$ splits. To be more
precise, the restriction of the data 2) to $C$ gives a sequence
\begin{equation}
  \label{Fc}
  0 \longrightarrow \theta \dual(-1) \longrightarrow F_C \longrightarrow \theta
\dual(-2)\longrightarrow0
\end{equation}
which has to split for the existence of
the data 3). This is the main difficulty of
this description of $U_N$. So let's first globalize data 2).
                                
The group $P=^{GLH}/_{\pm{}Id}$ acts on $H$ by $^t p.h.p$, and it also
acts on the previous data according to the exact sequence (\ref{reseau}).
Furthermore, 2 elements of $U_N$ are isomorphic if and
only if they are in the same orbit of this action. 

Our reference about properties and existence of the moduli space of
those theta characteristic will be [So]. A net of quadrics will be
called semi-stable (stable) if and only if, for all non zero totally
isotropic (for the net) subspace $L$ we have: \hbox{$\dim L+\dim
L^{\bot }\leq \dim H$} ($< \dim H$). Let $\Theta $ be the quotient
\hbox{$(V\otimes S_2H \dual)^{ss}/P$}, and $\Theta ^s_{lf}$ be the stable
points which represent locally free sheaves, and
$\Theta _0$ be the dense open subset of $\Theta $ made of isomorphic
classes of  $\theta $-characteristic which have a smooth support. The
space  $\Theta $ is irreducible of dimension $\frac{n(n+3)}2$, and
normal with rational singularities (Cf [So]
th\'eor\`eme 0.5), and $\Theta ^s_{lf}$ is smooth (Cf [So]
th\'eor\`eme 0.4) because when $\theta$ is locally free, the notion
of stable nets coincide with the stability of $\theta $ in the sheaf
meaning (without the locally free assumption, the last notion is stronger).

By construction, any element of $U_N$ has a non-jumping line through
$N$. This implies that the associated net is
semi-stable. Furthermore by definition of $U_N$, the
$\theta$-characteristic considered are locally free. So we have a morphism:
$$I_n \supset U_N \longrightarrow \Theta _{lf}$$
The key of the construction is to understand the fiber of this
morphism. We have thus to globalize the data 2) and then, to show
that for $n \leq 11$ this map is dominant, so we will
have to understand under what condition the data 3) exists.

\bigskip
Define by  
$G^{\prime }_{\theta}=G\big(2,H^0(\theta (2))\big)$  the Grassmann
manifold of 2-dimensional subspace of $H^0(\theta (2))$, and
$G^{\prime \prime }_{\theta}$ be the open subset of
$G^{\prime}_{\theta}$  made of pairs
of sections of $\theta (2)$ with disjoint zeros. Denote by
$K_{\theta}$ the tautological subbundle of $\OO_{G^{\prime
    \prime}_{\theta}} \otimes H^0(\theta(2))$.
The existence condition of the data 3) was already
identified in [E-S] \S{4.2} as the vanishing locus of a map:
$$\delta_{\theta} :\OO_{G^{\prime \prime }_{\theta}}\longrightarrow H^1(\OO%
_C(1))\otimes \stackrel{2}{\Lambda }K\dual_{\theta}$$
obtained by the following way:

\bigskip
We pick from data 2) the exact sequence:
\begin{center}
  $0\longrightarrow F\longrightarrow \OO_{\proj(V\dual%
)}\boxtimes K_{\theta}\longrightarrow \theta (2)\boxtimes \OO_{G^{\prime
\prime }_{\theta}}\longrightarrow 0$
\end{center}
whose restriction to $C\times{}G^{\prime \prime }$ (where $C$ is
$\theta$ 's support) gives:
\begin{center}
  $0\longrightarrow \OO_C(1)\boxtimes \OO_{G^{\prime \prime
}_{\theta}}\longrightarrow F(\theta (2)\boxtimes \OO_{G^{\prime \prime
}_{\theta}})\longrightarrow \OO_C\boxtimes \stackrel{2}{\Lambda }%
K_{\theta}\longrightarrow 0$
\end{center}
Pushing down this sequence twisted by $\OO_C\boxtimes
\stackrel{2}{\Lambda }K\dual_{\theta}$  with the second projection of
$C\times{}G^{\prime \prime }_{\theta}$ direct image's functor, we obtain the
boundary $\delta_{\theta} $.

\bigskip
Denote by $Z_{\delta_{\theta}} $ the vanishing locus of $\delta_{\theta} $. Then,
according to [E-S], the fiber $U_{N,\theta}$ is the product of
$Z_{\delta_{\theta}} $ by $\proj_4$. 
$$I_n \supset U_N \longrightarrow \Theta_{lf} \hspace{2cm} \forall \theta
\in \Theta_{lf},U_{N,\theta}=\proj_4 \times Z_{\delta_{\theta}}$$ 

Let $\Theta _0$ be the subspace of $\Theta $ made of isomorphic
classes of $\theta$-characteristic with smooth support. 
Ellingsrud and Stromme had also remarked ([E-S] \S{4.2}) that for any
$\theta$ in $\Theta_0$, the
bundle $\stackrel{2}{\Lambda }K_{\theta}\dual$ generates
$Pic G_{\theta} ^{\prime \prime }$, and that for $n
\leq 11$ $\Theta _0$ is included in $U_N$ 's image. This image is thus
irreducible, normal and $\frac{n(n+3)}2$ dimensional. 

Furthermore, if we denote by $U_N ^s$ the elements of $U_N$ which have
a stable image in $\Theta $, we obtain that the image of $U_N ^s$ is
irreducible and smooth. So we need to prove the following lemma: 

\begin{lemma}
\label{recstable}For $n=4$ and $5$, the $U_N^s$ make a covering of
I$_n$ when $N$ travels through $\proj_3$ 
\end{lemma}

We have to understand the non stable nets which could occur in the
image of $U_N$. Let's recall from [So] that $\theta $ is not stable if
there exists a totally isotropic (for the net $m$ associated to
$\theta$) subspace $L\subset H$ such that $\dim L+\dim L^{\perp }\geq
n$. As there exists a smooth quadric in the net, $L$ has to be at most
2 dimensional. The case $\dim L=1$ has already been avoided in [E-S]
\S{5.10} using Barth's condition $\alpha 2$ of [Ba3]. So we have just
to study the case $\dim L=2$.

If there exists such a 2-dimensional space, then all the quadrics of
the net can be written in the following way: $\sum_{i=1}^{i=3}Y_i\left( 
\begin{array}{cl}
  \begin{picture}(80,100)(0,20)
    \put(15,7){0}
    \put(0,-35){\line(1,0){80}}
    \put(0,100){\line(1,0){130}}
    \put(80,20){\line(1,0){50}}
    \put(80,-35){\line(0,1){55}}
    \put(0,-35){\line(0,1){135}}
    \put(130,20){\line(0,1){80}}
   \end{picture}
& \hfill ^{^\tau A_i} \\ 
_{A_i}\hfill & B_i
\end{array}
\right) $ where $A_i$ is a $2\times{}2$ matrix, and where $B_i$ is a
$2\times{}2$ matrix in case $n=4$ and a $3\times{}3$ matrix in case
$n=5$. Their determinant must then be a multiple of $(\det (
\sum_{i=1}^{i=3} Y_iA_i))^2$. The curve of jumping lines of any
preimage of this net must contain this double conic. Let $
A_N=$\hbox{\tiny $\left(\begin{array}{cc}
  a_N & b_N \\
  c_N & d_N
\end{array}\right)$}$=\sum_{i=1}^{i=3}Y_iA_i$, and first prove the
following lemma before continuing the \ref{recstable}:
\begin{lemma}\label{conicsing}
  If the conic $det(A_N)$ is singular in $p$, then either $p$ represents a
  multi-jumping line or $N$ is in the vanishing locus of a section of $E_H$
  for some jumping plane $H$ 
\end{lemma}
The first case is obtained when $a_N,b_N,c_N,d_N$ represent lines
containing a same point $p$, and the other cases are reduced to the 2
cases where
$a_N,b_N,d_N$ or $b_N,c_N,d_N$ is a basis of $V$, because we are
allowed to change the basis of $L$. Computing the remaining coefficient in
this basis, and using that $A_N$ is singular, we are reduced by
changing the basis of $L$ to the case where $a_N$ and $c_N$ are
proportional. But in this case, $a_N$ divides all the $(n-1)\times (n-1)$
minors except one. It means that the line $a_N$ cuts the multi-jumping
lines in length $n-1$. The computation of \ref{pinceaux} also proves
that this is not possible in a stable plane, and the
\ref{pinceauxsauteurfini} that it is impossible in a jumping 
plane $H$ where $N$ is not in the vanishing locus of the section of
$E_H$, because there is no jumping line of order $n$ through $N$ when
$a_N,b_N,d_N$ or $b_N,c_N,d_N$ is a basis of $V$. This proves the
lemma \ref{conicsing}. \fin

\bigskip
To obtain the \ref{recstable}, we have now to prove that there exists
no instanton such that for 
every $N$ in $\proj _3$ the net has such a 2 dimensional totally
isotropic space $L_N$. If there was such a
bundle, then, its hypersurface of jumping line would contain a
double quadratic complex. Let $Q$ be this quadratic complex, and
$\Sigma$ be the set of $N\in \proj_3$ such that the $\alpha$-plane
$\alpha_N$ doesn't cut $Q$ in a smooth conic,  and let
$\Sigma^{\prime}$ be the set of singular rays (i.e: the points of $G$
which are singularities of $Q\cap\alpha_N$ for some $N$). As $\Sigma$
is the degeneracy locus of a square matrix, it is at least 2
dimensional. We can show that for any $Q$, the set $\Sigma^{\prime}$
is at least 2 dimensional too: let $\sum_{i=0}^5 x_i^2=0$ and
$\sum_{i=0}^5 k_ix_i^2=0$ be the equations of $G$ and $Q$, and consider
the third quadric $Q^{\prime}: \sum_{i=0}^5 k_i^2x_i^2=0$, then
$ G\cap Q\cap Q^{\prime} \subset \Sigma^{\prime}$ because $(x_i)\in
 G\cap Q\cap Q^{\prime} $ implies that $(k_ix_i)$ is a line meeting the
line $(x_i)$ in some point $N$, and for any $(z_i)$ in the
$\alpha$-plane $\alpha_N$, the equation $(z_i)$ cut $(k_ix_i)$
which is $\sum_{i=0}^5 k_i x_iz_i=0$ means also $(z_i)\in T_{(x_i)}Q$. 

Then we can conclude using the lemma \ref{conicsing}, proposition
\ref{c2=5courbe}, and Coanda's theorem ([Co]), remarking in the cases
where $N$ is on the jumping section of $E_H$ for some jumping-plane
$H$, that $E$ can't be a special t'Hooft bundle, otherwise there
would exist a $n$-jumping line through $N$ which as been excluded in
the demonstration of \ref{conicsing}. \fin

\subsection{Description of the boundary $\delta_{\theta}$}

We'd like now to understand more explicitly the fiber of
$Z_{\delta_{\theta}}$ over some $\theta $. In other words, we need
  a description of the splitting condition of the restriction of $F$
  to $C$. The following will consist of 2 descriptions of this condition. The
  first one is explicit in function of the net of quadrics, and enables
  to prove that the fiber $Z_{\delta_{\theta }}$ is the
  complete intersection of $G^{\prime \prime}_{\theta}$ with some hyperplanes
  This could have been conjectured from Ellingsrud and Stromme viewpoint
  because of the map $\delta$ of the previous subsection, but we can't
  conclude directly from this because there is some torsion in $Pic
  G^{\prime\prime}_{\theta}$ .
  Furthermore, we will need a second description to understand the ramification
  of this morphism in term of vector bundle and to get informations on
  I$_n$.
  
\begin{itemize}
\item First description
\end{itemize}
Let $K$ be the 2-dimensional vector space generated by 2 sections of
$\theta (2)$ with disjoint zeros. Let $F$ be the vector bundle over
the exceptional plane $Y= \proj (V\dual)$ defined by the exact sequence:
\begin{equation}\label{F}
0 \longrightarrow F \longrightarrow K \otimes \OO_Y \longrightarrow
\theta (2) \longrightarrow 0
\end{equation}
Let $S$ be the kernel of the map $H\dual \otimes \theta (3)
\rightarrow \theta ^2(4)$ obtained when twisting by $\theta (2)$ the
exact sequence (\ref{reseau}) of the net of quadric.

As ${\cal E}xt^1(\theta (2),\OO_Y)=\theta (1)$, we can compute the
kernel of the restriction of (\ref{F}) to $C$ by dualizing (\ref{F}).
This yields to $Tor_1(\theta (2),\theta (2))=\OO_C(1)$. We can now
deduce from (\ref{reseau}) the exact sequence:
$$0\longrightarrow \OO_C(1)\longrightarrow H\otimes \theta
(2)\longrightarrow H\dual\otimes \theta (3)\longrightarrow \theta
^2(4)\longrightarrow 0$$
Taking there global sections, we get the following diagram where
$\Sigma $ is by definition the cokernel of \hbox{$H\otimes H^0(\theta
  (2))\rightarrow H\dual\otimes H^0(\theta (3))$}.

\begin{center}
{\begin{picture}(3000,1100)
\setsqparms[1`1`0`1;800`300]
\putsquare(1100,550)[H\otimes H^0(\theta (2))`\phantom{H^0(S)}`H\otimes H^0(\theta (2))`\phantom{H^{{\rm v}}\otimes H^0(\theta (3))};`\wr``b]
\setsqparms[1`1`1`1;600`300]
\putsquare(1900,550)[H^0(S)`H^1(\OO_C(1))`H^{{\rm v}}\otimes H^0(\theta (3))`\Sigma;```]
\setsqparms[0`1`1`1;600`300]
\putsquare(1900,250)[``H^0(\theta ^2(4))`H^0(\theta ^2(4));```\sim]
\putmorphism(400,850)(1,0)[{H^0(\OO_C(1))}`\phantom{H\otimes H^0(\theta (2))}`]{700}1a
\putmorphism(0,850)(1,0)[0`\phantom{H^0(\OO_C(1))}`]{400}1a
\putmorphism(2500,850)(1,0)[\phantom{H^1(\OO_C(1))}`0`]{400}1a
\putmorphism(2500,550)(1,0)[\phantom{\Sigma}`0`]{400}1a
\putmorphism(1900,250)(0,-1)[`0`]{250}1a
\putmorphism(2500,250)(0,-1)[`0`]{250}1a
\putmorphism(1900,1100)(0,-1)[0``]{250}1a
\putmorphism(2500,1100)(0,-1)[0``]{250}1a
\end{picture}}
\end{center}
On the other hand, we can compute the following diagram displaying the
sequence (\ref{F}) horizontally and the sequence (\ref{reseau})
vertically:

\begin{center}
{\begin{picture}(2200,1100)(0,0)
\setsqparms[1`1`1`1;600`300]
\putsquare(500,550)[F\otimes H`K\otimes H\otimes\OO_Y`F\otimes H^{{\rm v}}(1)`K\otimes H^{{\rm v}}(1);```]
\setsqparms[1`0`1`1;600`300]
\putsquare(1100,550)[\phantom{K\otimes H\otimes\OO_Y}`H\otimes \theta(2)`\phantom{K\otimes H^{{\rm v}}(1)}`H^{{\rm v}}\otimes \theta(3);```]
\setsqparms[0`1`1`1;600`300]
\putsquare(500,250)[``F\theta(2)`K\otimes \theta(2);```]
\setsqparms[0`0`1`1;600`300]
\putsquare(1100,250)[``\phantom{K\otimes \theta(2)}`\theta^2(4);```]
\putmorphism(0,850)(1,0)[0`\phantom{F\otimes H}`]{500}1a
\putmorphism(0,550)(1,0)[0`\phantom{F\otimes H^{{\rm v}}(1)}`]{500}1a
\putmorphism(1700,850)(1,0)[\phantom{H\otimes\theta(2)}`0`]{500}1a
\putmorphism(1700,550)(1,0)[\phantom{H^{{\rm v}}\otimes\theta(3)}`0`]{500}1a
\putmorphism(1700,250)(1,0)[\phantom{\theta^2(4)}`0`]{500}1a
\putmorphism(1100,250)(0,-1)[`0`]{250}1a
\putmorphism(1100,1100)(0,-1)[0``]{250}1a
\putmorphism(1700,250)(0,-1)[`0`]{250}1a
\end{picture}}
\end{center}
We obtain by passing there to global sections the following:

\begin{center}
{\begin{picture}(2600,1600)
\setsqparms[1`1`1`1;800`400]
\putsquare(500,250)[H^{{\rm v}}\otimes V\otimes K`H^0(\theta (3))\otimes H^{{\rm v}}`H^0(\theta (2))\otimes K`\Sigma ;a`c``]
\putsquare(1300,250)[\phantom{H^0(\theta (3))\otimes H^{{\rm v}}}`H^1(F(1))\otimes H^{{\rm v}}`\phantom{\Sigma}`H^1(F\theta (2));```f]
\putsquare(500,650)[H\otimes K`H^0(\theta (2))\otimes H`\phantom{H^{\rm v}\otimes V\otimes K}`\phantom{H^0(\theta (3))\otimes H^{{\rm v}}};``b`]
\setsqparms[1`1`1`0;800`400]
\putsquare(1300,650)[\phantom{H^0(\theta (2))\otimes H}`H^1(F)\otimes H``;```]
\putmorphism(2100,650)(1,0)[\phantom{H^1(F(1))\otimes H^{{\rm v}}}`0`]{500}1a
\putmorphism(2100,1050)(1,0)[\phantom{H^1(F)\otimes H}`0`]{500}1a
\putmorphism(500,250)(0,-1)[`0`]{250}1a
\putmorphism(1300,250)(0,-1)[`0`]{250}1a
\putmorphism(500,1350)(0,-1)[0``]{300}1a
\putmorphism(2100,1350)(0,-1)[H^0(F\theta (2))``]{300}1a
\putmorphism(2100,1600)(0,-1)[0``]{250}1a
\end{picture}}
\end{center}
The splitting condition of sequence (\ref{Fc}) is just
$f(g(H^0(\OO_C)))=0$ , where $f$ and $g$ are constructed in the
following diagram obtained by taking global sections in the sequence
(\ref{F}) twisted by $\theta (2)$.

\begin{center}
{\begin{picture}(2500,550)
\setsqparms[1`1`1`1;600`300]
\putsquare(1000,0)[H^0 F\theta(2)`H^0\OO_C`H^0
F\theta(2)`K\otimes H^0\theta(2); `\wr`g`]
\setsqparms[1`0`1`1;600`300]
\putsquare(1600,0)[\phantom{H^0\OO_C}`H^1(\OO_C(1))`\phantom{K\otimes H^0\theta(2)}`\Sigma;```f]
\putmorphism(1600,550)(0,-1)[0``]{250}1a
\putmorphism(2200,550)(0,-1)[0``]{250}1a
\putmorphism(0,300)(1,0)[0`\phantom{H^0(\OO_C(1))}`]{400}1a
\putmorphism(400,300)(1,0)[H^0(\OO_C(1))`\phantom{H^0 F\theta(2)}`]{600}1a
\end{picture}}
\end{center}
So, the condition that the 2-dimensional vector space $K\subset H^0(\theta
(2))$ gives an instanton, which is also
the vanishing of $H^0(\OO_C)$ in $\Sigma $, can be translated if $W$
denote the preimage by $c$ of
$g(H^0(\OO_C))$ in $H\dual \otimes V \otimes K$, by the condition:
$a(W)\subset {\rm Im}b$. 

\bigskip
We will now identify this condition with an explicit map
\hbox{$\stackrel{2}{\Lambda }H^0(\theta (2))
\stackrel{\beta}{\longrightarrow} H^1(\OO_C(1))$}. To construct $\beta
$, let's first consider the Eagon-Northcott complexes associated to
the sequence (\ref{reseau}) twisted by $-1$. This gives for the second
symmetric power the following sequence:
$$ 0\rightarrow (\stackrel{2}{\Lambda}H)(-2) \rightarrow H\dual 
\otimes H(-1) \kercoker{A\dual} S_2H\dual \rightarrow \theta ^2(2)
\rightarrow 0$$
We obtain by dualizing it the 2 following short exact sequences:
\begin{center}
$0 \rightarrow S_2H \otimes \OO_Y\rightarrow A \rightarrow
 \OO _C(1)
 \rightarrow 0$ and $0\rightarrow A \rightarrow H\dual (1) \otimes H
\stackrel{e}{\rightarrow} (\stackrel{2}{\Lambda}H)\dual (2)
\rightarrow 0$
\end{center}
Let $S_2H \otimes \OO_Y \stackrel{d}{\rightarrow}  H\dual (1) \otimes
H $ be the composition of the 2 previous injections. Then $d$ and $e$
are illustrated in the following diagram not exact in the middle:
\begin{center}
{\begin{picture}(2500,300)
\setsqparms[1`1`1`1;800`300]
\putsquare(500,0)[S_2 H\otimes\OO_Y`H\otimes H^{{\rm v}}(1)`S_2 \left(H^{\rm v}(1)\right)`H^{{\rm v}}(1)\otimes H^{{\rm v}}(1);d`S_2 m`m\otimes 1`]
\setsqparms[1`0`1`1;800`300]
\putsquare(1300,0)[\phantom{H\otimes H^{{\rm v}}(1)}`(\buildrel{2}\over{\Lambda}H^{{\rm v}})(2)`\phantom{H^{{\rm v}}(1)\otimes H^{{\rm v}}(1)}`(\buildrel{2}\over{\Lambda}H^{{\rm v}})(2);e``\wr`]
\putmorphism(0,0)(1,0)[0`\phantom{S_2 H^{\rm v}(1) \otimes\OO_Y}`]{500}1a
\putmorphism(0,300)(1,0)[0`\phantom{S_2 H\otimes\OO_Y}`]{500}1a
\putmorphism(2100,0)(1,0)[\phantom{(\buildrel{2}\over{\Lambda}H^{\rm v})(2)}`0`]{500}1a
\putmorphism(2100,300)(1,0)[\phantom{(\buildrel{2}\over{\Lambda}H^{{\rm v}})(2)}`0`]{500}1a
\end{picture}

}
\end{center}
The display of the monad $(d,e)$ made by the first line is: 

\begin{center}

\begin{picture}(2000,1100)
\setsqparms[1`1`1`1;600`300]
\putsquare(1100,250)[H\otimes H^{\rm
v}(1)`B`(\buildrel{2}\over{\Lambda}%
H^{\rm v})(2)`(\buildrel{2}\over{\Lambda}H^{\rm v})(2);`e``\sim ]
\setsqparms[1`1`1`1;600`300]
\putsquare(500,550)[S_2 H \otimes \OO_Y`A`S_2
H \otimes\OO_Y`\phantom{H\otimes H^{\rm
v}(1)};`\wr ``d]
\setsqparms[1`0`1`0;600`300]
\putsquare(1100,550)[\phantom{A}`\OO_C(1)`\phantom{H\otimes H^{\rm
v}(1)}`\phantom{B};```]
\putmorphism(0,550)(1,0)[0`\phantom{S_2 H \otimes\OO_Y}`]{500}1a
\putmorphism(0,850)(1,0)[0`\phantom{S_2 H \otimes\OO_Y}`]{500}1a
\putmorphism(1700,550)(1,0)[\phantom{B}`0`]{300}1a
\putmorphism(1700,850)(1,0)[\phantom{\OO_C(1)}`0`]{300}1a
\putmorphism(1100,250)(0,1)[\phantom{(\buildrel{2}\over{\Lambda}H^{\rm v})(2)}`0`]{250}{1}l
\putmorphism(1700,250)(0,1)[\phantom{(\buildrel{2}\over{\Lambda}H^{\rm v})(2)}`0`]{250}{1}l
\putmorphism(1100,1100)(0,1)[0`\phantom{A}`]{250}1l
\putmorphism(1700,1100)(0,1)[0`\phantom{\OO_C(1)}`]{250}1l
\end{picture}

\end{center}
The last column of this display gives when passing to global sections
a boundary map \hbox{$\stackrel{2}{\Lambda }H\dual\otimes S_2V
\stackrel{\beta_0}{\rightarrow} H^1(\OO_C(1))$} which vanishes on the
image of the map \hbox{$H\otimes H\dual \otimes V \stackrel{b^{\prime
    }}{\rightarrow } \stackrel{2}{\Lambda }H\dual \otimes S_2V$}
obtained from $e$ by taking global sections. Using now the surjection
of (\ref{reseau}) we obtain the commutative diagram:
 \begin{center}
\begin{picture}(800,300)
\setsqparms[1`1`1`1;800`300]
\putsquare(0,0)[H\otimes H^{\rm v}\otimes V`H^{\rm v}\otimes H^{\rm v}\otimes
S_2 V`H\otimes H^0 \theta (2)`H^{\rm v}\otimes H^0 \theta (3);b^{\prime}```b]
\end{picture}
\end{center}
which identify the conditions $a^{\prime }(W)\subset Im b^{\prime}$
and $a(W) \subset Im b$, where $a^{\prime }$ is defined in the diagram
below. As we had already identified the splitting condition with $a(W)
\subset Im b$, we can from the following diagram construct a map
$\beta$ giving the corollary \ref{scindage}.
\begin{center}
{\begin{picture}(2000,550)
\setsqparms[1`-1`-1`1;800`300]
\putsquare(0,0)[H\otimes H^{{\rm v}}\otimes V`\buildrel{2}\over{\Lambda}H^{{\rm v}}\otimes S_2 V`H\otimes H^{{\rm v}}\otimes V`\buildrel{2}\over{\Lambda}(H^{{\rm v}}\otimes
V);b^{\prime }`\wr`a^{\prime }`]
\setsqparms[1`0`-1`1;800`300]
\putsquare(800,0)[\phantom{\buildrel{2}\over{\Lambda}H^{{\rm v}}\otimes S_2 V}`H^1 \OO_C (1)`\phantom{\buildrel{2}\over{\Lambda}(H^{\rm v}\otimes
V)}`\buildrel{2}\over{\Lambda}H^0\theta(2);\beta _0``\beta`c]
\put(800,350){\vector(0,1){150}}
\put(775,525){0}
\putmorphism(1600,0)(1,0)[\phantom{\buildrel{2}\over{\Lambda}H^0\theta(2)}`0`]{500}1a
\putmorphism(1600,300)(1,0)[\phantom{H^1 \OO_C(1)}`0`]{500}1a
\end{picture}
}
\end{center}
\begin{corollary}
\label{scindage}The splitting condition of the sequence (\ref{Fc}) is
given by the vanishing of a surjective map $\stackrel{2}{%
\Lambda }(H^0\theta (2))\stackrel{\beta }{\rightarrow }H^1(\OO_C(1))$
\end{corollary}
This gives the following description of $Z_{\delta_{\theta}}$ where we
    denote by $Z_s$ the vanishing locus of a 
    section $s$ of $\theta (2)$:
$$Z_{\delta_{\theta}}=\proj(\{s\wedge s^{\prime }|s,s^{\prime }\in
H^0\theta (2),\;Z_s \cap Z_{s^{\prime}}= \emptyset \;\text{and}\;\beta (s\wedge s^{\prime })=0\})$$
\begin{itemize}
\item Second description
\end{itemize}
Let $E$ be an $n$-instanton, $\theta $ its associated $\theta
$-characteristic, and $s$, $s^{\prime}$ be 2 sections of $\theta (2)$
associated to $E$. The purpose of this second description is to
understand the possible singularities of
$Z_{\delta_{\theta}}$ at $s\wedge s^{\prime}$ in terms of the vector
bundle $E$. This part will be divided in 3 steps. First, we will
construct a new map $\beta_E $, then we'll have to identify the
vanishing locus of $\beta$ and $\beta_E$, and in the third step, this will
enable us to understand the particular bundles that give singularities of
$Z_{\delta_{\theta}}$.

\bigskip
1) We will construct here a map \hbox{$\beta_E:\ H^0\theta (2) \oplus H^0\theta
(2) \rightarrow H^1 \OO_C(1)$}. Let's take the exact sequence of
restriction to the exceptional divisor $x=\tau-\sigma$ twisted by
$E(\tau)$, where $\tau=p^* \Op(1)$ and $\sigma= q^* \OO_Y(1)$.
\begin{center}
$0\longrightarrow p^{*}E(\sigma )\longrightarrow p^{*}E(\tau
)\longrightarrow p^{*}E(\tau )_{\mid x}\longrightarrow 0$
\end{center}
apply it to the functor $q_*$ to get the following exact sequence where
we recall that $F=q_*p^*E$, and that the 2 dimensional vector space in
the right of the sequence is naturally the fiber of $E$ at $N$.
\begin{equation}
  \label{Etau}
  0\longrightarrow F(1)\longrightarrow q_{*}p^{*}E(\tau)
  \longrightarrow  2 \OO_Y \longrightarrow 0
\end{equation}
Let $2H^0\theta (2)\stackrel{\delta_E}{\rightarrow} H^1F\theta (3)$ be
the boundary map obtained when taking global sections in (\ref{Etau})
twisted by $\theta (2)$. Taking global sections in (\ref{Fc}) twisted
by $\theta (3)$ gives us a map $H^1F\theta (3)\rightarrow
H^1\OO_C(1)$. Define $\beta _E$ as the composition of $\delta_E$ and
this map:
$$ \beta_E:\ 2H^0\theta (2) \stackrel{\delta_E}{\longrightarrow}
H^1F\theta (3) \longrightarrow H^1(\OO_C(1))$$

\bigskip
2) To link $\beta$ and $\beta _E$, we will prove here that for every
surjection \hbox{$2\OO_Y \stackrel{(\sigma ,\sigma ^{\prime})}{
    \longrightarrow } \theta (2)\rightarrow 0$} where  $\sigma ,\sigma
^{\prime }$ are such that $\beta(\sigma\wedge \sigma ^{\prime })=0$,
then $\beta _E(-\sigma^{\prime },\sigma)=0$. 

For all that, apply to the sequence \ref{Etau} the functor
$Hom(*,\theta (2))$. It gives a map $Hom(2\OO_Y,\theta (2))
\rightarrow Hom(q_{*}p^{*}E(\tau ),\theta (2))$, which enable to
pull back any $\phi: 2\OO_Y \stackrel{(\sigma ,\sigma ^{\prime})}{
    \longrightarrow } \theta (2)\rightarrow 0$ by a map
  $q_*p^*E(\tau)\stackrel{\psi}{\rightarrow} \theta(2)$ making the
  following diagram commutative:\\
\vbox{\begin{equation}\label{phi}
{\begin{picture}(1600,550)(0,550)
\setsqparms[1`-1`0`1;500`300]
\putsquare(300,250)[F (1)`\phantom{q_* p^* E(\tau )}`F(1)`\phantom{N};`\wr ``]
\setsqparms[1`-1`-1`1;500`300]
\putsquare(800,250)[q_* p^* E(\tau)`2\OO_{\proj_2}`N`F^{\prime };```]
\setsqparms[1`-1`-1`0;500`300]
\putsquare(800,550)[\theta (2)`\theta (2)`\phantom{mmmm}`\phantom{mmm};\sim `\psi`\phi`]
\setsqparms[0`-1`-1`0;500`250]
\putsquare(800,850)[0`0`\phantom{mm}`\phantom{mm};```]
\putmorphism(0,250)(1,0)[0`\phantom{F(1)}`]{300}1a
\putmorphism(0,550)(1,0)[0`\phantom{F(1)}`]{300}1a
\putmorphism(1300,850)(1,0)[\phantom{\theta (2)}`0`]{300}1a
\putmorphism(1300,250)(1,0)[\phantom{F^{\prime }}`0`]{300}1a
\putmorphism(1300,550)(1,0)[\phantom{2\OO_C}`0`]{300}1a
\setsqparms[0`-1`-1`0;500`250]
\putsquare(800,0)[``0`0;```]
\end{picture}
}
\end{equation}
\vskip 5.5em}
where $N$ and $F^{\prime}$ are by definition the kernel of $\psi $ and
$\phi$, and where the middle line is the sequence \ref{Etau}.
Restricting (\ref{phi}) to $C$ gives the 2 following diagrams where
the kernel of the restriction of $\psi$ to $C$ is noted $M$:\\
\vbox{
\begin{equation}\label{phiC1}
{\begin{picture}(1600,550)(0,550)
\setsqparms[1`-1`0`1;500`300]
\putsquare(300,250)[F_C (1)`\phantom{q_* p^* E(\tau )_C}`F_C (1)`\phantom{M};`\wr ``]
\setsqparms[1`-1`-1`1;500`300]
\putsquare(800,250)[q_* p^* E(\tau)_C`2\OO_{\proj_2}`M`\theta
^{\rm v}(-2);```]
\setsqparms[1`-1`-1`0;500`300]
\putsquare(800,550)[\theta (2)`\theta (2)`\phantom{mmmm}`\phantom{mmm};\sim `\psi`\phi`]
\setsqparms[0`-1`-1`0;500`250]
\putsquare(800,850)[0`0`\phantom{mm}`\phantom{mm};```]
\putmorphism(0,250)(1,0)[0`\phantom{F_C(1)}`]{300}1a
\putmorphism(0,550)(1,0)[0`\phantom{F_C(1)}`]{300}1a
\putmorphism(1300,850)(1,0)[\phantom{\theta (2)}`0`]{300}1a
\putmorphism(1300,250)(1,0)[\phantom{\theta ^{\rm v}(-2)}`0`]{300}1a
\putmorphism(1300,550)(1,0)[\phantom{2\OO_C}`0`]{300}1a
\setsqparms[0`-1`-1`0;500`250]
\putsquare(800,0)[``0`0;```]
\end{picture}
}
\end{equation}
\vskip 6em}\\
\vbox{
\begin{equation}\label{phiC2}
{\begin{picture}(1600,550)(0,550)
\setsqparms[1`-1`-1`1;500`300]
\putsquare(300,550)[F_C (1)`M`F_C (1)`N_C;`\wr ``]
\setsqparms[1`0`-1`1;500`300]
\putsquare(800,550)[\phantom{M}`\theta ^{\rm v}(-2)`\phantom{N_C}`F^{\prime }_C;```]
\setsqparms[0`-1`-1`1;500`300]
\putsquare(800,250)[\phantom{N_C}`\phantom{F^{\prime }_C}`\theta^{\rm v}(-1)`\theta^{\rm v}(-1);```\sim]
\setsqparms[0`-1`-1`0;500`250]
\putsquare(800,850)[0`0`\phantom{M}`\phantom{\theta ^{\rm v}(-2)};```]
\putmorphism(0,550)(1,0)[0`\phantom{F_C(1)}`]{300}1a
\putmorphism(0,850)(1,0)[0`\phantom{F_C(1)}`]{300}{1}a
\putmorphism(1300,850)(1,0)[\phantom{\theta ^{\rm v}(-2)}`0`]{300}1a
\putmorphism(1300,550)(1,0)[\phantom{F^{\prime }_C}`0`]{300}1a
\setsqparms[0`-1`-1`0;500`250]
\putsquare(800,0)[``0`0;```]
\end{picture}
}
\end{equation}
\vskip 5.5em}
Construct a map $0 \rightarrow \theta\dual \stackrel{i}{\rightarrow}M$
by composition of the injection $0 \rightarrow \theta\dual \rightarrow
F_C(1)$ of sequence (\ref{Fc}) with the injection of $F_C(1)$ in $M$
of the diagram (\ref{phiC2}). Denote by $M^{\prime}$ the cokernel of
$i$, and by $E^{\prime}$ the cokernel of the injection $0 \rightarrow
\theta\dual \rightarrow q_*p^*E(\tau)_C$ obtained by composing $i$
with the injection of $M$ in $q_*p^*E(\tau)_C$ of the diagram
(\ref{phiC1}). Now using the fact that we have an injection of the
cokernel of (\ref{Fc}) in 
$M^{\prime}$ and in $E^{\prime}$, we obtain from (\ref{phiC1}) the
following diagram:\\ 
\vbox{
\begin{equation}\label{Mprime}
{\begin{picture}(1600,550)(0,550)
\setsqparms[1`-1`0`1;500`300]
\putsquare(300,250)[\theta ^{\rm v}(-1)`E^{\prime }`\theta ^{\rm v}(-1)`M^{\prime };`\wr ``]
\setsqparms[1`-1`-1`1;500`300]
\putsquare(800,250)[\phantom{E^{\prime }}`2\OO_{\proj_2}`\phantom{M^{\prime }}`\theta ^{\rm v}(-2);```]
\setsqparms[1`-1`-1`0;500`300]
\putsquare(800,550)[\theta (2)`\theta (2)`\phantom{E^{\prime }}`\phantom{2\OO_{\proj_2}};\sim ``\phi`]
\setsqparms[0`-1`-1`0;500`250]
\putsquare(800,850)[0`0`\phantom{mm}`\phantom{mm};```]
\putmorphism(0,250)(1,0)[0`\phantom{\theta ^{\rm v}(-1)}`]{300}1a
\putmorphism(0,550)(1,0)[0`\phantom{\theta ^{\rm v}(-1)}`]{300}1a
\putmorphism(1300,250)(1,0)[\phantom{\theta ^{\rm v}(-2)}`0`]{300}1a
\putmorphism(1300,550)(1,0)[\phantom{2\OO_C}`0`]{300}1a
\setsqparms[0`-1`-1`0;500`250]
\putsquare(800,0)[``0`0;```]
\end{picture}
}
\end{equation}  
\vskip 5.5em}
We can now notice that $\beta _E$ is just the boundary obtained when
taking global section in the middle line of (\ref{Mprime}) twisted by
$\theta(2)$. As the right column of (\ref{Mprime}) twisted by $\theta
(2)$ gives an injection on $H^0(\OO_C)$ in $2H^0(\theta(2))$ of image
$(-\sigma^{\prime}, \sigma)$, we obtain that $\beta_E
(-\sigma^{\prime}, \sigma)=0$ if and only if $M^{\prime}$ splits. But
$F_C$ split because $\beta(\sigma^{\prime}\wedge \sigma)=0$ by
hypothesis, so we have after the choice of a section $\alpha$ of
(\ref{Fc}) the diagram:
\begin{center}
{\begin{picture}(1550,1100)
\setsqparms[1`-1`-1`1;500`300]
\putsquare(300,550)[\theta ^{\rm v}`M`F_C (1)`N_C;i```]
\setsqparms[1`-1`-1`1;450`300]
\putsquare(800,550)[\phantom{M}`M^{\prime }`\phantom{N_C}`F^{\prime }_C;```]
\setsqparms[0`-1`-1`1;500`300]
\putsquare(300,250)[\phantom{F_C (1)}`\phantom{N_C}`\theta ^{\rm v}(-1)`\theta ^{\rm
v}(-1);`\alpha ``\sim]
\setsqparms[0`-1`-1`0;500`250]
\putsquare(300,850)[0`0`\phantom{mm}`\phantom{mm};```]
\putmorphism(0,550)(1,0)[0`\phantom{F_C (1)}`]{300}1a
\putmorphism(0,850)(1,0)[0`\phantom{\theta ^{\rm v}}`]{300}1a
\putmorphism(1250,550)(1,0)[\phantom{F^{\prime }_C}`0`]{300}1a
\putmorphism(1250,850)(1,0)[\phantom{M^{\prime }}`0`]{300}1a
\setsqparms[0`-1`-1`0;500`250]
\putsquare(300,0)[``0`0;```]
\end{picture}
}
\end{center}
which proves that $M^{\prime }\simeq F_C^{\prime }$, and that $\beta$
and $\beta_E$ have the same vanishing locus. \fin

\bigskip
3) Identification of $Z_{\delta_{\theta}}$ 's singularities.\\
Let $E$ be an instanton, and $s\wedge s^{\prime }$ its associated
element of $Z_{\delta_{\theta}}$. Denote by $K$ the vector space
$Vect(s,s^{\prime})\subset H^0(\theta(2))$.

The bundle $E$ gives a singularity of $G(2,H^0(\theta(2)))\cap
\ker\beta$ if and only if the Zariski tangent space to
$G(2,H^0(\theta(2)))$ at $s\wedge s^{\prime }$ is included in one of
the hyperplanes given by $\beta $. In other words, it means that $K$
is in the kernel of one of the skew map given by $\beta $. As
$H^1(\OO_C(1))=S_{n-4}V\dual$, it means  that there is some $p\in
S_{n-4}V$ such that $K\otimes p$ is in the kernel of \hbox{$K \otimes
  S_{n-4}V \stackrel{\beta_E}{\rightarrow} H^0(\theta(2))\dual $}.

But we can understand this kernel explicitly with the following
diagram, where the first 2 lines are obtained from the sequences
(\ref{F}) and (\ref{Fc}), and the last from sequence (\ref{Etau})
twisted by $\OO_Y(n-4)$:\\
\vbox{
\begin{equation}\label{diagsing}
{\begin{picture}(3620,425)(0,425)
{\footnotesize
\setsqparms[1`-1`-1`1;600`300]
\putsquare(2700,550)[H^2 F(-3)`H^2 (2\OO_C (-3))`H^1(F_C (n-3))`(H^0 \theta(2))^{\rm v};```]
\setsqparms[1`-1`0`1;700`300]
\putsquare(2000,550)[H^1 \theta(-1)`\phantom{H^2 F(-3)}`H^1 \theta^{\rm v}(n-4)`\phantom{H^1(F_C (n-3))};`\wr ``]
\putmorphism(-200,250)(1,0)[0`H^0 (q_* p^* E(n-3)\sigma)`]{500}1a
\putmorphism(300,250)(1,0)[\phantom{H^0 (q_* p^* E(n-3)\sigma)}`H^0 (q_* p^* E(\tau +(n-4)\sigma))`]{950}1a
\putmorphism(1250,250)(1,0)[\phantom{H^0 (q_* p^* E(\tau +(n-4)\sigma))}`K \otimes S_{n-4}V`]{800}1a
\putmorphism(2050,250)(1,0)[\phantom{K \otimes S_{n-4}V}`H^1 (F(n-3))`\delta_E]{650}1a
\putmorphism(2700,250)(0,1)[\phantom{MM}`0`]{250}{-1}l
\putmorphism(2700,550)(0,1)[\phantom{MM}`\phantom{MM}`]{300}{-1}l
\putmorphism(1500,850)(1,0)[0`\phantom{H^1 \theta(-1)}`]{500}{1}a
\putmorphism(3300,550)(1,0)[\phantom{(H^0 \theta(2))^{\rm v}}`0`]{320}{1}a
\putmorphism(1500,550)(1,0)[0`\phantom{H^1 \theta^{\rm v}(n-4)}`]{500}{1}a
\bezier{0}(2050,300)(2050,400)(2300,400)
\put(2300,400){\line(1,0){750}}
\bezier{0}(3050,400)(3300,400)(3300,500)
\put(3300,500){\vector(0,1){20}}
\put(3150,370){$\beta_E$}
}
\end{picture}
}
\end{equation}
\vskip 3.0em}
This proves that the maps $\beta_E$ and $\delta _E$ have the same
kernel, which is nothing else than the cokernel of $H^0(q_{*}p^{*}E(n-3)\sigma
)\stackrel{j}{\rightarrow} H^0(q_{*}p^{*}E(\tau +(n-4)\sigma ))$.
\begin{corollary}\label{ramification}
  The ramification locus of the map $U_n^s \rightarrow \Theta^s_{lf}$
  is made of special t'Hooft bundle for $c_2=4$.\\
  For $c_2=5$, there are 3 possibilities. If the Zariski tangent space
  of $Z_{\delta_{\theta}}$ is 15 (respectively 16) dimensional, then
  this point arise from a t'Hooft bundle (respectively special
  t'Hooft). If it is 14 dimensional, then the bundle twisted by 2 have
  2 sections vanishing at $N$.
\end{corollary}
{\bf NB:} For $c_2=5$, a bundle such that $h^0\ideal_N \otimes
E(2)\geq 2$ doesn't necessarily give a singularity of $Z_{\delta_{\theta}}$.

The result \ref{ramification} is clear for $c_2=4$ because we have in this situation
$h^0(q_*p^*E(\tau)) \geq 2$.

When $c_2=5$, any singularity of $Z_{\delta_{\theta}}$ is such that $K
\otimes p$ is in the cokernel of $j$ for some $p\in V$, so
$h^0(q_*p^*E(\tau +\sigma)) \geq 2$. Although this will be enough to
prove the smoothness of I$_5$, we will go further to understand the
ramification, and also because it seems it is not enough to obtain the
connexity. The cases of a 15 or 
16 dimensional Zariski tangent space, have the following interpretation in term of the
cokernel of $H^0q_*p^*E(\sigma) \stackrel{j^{\prime}}{\rightarrow} H^0q_*p^*E(\tau)$:

Take global sections in (\ref{Etau}) to obtain  like in the
diagram (\ref{diagsing}) a
map $\beta ^{\prime }:K \rightarrow H^1\theta\dual(-1) = H^0(\theta
(3))\dual $ whose kernel is the cokernel of $j^{\prime}$. In fact
$\beta$ is the composition of $\beta^{\prime}$ with the injection
\hbox{$H^0(\theta(3))\dual \rightarrow (H^0(\theta(2))\otimes
  V)\dual$}. But the Zariski tangent space of $Z_{\delta_{\theta}}$ at
$s \wedge s^{\prime}$ is
15 (resp 16) dimensional if and only if there is a 2 (resp 3)
dimensional subspace $W$ of $V$ such that the map \mbox{$K\otimes H^0(\theta
(2))\otimes W\rightarrow \complex $} induced by $\beta $ is zero.

-If dim$W=3$, then $\beta^{\prime}=0$ and thus $h^0(q_*p^*E(\tau))=2$.

-If dim$W=2$, then we have to find an element $s_0 \in K$ such that
the map $\beta^{\prime}(s_0): H^0(\theta(3)) \rightarrow \complex$ is
zero.
Let $(p,q)$ be a basis of $W$, and compute the dimension of
\hbox{$\{H^0(\theta (2))\otimes p\}\cap \{H^0(\theta (2))\otimes q\}$}
in $H^0(\theta(3))$. Indeed, as \hbox{$H^0\theta(3)\otimes K
\stackrel{\beta^{\prime}}{\rightarrow} \complex $} is zero on the image
of $H^0(\theta(2)) \otimes W \otimes K$ in $H^0(\theta(3)) \otimes K$,
if $\{H^0(\theta (2))\otimes p\}+\{H^0(\theta (2))\otimes q\}$ is 0 or
1 codimensional in $H^0\theta(3)$, then there is such an $s_0$. So
assume that it is not the case, then $\dim \{H^0(\theta (2))\otimes
p\}\cap \{H^0(\theta (2))\otimes q\}\geq 7$, and this intersection would
contain 2 elements independent relatively to $H^0(\theta(1))\otimes
pq$. So there would exist 2 sections $s,\sigma$ of $\theta (2)$
(independent relatively to $H^0(\theta(1))\otimes q$), and
$s^{\prime},\sigma^{\prime}$ (independent relatively to
$H^0(\theta(1))\otimes p$) such that $s.p=s^{\prime}.q$ and
$\sigma.p=\sigma^{\prime}.q$. As those sections are not in
$H^0(\theta(1))\otimes pq$, the lines $p,q$ have to cut one another in a point
$P\in C$, and the 4 remaining points of $q\cap C$ have to be zeros of
$s$ and $\sigma$. As $\theta(2)$ is generated by its global sections,
those vanishing in a given point are always an hyperplane of
$H^0(\theta(2))$. As the intersection of the 5 hyperplanes associated
to $q\cap C$ is $H^0(\theta(1))\otimes q$, the intersection of 4 of
them can't be 7 dimensional, which contradicts the hypothesis.   
\begin{remark}\label{ngone}
  For $c_2=n$, if there is a section of $\theta (2)$ which is in the cokernel of
  $H^0q_*p^*E(\sigma) \stackrel{j^{\prime}}{\rightarrow}
  H^0q_*p^*E(\tau)$, then this section vanishes on the vertices of the
  complete $(n+1)$-gone (inscribed in the curve of jumping lines)
  obtained from the lines through $N$ which are bisecant to the zero
  locus of the associated section of $E(1)$  
\end{remark}
Let $E$ be an $n$-instanton such that there is a section $s$ of
$\theta (2)$ in the cokernel of $j^{\prime}$. Let $\Im $ be the ideal
of the vanishing locus of the section of $E(1)$ coming from $s$. The
bundle $F=q_*p^*E$ is also $q_*p^*\Im(\tau)$, and we can understand
$s$ with the following diagram, where the first column is the natural
evaluation: \hbox{$ (q_{*}p^{*}E) \otimes (q_{*} \OO_{
    \widetilde{\proj}\!_3} (\tau))\rightarrow q_{*} p^{*}E(\tau )$}.  

\begin{center}
{\begin{picture}(2200,900)
\setsqparms[1`-1`-1`1;900`300]
\putsquare(800,0)[q_* p^* E(\tau )`q_{*}p^{*}\Im (2\tau )`F\oplus F(1)`q_{*}p^{*}\Im (\tau )\otimes q_{*}(\OO_{\widetilde{\proj}\!_3}(\tau ));```\sim ]
\setsqparms[0`-1`-1`1;900`300]
\putsquare(800,600)[0`0`\theta (2)`L;```]
\put(800,375){\vector(0,1){150}}
\put(1700,375){\vector(0,1){150}}
\putmorphism(300,300)(1,0)[\OO_Y`\phantom{q_* p^* E(\tau )}`]{500}1a
\putmorphism(1700,300)(1,0)[\phantom{q_{*}p^{*}\Im (2\tau )}`0`]{500}1a
\putmorphism(0,300)(1,0)[0`\phantom{\OO_Y}`]{300}1a
\put(400,400){\vector(2,1){250}}
\put(500,500){$s$}
\end{picture}}
\end{center}
As $R^1q_*p^*\Im(2\tau) =0$ because there are no multi-jumping lines
through $N$, we have for any line $D$ containing $N$ a
surjection from $q_*p^*\Im (2\tau)_{\{d\}}$ onto $H^0(\Im_D(2\tau))$
where $\Im_D=\Im \otimes \OO_D$. But the map \hbox{$q_{*}p^{*}\Im
  (\tau )_{\{d\}}\otimes q_{*}(\OO_{\widetilde{\proj}\!_3}(\tau ))
  \rightarrow H^0(\Im _D(2\tau)) $} has to be zero if $D$ is bisecant
to the scheme defined by $\Im$, thus the support of $L$ must contain
the points corresponding to those lines. Computing the degree, we can
conclude that $s$ vanish on the $(n+1)$-gone whose vertices are the
projections from $N$ of the previous bisecant.
\begin{corollary}
  When $c_2=4$ the fiber $Z_{\delta_{\theta}}$ is singular if and only
  if the support of $\theta $ is a Lur\"oth quartic. In this situation,
  the singularity $s\wedge s^{\prime}$ is such that the zeros of $s$ and
  $s^{\prime}$ give the pencil of complete pentagons inscribed on the
  quartic.
\end{corollary}

\subsection{The normality condition when $c_2 =5$; applications}
We want here to understand when $Z_{\delta_{\theta}}$ is not regular
in codimension 1. When $c_2=5$, the fiber $Z_{\delta_{\theta}}$ is an
open subset of the Grassmannian $G(2,10)$ cut by 3 hyperplanes. Let
$A_i$ be the skew forms associated to those hyperplanes, and $a_i$ be their
skew linear maps. The key result of \ref{ramification} and \ref{ngone}
is that any $s$ in the kernel of all the $a_i$ gives a $(n+1)$-gone
inscribed in $\theta$'s support. So we want here to prove the
following:

\begin{proposition}
  \label{fibrenonorm}When $c_2=5$, if $G(2,H^0\theta(2)) \cap \ker
  \beta$ is singular in codimension 1 or have an excess  dimension, then
  the $a_i$ have at least a 4-dimensional common kernel 
\end{proposition}

Both singularity in codimension 1 and excess dimension imply that
there is a 12 dimensional subscheme $S$ of $G(2,H^0\theta(2)) \cap \ker\beta$
such that the intersection is not transverse at any points of $S$. So,
any
$s\wedge s^{\prime }$ of $S$ is necessarily in some $\ker
(\sum_{i=1}^{i=3}\lambda_i a_i)$, then $S$ is in
$\bigcup\limits_{\lambda _i}G(2,\ker (\sum_{i=1}^{i=3}\lambda _ia_i))$
. As $S$ is 12 dimensional, at least one of those Grassmann manifold is
of dimension 10 or more. So this Grassmann manifold has to be at least
12 dimensional, and for this value of $(\lambda _1,\lambda _2,\lambda
_3)$ one have $\dim \ker(\sum_{i=1}^{i=3}\lambda _ia_i)=8$ (none of
those maps are zero according to the \ref{scindage}). We can
assume that this occur for $(1,0,0)$, so $\dim \ker a_1=8$. This
yields to the following discussion:

a) There is no $(\lambda _1,\lambda _2,\lambda _3) \neq (1,0,0)$ such
that $\dim \ker \sum_{i=1}^{i=3}\lambda _ia_i)=8$. Then the union of
those Grassmann manifolds is the union of $G(2,\ker a_1)$ with
something of dimension at most 8+2. So we must have $G(2,\ker
a_1)=S$, but it implies that all the $A_i$ are zero on
$\stackrel{2}{\wedge }\ker a_1$. The $a_i$ have then at least a 4
dimensional common kernel.

b) Otherwise, we can assume that $\dim \ker a_2=8$. If $a_3$ had also
an 8 dimensional kernel, then the $a_i$ would easily have a 4
dimensional common kernel, so let's assume that if $\lambda_3 \neq 0$
then $\dim \ker (\sum_{i=1}^{i=3}\lambda _ia_i) \leq 6$. So
\hbox{$\dim \bigcup\limits_{\lambda _i,\lambda_3\neq 0}G(2,\ker
(\sum_{i=1}^{i=3}\lambda _ia_i))\leq 10$}, then $S\subset
\bigcup\limits_{\lambda _i}G(2,\ker(\sum_{i=1}^{i=2}\lambda _ia_i))$.
So $S\cap G(2,\ker a_1)$ is at most 1 codimensional in $G(2,\ker
a_1)$, and it is in the vanishing of $A_{2_{|\stackrel{2}{\wedge }\ker a_1}}$
and of $A_{3_{|\stackrel{2}{\wedge }\ker a_1}}$. So those 2 skew forms
have to be proportional on $\ker a_1$, then we can assume that
$A_{3_{|\stackrel{2}{\wedge }\ker a_1}}=0$ so the $a_i$ have again a 4
dimensional common kernel as claimed.   

\begin{theorem}
  \label{irred5}
  The moduli space I$_5$ is irreducible of dimension 37.
\end{theorem}
Let's first show that for a general $\theta$ in the image of $U_N^s$,
then  $Z_{\delta_{\theta}}$ is irreducible. As $Z_{\delta_{\theta}}$
is open in $G^{\prime}_{\theta}\cap \ker \beta$,it is enough to show
that $G^{\prime}_{\theta} \cap \ker \beta$ is irreducible, where
$G^{\prime}_{\theta}=G(2,H^0(\theta(2)))$. But any
degeneracy locus of $3\OO_{G^{\prime}_{\theta}} \rightarrow
\OO_{G^{\prime}_{\theta}}(1)$ is connected (Cf [A-C-G-H] p311), and complete
intersection. So it satisfies Serre's condition (S2), and in our
situation, the \ref{fibrenonorm} shows that if it is not regular in
codimension 1 (R1) or if it is excess dimensional, then the support of
$\theta$ is a Darboux pentic according to \ref{ngone} (in fact it
would be 4 times Darboux if this had a sense!). Anyway $\theta$ can't
be general in $\Theta$. Then if $\theta$ is not Darboux,
$G^{\prime}_{\theta}\cap \ker \beta$ is 
normal, connected and 13-dimensional, so it is irreducible and
$Z_{\delta_{\theta}}$ too.

But we proved in the \ref{recstable} that the $(U_N^s)_{N \in
  \proj_3}$ cover $I_n$.  Furthermore, the
basis of the fibration $U_N^s\rightarrow \Theta ^s_{lf}$ is
irreducible and smooth according to [So], and the fiber is $\proj_4
\times Z_{\delta_{\theta}}$ which is irreducible, normal and 17 dimensional
when $\theta$ is not Darboux. Furthermore, for any
$\theta$ in the image of $U_N^s$, $Z_{\delta_{\theta}}$ is at most 15
dimensional because 
$G^{\prime}_{\theta}$ is not include in a hyperplane. As $U_N$ is
open in I$_5$, all its irreducible components have dimension at least
37. But the Darboux
$\theta$ are 3 codimensional in $\Theta^s_{lf}$ which is 20
dimensional, so the preimage of the Darboux $\theta$ is too small to
make an irreducible component of $U_N^s$. Then $U_N^s$ has to be
irreducible, and 37 dimensional. 
We can conclude using the \ref{recstable} that it is also the case for
I$_5$.
\begin{theorem}
  \label{lissite}
  The moduli space of mathematical instanton with $c_2=5$ is smooth
\end{theorem}
As $U_N^s$ is a fibration  over $\Theta^s_{lf}$ which is smooth, any
bundle which is not in the ramification of this morphism is a smooth
point of I$_5$. So we have to check that a bundle $E$ which is in all
the ramifications of $U_N^s \rightarrow \Theta^s_{lf}$ when $N$ fills
$\proj_3$ is a smooth point of I$_5$. According to the
\ref{ramification} we must have $h^0(\ideal_NE(2))\geq 2$ for every
$N$. So $h^0E(2) \geq 4$, and either $E$ is a t'Hooft bundle, or $E$ belongs to the
family described in the \ref{famille}, and both are smooth points of
the moduli space (Cf \ref{lisse}).  
\subsection{Residual class}\label{appendiceresiduel}

Let $E$ be an $n$-instanton, and assume that $E$ has a 2 parameter
family $S$ of multi-jumping lines. Furthermore, we assume here that $S$
is irreducible, and that the residual
scheme of $S$ in the scheme of multi-jumping lines $M$, is a curve (may
be empty) not drawn on $S$. Eventually assume that $S_{red}$ is locally
complete intersection.\\
\begin{picture}(800,0)(-60,800)
\setsqparms[1`1`1`1;400`400]
\putsquare(0,0)[\widetilde{S_{red}}`\widetilde{G}`S_{red}`G;j`g`f`i]
\setsqparms[0`0`1`0;400`400]
\putsquare(0,400)[`{\proj(H^1E\dual\otimes 
\OO_{\widetilde{G}})}`\phantom{\widetilde{G}}`;``\pi`]
\end{picture}
\hfill
\begin{minipage}[t]{10.5cm}

 In order to define and compute the class of $C$ in the chow ring
of $G$ (graduated by the dimension), we will have to work in the blow
up $\widetilde{G}$ of $G$ along $S_{red}$. Let $\widetilde{S_{red}}$
be the exceptional divisor, and $x$ be the class of
$\widetilde{S_{red}}$ in $A_3\widetilde{G}$, we have: \hfill
{$A_k\widetilde{G}=(A_k\widetilde{S_{red}}\oplus 
A_kG)/\alpha (A_{k}S_{red})$}\hfill\phantom{.}\\
where $\forall y\in A_{k}S_{red}$, $\alpha 
(y)=(c_1(g^{*}N/{\cal O}_N(-1)\cap g^{*}y,-i_{*}y)$, and where $N$ is
the normal bundle of $S_{red}$ in $G$. Let's define:\hfill
$\zeta=c_1(j^*\OO_{\widetilde{G}} (\widetilde{S_{red}})
=c_1({\OO}_N(-1))$ \hfill \phantom{.}
\end{minipage} 
The multiplicative structure of $A_*\widetilde{G}$ is given by the
following rules:\\
\begin{center}
$\left\{ 
\begin{array}{lcl}
f^{*}\gamma .f^{*}\gamma ^{\prime } & = & f^{*}\gamma \gamma ^{\prime } \\ 
j_{*}\widetilde{\sigma }.j_{*}\widetilde{\sigma ^{\prime }} & = &j_{*}(\zeta 
\widetilde{\sigma }\widetilde{\sigma ^{\prime }}) \\ 
f^{*}\gamma .j_{*}\widetilde{\sigma } & = & j_{*}((g^{*}i^{*}\gamma ).\widetilde{%
\sigma })
\end{array}
\right. $  
\end{center}
where $\gamma \in A_{*}G$ and $\widetilde{\sigma },\widetilde{\sigma
^{\prime }}\in A_{*}\widetilde{S_{red}}$. So we have:
$x^2=j_{*}\zeta =(g^{*}c_1N,-i_{*}[S_{red}])$, because $\zeta =-c_1(g^{*}N/%
{\cal O}\!_N(-1))+g^{*}c_1N$ and similarly 
$x^3=j_{*}\zeta ^2=((c_1N)^2-c_2(N),-i_{*}(c_1(N))$, and
\begin{center}
$A_{*}\widetilde{S_{red}}=A_{*}S_{red}[\zeta ]/(\zeta ^2-c_1(N)\zeta
+c_2(N)) $.
\end{center}
One has the exact sequence: $H^1E(-1)\otimes Q\dual \stackrel{\sigma
}{ \longrightarrow }H^1E\otimes {\OO}_G \rightarrow
R^1q_{*}p^{*}E\rightarrow 0$, where $Q$ is the tautological quotient
bundle of $H^0(\OO_{\proj\!_3}(1))\dual \otimes \OO_G$. According to
[G-P], the scheme of multi-jumping lines $M$ is the locus where
$\sigma$ has rank at most $2n-2$. Denote by $U$ the universal
subbundle of $\proj(H^1E\dual \otimes \OO-{\widetilde{G}})$ and
$\pi$ the projection on $\widetilde{G}$. The degeneracy locus of
$f^*\sigma $ is equal to the one of $f^*\sigma \dual$, and we can
compute $f^*[M]$ in function of the vanishing locus $Z_s$ of some
section $s$ of $U\dual \otimes \pi^{*}(H^1E\dual \otimes
\OO_{\widetilde{G}})$ arising from the sequence:
\begin{center}
$0\longrightarrow U\longrightarrow \pi ^{*}(H^1E\dual \otimes \OO_{
  \widetilde{G}}) \stackrel{\pi ^{*}f^{*}\sigma \dual}{\longrightarrow
  }\pi ^{*}(H^1E(-1) \otimes f^{*}Q\dual)\dual$
\end{center}
Indeed, one has $f^{*}M=\pi _{*}(Z_s)$, and $s$ vanishes on the
divisor $\pi^*\widetilde{S}$ so it gives a regular section
$s^{\prime}$ of $U\dual \otimes \pi ^{*}(H^1E\dual \otimes
\OO_{\widetilde{G}}) (-mx)$ where $S=mS_{red}$ in $A_*G$. Computing
this class as in [Fu] Ex 14.4 gives:
\begin{center}
$Z_{s^{\prime }}=\sum\limits_{i=0}^{2n}(-1)^ic_{2n-i}(U\dual \otimes
\pi ^{*}f^{*}F)(\pi ^{*}mx)^i$, where $F=H^1E(-1)\dual \otimes Q$.
\end{center}
According to Josefiak-Lascoux-Pragacz, [Fu]Ex14.2.2, one has:
\begin{center}
$\pi _{*}(c_{2n-i}(U\dual \otimes \pi^{*}f^{*}F)) =c_{2n-i-(2n-2)+1}
(F-H^1E\dual \otimes \OO_{\widetilde{G}}))=c_{3-i}(F)$
\end{center}
But the Chern polynomial of $F$ is:
\begin{center}
$c_Y(F)=1+(nt)Y+[\binom{n+1}2t^2-nu]Y^2+[2\binom{n+1}3tu]Y^3$
\end{center}
where $t$ is the class of a hyperplane section of $G$, and $u$ is
represented by the lines in a plane. (Let's recall that the Chow ring
of $G$ is generated by $t$ and $u$ with the relations: $t^3=2tu$;
$u^2=t^2u$. Furthermore, $c_2Q=t^2-u$ and $[S]=\alpha (t^2-u) +\beta
u$. So we have: $\pi _{*}(Z_{s^{\prime
}})=\sum\limits_{i=0}^{2n}c_{3-i}(F)(-mx)^i=\sum\limits_{i=0}^3c_{3-i}(F)(-mx)^i$\\
  $\pi
  _{*}(Z_{s^{\prime}})=\big( m^3[c_2N-(c_1N)^2]+m^2(c_1N.n.i^{*}t)-m.i^{*}[\binom{n+1}2t^2-nu],$\hfill \phantom{.}\\
\phantom{.}\hfill $2\binom{n+1}3tu-nm^2t.i_{*}[S_{red}]+m^3i_{*}c_1N \big)$\\
\begin{example}
  \label{exresiduel}
  If $S$ is a smooth congruence of bidegree $(\alpha,\beta)$, then the
  residual class is:
  $\bigl(\frac{\alpha ^2+\beta
^2-(n^2-7n+13)\beta -(n^2-5n+13)\alpha +(2\pi -2)(2n-12)-12\chi
  (\OO_S)} 2\overline{p},$\\
\phantom{.}\hfill${{[2\binom{n+1}3-\!\!(n-3)(\alpha \!+\!\beta )\!+\!2\pi\!\! -\!\!2]tu \bigr)}
  }$\\
where $\pi$ is the genus of a hyperplane section of $S$, and
  $\overline{p}$ is the class of a point in $A_0S$.
\end{example}
NB: the term of the right is still valid when $S$ is just locally
complete intersection, and that is the non vanishing of this term in
the required situations that had been used in the \ref{c2=5courbe}.

Let $c_i\Omega $ be the Chern classes of the cotangent bundle of $S$.
If $C$ is a general hyperplane section of $S$, the normal bundle of
$C$ is $N_C=N_{|C} \oplus \OO_C(1)$, where $N$ is still the normal
bundle of $S$. But $c_1N_C=c_1(\Omega_G\dual)_{| C}+c_1(\omega
_C)=4(\alpha +\beta )+2\pi -2$, so we have: $t.c_1N=[3(\alpha +\beta
)+2\pi -2].\overline{p}$.

On the other hand, one has $c_1N^2-c_2N=[5(\alpha +\beta )+8(\pi
-1)+c_2\Omega ].\overline{p}$, which gives the formula using
Hirzebruch-Riemann-Roch theorem and the following relation satisfied
by every smooth congruence (Cf [A-S]):
\begin{center}
$\alpha ^2+\beta ^2=3(\alpha +\beta )+4(2\pi -2)+2(c_1\Omega )^2-12\chi (\OO_S)$
\end{center}
\subsection{The family of $n$-instantons with $h^0(E(2)) \geq 4$}
  \begin{proposition}\label{famille}
  Let $E$ be any $n$-instanton with $n \geq 5$, $h^0(E(1))=0$ and \hbox{$h^0
 (E(2)) \geq 4$}, then $h^0(E(2))=4$ and there is a section of $E(2)$
 vanishing on the union of $2$ curves of arithmetic genus $0$ cutting
 each other in length 2. Furthermore, one of these curves has to be a
 skew cubic (not necessarily integral), and the other one may be
 chosen smooth and has bidegree $(1,n)$ in a  smooth fixed quadric.
\end{proposition}
We can construct from the hypothesis $h^0 (E(2)) \geq 4$ a map \hbox{$4 \OO
_{\proj_3} \rightarrow E(2)$} whose kernel and cokernel will be
denoted by $E^{\prime}(-2)$ and $\LL$. So we have the exact sequence:
$$0 \longrightarrow E^{\prime}(-2)\longrightarrow 4 \OO _{\proj_3}
\longrightarrow E(2) \longrightarrow \LL\longrightarrow 0$$
The key is that the support of $\LL$ has to contain a
quadric surface. So we have to eliminate in the following the other cases.

\bigskip  
\noindent\underline{First case:} dim $supp \LL \leq 1$
\begin{lemma}
  \label{c2E'}
  If dim $supp \LL \leq 1$, then we have
  $c_2(E^{\prime})=8-n-d_\LL$, where 
  $d_\LL$ is the degree of the sheaf $\Ext^2(\LL,\OO_{\proj_3})$.
\end{lemma}
NB: In this case $d_\LL$ is also the degree of $\LL$, but we'd like
to keep the definition of $d_\LL$ in the other cases.

Let $a$ be such that $E\primedual(-a)$ has a section. Choosing one
enables to build the following diagram where the middle line is obtained by
dualizing the previous exact sequence, and where $X$ is the vanishing
locus of the chosen section ($X$ could be empty):\\
\vbox{

\begin{equation}\label{diag1}
  \begin{picture}(2550,500)(0,500)
    \setsqparms[1`1`1`1;650`250]
    \putsquare(800,250)[4\Op (-2)`E\primedual`
    k\Op(-2)`{\mathcal{J}}_X (-a);```]
    \setsqparms[1`0`1`1;600`250]
    \putsquare(1450,250)[\phantom{E\primedual}`\Ext^2(\LL ,\Op (-2)) `\phantom{{\mathcal{J}}_X (-a)}`A;```]  
    \setsqparms[1`1`1`0;650`250]
    \putsquare(800,500)[(4-k)\Op (-2)`\Op (a)`\phantom{4\Op
    (-2)}`\phantom{E\primedual};```] 
    \setsqparms[1`0`1`0;600`250]
    \putsquare(1450,500)[\phantom{\Op (a)}`B`\phantom{E\primedual}`
    \phantom{\Ext^2(\LL ,\Op (-2))};```]
    \putmorphism(0,500)(1,0)[0`E(-4)`]{300}1a
    \putmorphism(300,500)(1,0)[\phantom{E(-4)}`\phantom{4\Op
    (-2)}`]{500}1a
    \putmorphism(2050,500)(1,0)[\phantom{\Ext^2(\LL
    ,\Op (-2))}`0`]{550}1a
  \putmorphism(800,250)(0,1)[`0`]{250}1l
  \putmorphism(1450,250)(0,1)[`0`]{250}1l
  \putmorphism(2050,750)(1,0)[\phantom{B}`0`]{550}1a
  \putmorphism(2050,250)(1,0)[\phantom{A}`0`]{550}1a
  \putmorphism(2050,250)(0,1)[`0`]{250}1l
  \putmorphism(800,1000)(0,1)[0``]{250}1l
  \putmorphism(1450,1000)(0,1)[0``]{250}1l
    \end{picture}
  \end{equation}
 \vskip 5em}
We can first take $a=-2$. So there arise $k=3$ quartic surfaces
containing $X$, and $B$ has to vanish. The kernel of the last line
of (\ref{diag1}) is in this case $E(-4)$. Choose 2 of those quartics,
and denote by $\Gamma $ their complete intersection, and by $\ideal
_{X|\Gamma}$ the ideal of $X$ in $\Gamma$. The last line of
(\ref{diag1}) and the resolution of $\Gamma$ give in the following a
section $s$ of $E(2)$ whose vanishing locus will be noted $Z$.\\   
\vbox{
\begin{equation}\label{diag2}
  \begin{picture}(2150,550)(50,500)
    \setsqparms[1`1`1`1;500`250]
    \putsquare(300,250)[E(-2)`3\Op `{\mathcal{J}}_Z `\Op;```]
    \setsqparms[1`1`1`0;500`250]
    \putsquare(300,500)[\Op(-4)`2\Op`\phantom{E(-2)} `\phantom{3\Op};`s``]
    \setsqparms[1`0`1`1;500`250]
    \putsquare(800,250)[\phantom{3\Op}`{\mathcal{J}}_X (4)`
    \phantom{\Op}` {\mathcal{J}}_{X|\Gamma} (4);```]
    \setsqparms[1`0`1`0;500`250]
    \putsquare(800,500)[\phantom{2\Op}`{\mathcal{J}}_{\Gamma}
    (4)` ` ; ` ` ` ]
    \setsqparms[1`0`1`1;600`250]
    \putsquare(1300,250)[\phantom{{\mathcal{J}}_{\Gamma}
    (4)}`\Ext^2(\LL ,\Op ) `
    \phantom{{\mathcal{J}}_{X|\Gamma}  (4)}`\Ext^2(\LL
    ,\Op );``\wr`]
    \setsqparms[1`0`1`0;600`250]
    \putsquare(1300,500)[\phantom{{\mathcal{J}}_{\Gamma} (4)}`0``;```]
     \putmorphism(-50,250)(1,0)[0`\phantom{{\mathcal{J}}_Z
     }`]{300}1a
     \putmorphism(-50,500)(1,0)[0`\phantom{E(-2)}`]{350}1a
     \putmorphism(-50,750)(1,0)[0`\phantom{\Op (-4)}`]{350}1a
     \putmorphism(1900,250)(1,0)[\phantom{\Ext^2({
         \mathcal{L}} ,\Op )}` 0`]{450}1a
     \putmorphism(1900,500)(1,0)[\phantom{\Ext^2({
         \mathcal{L}} ,\Op )}` 0`]{450}1a
     \putmorphism(300,250)(0,1)[`0`]{250}1l
     \putmorphism(800,250)(0,1)[`0`]{250}1l
     \putmorphism(1300,250)(0,1)[`0`]{250}1l
     \putmorphism(1900,250)(0,1)[`0`]{250}1l
     \putmorphism(300,1000)(0,1)[0``]{250}1l
     \putmorphism(800,1000)(0,1)[0``]{250}1l
     \putmorphism(1300,1000)(0,1)[0``]{250}1l
  \end{picture}
  \end{equation}
  \vskip 5em}
but the third column of (\ref{diag2}) proves that $X$ is linked by the 2
quartics to the support of $\ideal _{X|\Gamma}$, and the last line
implies that this support is the union of
$Z$ with the 1-dimensional part of \hbox{$\Ext^2(\LL,\Op)$} 's support
counted rank of $\Ext^2(\LL,\Op)$ times. So it
proves lemma \ref{c2E'}. \fin

We can remark now that the degeneracy class of \hbox{$4\Op \rightarrow
  E(2)$} is negative when $n \geq 5$, so the support of $\LL$ has to
  contain at least a one dimensional component. So we have from lemma
  \ref{c2E'} $c_2 E^{\prime}\leq 2$.

\bigskip
Let's use again the diagram (\ref{diag1}), but this time with the
biggest $a$ such that $E\primedual (-a)$ has a section. We will now
study the possible cases recalling that $A$'s support is at most 1
dimensional.

\begin{itemize}
\item If $k=2$

  Then the kernel of the first line of (\ref{diag1}) is $\Op (-a-4)$,
  and it is injected into $E(-4)$. So it gives a section of $E(a)$,
  and by the hypothesis made on $E$ and $a$, we have $a\geq 2$, and
  then $X$ is empty because $A$'s support is at most 1 dimensional. So
  the middle column of (\ref{diag1}) splits. (in this situation we
  will say in the following $E\primedual$ splits).

\item If $E\primedual$ splits
  
  Then the map from $B$ to $\Ext^2(\LL,\Op(-2))$ of (\ref{diag1}) is
  an injection, so $B$'s support has also to be at most 1 dimensional.
  The only possible cases are thus:

\begin{itemize}
\item $k=1,a=2,X=\emptyset$, but it implies  $E\primedual= \Op (2)
  \oplus \Op (-2)$, so it contradicts the independence of the  4
  chosen sections of $E(2)$.

\item $k=2$, then $A$ and $B$'s supports are 1 dimensional of degree
  $(a-2)^2-c_2(E^{\prime})-a^2$ and $(a+2)^2$. As $E\primedual$ is
  splited, the degree of $\Ext^2(\LL,\Op)$ is the sum of $A$ and
  $B$'s degree. Using $c_2 E^{\prime}\leq 3-d_\LL$,
  we obtain the contradiction: $d_\LL\geq a^2+5+d_\LL$.

\item $k=3,a=-2$ but it  is impossible because $E\primedual $ is
  reflexive and \mbox{$c_1(E\primedual)=0$}; $c_2(E\primedual ) \leq 2$, so
  $E\primedual (1)$ must have a section. 
\end{itemize}

So the only remaining case is:
\item $k\geq 3,X\neq\emptyset$

The curve $X$ can't be in 3 independent planes, so $a \geq 0$. In
other words, it means that $E^{\prime}$ is semi-stable.

\begin{itemize}
\item If $a=0$, then $X$ has degree $c_2(E^{\prime})$, which must be 1
  or 2. So there is a plane $H$ having a curve in its intersection
  with $X$. This intersection is in fact the union of a curve of
  degree $b=1$ or $2$ with a scheme $X^{\prime}$ which is the
  vanishing of a section of $E_H(-b)$. The last line of (\ref{diag1})
  has  kernel $E(-4)$ because $k\geq 3$, so it gives the following
  exact sequence:
    $$0\longrightarrow E_H(-4) \longrightarrow k \OO_H(-2) \longrightarrow
     {\mathcal{J}}_{X^{\prime}}(-b)\longrightarrow A^{\prime
        }\longrightarrow 0$$
  But $X^{\prime}$ has degree $c_2(E^{\prime})+b^2\geq 2$, so it lies
      in only one line, thus those 3 sections have to be proportional.
      Then we would have $E_H(-4)=2\OO_H(-2)$, which is impossible
      because an instanton doesn't have unstable planes.

    \item If $a=-1$, (i.e: $E\primedual $ stable)

      Then we have $c_2(E\primedual )>0$, and on the other hand,
      $\Ext^2(\LL,\Op)$ has a non empty 1 dimensional component, so we
      have $5\leq n \leq 7$ and $0<c_2(E\primedual ) \leq 2$, and
      $d_\LL \leq 2$. Take this time a general plane $H$ such that
      $E_H$ and $E\primedual _H$ are stable, and
      $d_{\LL}=h^0\Ext^2(\LL,\Op)$. Restrict to $H$ the second line of
      diagram (\ref{diag1}) to obtain:
    $$ 0\rightarrow E_H(-4)\rightarrow 4 \OO _H(-2)
    \kercoker{N} E\primedual_H \rightarrow
    \Ext^2(\LL,\Op(-2)) _H \rightarrow 0$$

    But $h^1E\primedual_H=0$ because $c_2(E\primedual ) \leq 2$, and
    $h^0N(1)=h^1E_H(-3)=h^0E_H=0$, so we have \hbox{$h^0\Ext^2(\LL,\Op(-1)) _H \geq
    h^0E\primedual _H(1)$}, which is greater or equal to 4, so it
    contradicts $d_\LL \leq 2$. 
  \end{itemize}
\end{itemize}

\bigskip
\noindent\underline{Second case:} The support of $\LL$ contain a 2
dimensional part which is a plane $H$.

We have the following  sequence, where this time $c_1(E\primedual)=-1$.
$$0 \longrightarrow E^{\prime}(-2)\longrightarrow 4 \OO _{\proj_3}
\kercoker{M} E(2) \longrightarrow \LL\longrightarrow 0$$
As previously, we obtain the following diagram:\\
\vbox{ 
\begin{equation}\label{diag1b}
  \begin{picture}(2600,550)(0,500)
    \setsqparms[1`1`1`1;650`250]
    \putsquare(800,250)[4\Op (-2)`E\primedual`
    k\Op(-2)`{\mathcal{J}}_X (-a-1);```]
    \setsqparms[1`0`1`1;600`250]
    \putsquare(1450,250)[\phantom{E\primedual}`\Ext^2(\LL ,\Op (-2)) `\phantom{{\mathcal{J}}_X (-a-1)}`A;```]  
    \setsqparms[1`1`1`0;650`250]
    \putsquare(800,500)[(4-k)\Op (-2)`\Op (a)`\phantom{4\Op
    (-2)}`\phantom{E\primedual};```] 
    \setsqparms[1`0`1`0;600`250]
    \putsquare(1450,500)[\phantom{\Op (a)}`B`\phantom{E\primedual}`
    \phantom{\Ext^2(\LL ,\Op (-2))};```]
    \putmorphism(-100,500)(1,0)[0`M\dual(-2)`]{350}1a
    \putmorphism(250,500)(1,0)[\phantom{M\dual(-2)}`\phantom{4\Op
    (-2)}`]{550}1a
    \putmorphism(2050,500)(1,0)[\phantom{\Ext^2(\LL
    ,\Op (-2))}`0`]{550}1a
  \putmorphism(800,250)(0,1)[`0`]{250}1l
  \putmorphism(1450,250)(0,1)[`0`]{250}1l
  \putmorphism(2050,750)(1,0)[\phantom{B}`0`]{550}1a
  \putmorphism(2050,250)(1,0)[\phantom{A}`0`]{550}1a
  \putmorphism(2050,250)(0,1)[`0`]{250}1l
   \putmorphism(800,1000)(0,1)[0``]{250}1l
  \putmorphism(1450,1000)(0,1)[0``]{250}1l
   \end{picture}
\end{equation}
\vskip 5em}

Working as in the first case, we prove that the union of a section of
$E\primedual (2)$ with $\Ext^2(\LL,\Op)$'s support (counted rank of
$\Ext^2(\LL,\Op)$ times) is linked by 2
cubic surfaces (sections of $\stackrel{2}{\Lambda}(E\primedual(2))$)
to a section of $M\dual(3)$. Furthermore, this section is linked by a
section of $E(2)$ to a (possibly empty) curve of $H$. But this plane
curve is of degree at most 2 because $E_H$ is necessarily semi-stable.
 Then, $c_2(M\dual(3))=c_2(M\dual)\geq n+2$ ($n\geq 5$), so we have
 $d_\LL+c_2(E\primedual)\leq 0$  and $d_\LL\leq 2$.

But as $c_1(E\primedual)=-1$ and $c_2(E\primedual)\leq 0$, we must have 
$h^0(E\primedual)\neq 0$, we can assume  $a\geq 0$ in diagram  (\ref{diag1b}).

\begin{itemize}
\item If $X=\emptyset$
  Then $E\primedual$ has to split, so $B$ injects itself in
  $\Ext^2(\LL,\Op (-2))$, thus $B$'s support has to be at most 1
  dimensional. As $A$ has also a 1 dimensional support (at most), we
  have only the following cases:

 \begin{itemize}
  \item $k=1,a=1$, so $E^{\prime}=\Op(-1)\oplus \Op(2)$, but it conflicts
    with the independence of the 4 chosen sections of  $E(2)$.

  \item $k=2$, then $B$ has a 1 dimensional support of degree
    \hbox{$(2+a)^2\geq 4$}, but it conflicts with  $d_\LL\leq 2$.    
  \end{itemize}

\item If $X\neq\emptyset$

  The facts that $A$ has an at most 1 dimensional support, and that
  $a\geq 0$ imply $k=2, X=\proj_1 ,a=0$. Then $X$ must have degree
  $c_2(E\primedual)$, which contradicts $c_2(E\primedual)\leq 0$.
\end{itemize}
 
\noindent\underline{The effective situation:}

We can conclude from the previous cases that there is a quadric $Q$ in $\LL$'s support
because it can't contain a cubic surface as there is no curve in 3
independent planes. So a section $C$ of $E(2)$ have a component $C_2$ in
the quadric $Q$, and another one $C_1$, which has to be skew because
$h^0E(1)=0$, and which lies on 3 quadrics. This curve $C_1$ must then
be a cubic curve of arithmetic genus 0.

Let's remark that $C_2$ has degree $n+1 \geq 6$, so $Q$ is normal
because $C$ can't have plane components of degree 3 or more otherwise
there would be an unstable planes, which is not possible for an instanton.

One has first to check that $C_1$ and $C_2$ have a 0-dimensional
intersection to get the \ref{famille}. The cubics  arising
when the section of $E(2)$ moves, are the vanishing locus of sections
of $M\dual(2)$, which is reflexive of rank 2. Take  a normal quadric
$Q^{\prime}$ containing $C_1$, then it must contain a pencil of
those cubics say $\lambda C_1+\mu C^{\prime}_1$. But if the associated
pencil of sections of $E(2)$ was 
such that $C_1 \cap C_2$ and $C^\prime_1 \cap C^\prime_2$ were 1
dimensional, then \hbox{$C_1 \cap C^{\prime}_1$} would have a 1
dimensional component, because $C_1 \cap C_2$ and $C^\prime_1 \cap
C^\prime_2$ lies in the fixed curve $Q \cap Q^{\prime}$, and this
1-dimensional component of $C_1 \cap C^{\prime}_1$ would be in the
singular locus of $Q^{\prime}$ which contradicts its normality.

So there is a section of $E(2)$ vanishing on $C=C_1 \cup C_2$, where
$C_1$ is a skew cubic of arithmetic genus 0 and $\dim C_1 \cap C_2=0$
. On the other hand, we 
have $\omega_C=\OO_C$, and the exact sequence of liaison (i=1 or 2):
$$0\longrightarrow \omega_{C_i} \longrightarrow \OO_C
\longrightarrow \OO_{C_{3-i}} \longrightarrow 0$$
gives when i=1 by restriction to $C_1$ that $C_1 \cap C_2$ has length
2, and when i=2 by
restriction to $C_2$ that $C_2$ has arithmetic genus 0. But a quadric
cone can't have curves of arithmetic genus 0 and degree greater or
equal to 4, so $Q$ is smooth and $C_2$ has bidegree $(1,n)$ in $Q$.

We'd like now to prove the smoothness of a general $C_2$. Denote by
$h$ the 3-dimensional subspace of $|C_2|$ induced by 
the 4 sections of $E(2)$. The base points of $h$ must be in the singular locus of
$Q\cup Q^{\prime}$ for any quadric $Q^{\prime}$ containing $C_1$,
hence those points are on $C_1$. But we showed that $C_1 \cap C_2$
was 0-dimensional, so the set of base points of $h$ is at most finite.
Furthermore, if $C_2$ is singular in some point $P$, then it must contain the
ruling of bidegree $(0,1)$ passing through $P$, so this ruling would
be a base curve of $h$ if the general curve was singular at $P$.
Hence, the generic element of $h$ is 
smooth and irreducible of bidegree $(1,n)$, and this is the
proposition \ref{famille}. \fin 

\begin{proposition}\label{multisautfamille}
   Let $E$ be a $n$-instanton with $n \geq 5$, $h^0(E(1))=0$ and \hbox{$h^0
 (E(2)) \geq 4$}, then $E$ has only a 1-dimensional scheme of
 multi-jumping lines.
\end{proposition}

The bundle $E$ belongs by hypothesis to the family \ref{famille}, so
denote by $Q$ the quadric which is the support of the cokernel of the map
given by the 4 sections to $E(2)$. Let $s$ be a section of $E(2)$ and
$Z_s=C_{2,s}\cup C_{1,s}$ its vanishing locus, where $C_{1,s}$ is the
rational cubic and $C_{2,s}$ is the curve of bidegree $(1,n)$ in $Q$.
If a line $d$ is a multi-jumping line then $E_d(2)=\OO_d(a+2)\oplus
\OO_d(-a+2)$ with $a\geq 2$. So a multi-jumping line $d$ meets $Z_s$ if
and only if $d\cap Z_s$ has length at least 4. So if $E$ has a 2
parameter family of multi-jumping lines, then infinitely many of them
would be 4-secant to $Z_s$ and up to a change of the section $s$, we can
assume that infinitely many  are not in $Q$. So those lines are
2-secant to $C_{2,s}$ and 2-secant to $C_{1,s}$ because they are not
in $Q$ and $C_{1,s}$ don't have trisecant. Denote by $\Sigma$ the
ruled surface of $\proj_3$ made with those lines.

a) If $C_{1,s}$ is smooth.\\
Then the basis of $\Sigma$ is a curve $\Gamma$ on the Veronese surface of
bisecant lines to $C_{1,s}$, and $\Sigma$ must contain $C_2$ because $C_2$
is irreducible. But we have a morphism of degree 2 from $C_2$ onto the
basis of $\Sigma$ hence $\Gamma$ is smooth and rational, so $\deg
\Sigma=2$ or $\deg \Sigma=4$, but $\Sigma$ contains $C_1\cup C_2$
so $\deg \Sigma \geq 4$ because $h^0E(1)=0$, and all the quartics
containing $C_1\cup C_2$ are some $Q\cup Q^{\prime}$, where
$Q^{\prime}$ is a quadric containing $C_1$, so $\Sigma$ can't be a
quartic containing $C_1\cup C_2$ and this case is impossible.

b) If $C_{1,s}$ is made of 3 lines containing a same point $N$
independent of $s$, then $C_{2,s}$ must contain $N$ because $Z_s$ is
locally complete intersection. As $C_{2,s}$ is smooth at $N$, any line
through $N$ meets 
$Z_s$ in length 2 around $N$, and $N\in Q$, so a line through $N$
don't meet $C_{1,s}$ in another point, and it
can't meet $C_{2,s}$ in 2 points distinct of $N$, so it can't be
4-secant to $C_{2,s}\cap C_{1,s}$. There must then
exist a plane $H$ containing 2 of the lines of $C_{1,s}$ and
infinitely many  multi-jumping lines meeting $C_{2,s}$. As
$h^0E_H(-1)=0$, $C_{2,s}\cap H$ must be 0-dimensional, and there would
be a point $P$ in $C_{2,s}\cap H$ such that every line of $H$ through $P$ is
bisecant to $C_{2,s}$. Hence 
$H$ is tangent to $Q$ at $P$, but one of the ruling of $Q$ through $P$
would be in $H$ and would be only 1 secant to $C_{2,s}$ because it has bidegree
$(1,n)$, which contradicts the definition of $P$. \fin

\begin{proposition}\label{lisse}
   Let $E$ be a $n$-instanton with $n \geq 5$, $h^0(E(1))=0$ and \hbox{$h^0
 (E(2)) \geq 4$}, then $E$ is a smooth point of the moduli space I$_n$.
\end{proposition}

From classical theory (cf [LeP]), one has to prove that $h^2(E\otimes
E)=0$, and if $C$ is the vanishing locus of a section of $E(2)$, then
$h^1(E_C(2))=0$ implies this result. As $E_C(2)$ is just the normal
bundle of $C$, we are considering the problem of the vanishing of
$h^1(N_C)$, so we will solve it as it was done for smoothing questions
in [H-H].

\bigskip
Let's recall from \ref{famille} that $C=C_1\cup C_2 $ where the $C_i$
have zero arithmetic genus, and where $C_1$ is a skew cubic and $C_2$
is on a smooth quadric $Q$. We want first to prove the vanishing of
$h^1(N_{C_i})$. As $C_2$ has bidegree $(1,n)$ on the smooth quadric
$Q$, we have $h^1(N_{C_2})=0$. (Cf for example the proof of prop 5.4
$\alpha$ in [H-H]). This will be true for $C_1$ even if it is 3
concurrent lines. Indeed $C_1$ is the vanishing of a section of $R(1)$
where $R$ is a rank 2 reflexive sheaf with $(c_1,c_2,c_3)=(0,2,4)$,
and the following exact sequences:
\begin{center}
$0 \rightarrow 2\Op(-1)\rightarrow 4\Op \rightarrow R(1)\rightarrow 0$
\&\nolinebreak  $0 \rightarrow 2\Op(-3)\rightarrow3\Op(-2)\rightarrow
\ideal_{C_1}\rightarrow 0$ 
\end{center}
imply that $h^1(R(1))=0$, $h^2(\ideal_{C_1}R(1))=0$ thus
$h^1(R_{C_1}(1))=0$, but \hbox{$R_{C_1}(1)=N_{C_1}$}.

So both $h^1N_{C_i}$ are zero, and we can now deduce from it that
$h^1N_C$ is also zero. We have the following sequence of liaison:
$$0\longrightarrow \OO_C \longrightarrow \OO_{C_1}\oplus \OO_{C_2}
\longrightarrow \OO_{C_1 \cap C_2} \longrightarrow 0$$
giving when twisted by $E(2)$:
$$0\longrightarrow N_C \longrightarrow E_{C_1}(2)\oplus E_{C_2}(2)
\longrightarrow E_{C_1 \cap C_2}(2) \longrightarrow 0$$
Furthermore, the inclusion ${\mathcal{J}}_C \subset
{\mathcal{J}}_{C_i}$ gives a map \hbox{$N_{C_i} \rightarrow N_C$}.
Then we have the following exact sequence defining $\LL_i$ and $M_i$:
$$ ? \longrightarrow N_{C_i} \kercoker{M_i} E_{C_i}(2)
\longrightarrow \LL_i \longrightarrow 0$$
As $h^1(N_{C_i})=0$ implies $h^1(M_i)=0$,
we have $
\left\{\begin{array}{l}
      h^1(E_{C_i}(2))=0\\
      H^0(E_{C_i}(2))\rightarrow H^0(\LL_i) \rightarrow 0
    \end{array}\right.$, which gives $h^1(N_C)=0$ when taking global
  sections in the following diagram.
\begin{center}
    \begin{picture}(2400,500)
      \setsqparms[1`1`1`1;800`250]
      \putsquare(300,250)[N_{C_1}\oplus N_{C_2}`E_{C_1}(2)\oplus
      E_{C_2}(2)`N_C`E_{C_1}(2)\oplus E_{C_2}(2);``\wr`]
      \setsqparms[1`0`1`1;800`250]
      \putsquare(1100,250)[\phantom{E_{C_1}(2)\oplus
      E_{C_2}(2)}`\LL_1 \oplus \LL_2
    `\phantom{E_{C_1}(2)  \oplus E_{C_2}(2)} `E_{C_1 \cap C_2}(2) ;```]
    \putmorphism(0,250)(1,0)[0`\phantom{N_C}`]{300}1a
    \putmorphism(1900,250)(1,0)[\phantom{E_{C_1 \cap C_2}(2)}`0`]{400}1a
    \putmorphism(1900,500)(1,0)[\phantom{\LL_1 \oplus
      \LL_2} `0`]{400}1a
    \putmorphism(300,250)(0,1)[`0`]{250}1l
    \putmorphism(1900,250)(0,1)[`0`]{250}1l
   \end{picture}
 \end{center}

We can now conclude that $E$ is a smooth point of the moduli space
because $N_C=E_C(2)$, and the previous vanishing gives $h^2(E\otimes
E)=0$ using the sequences:
\begin{center}
 $0 \rightarrow E(-2)\rightarrow E\otimes E\rightarrow \ideal_C
E(2)\rightarrow 0$ and $0 \rightarrow \ideal_C E(2)\rightarrow
E(2)\rightarrow N_C\rightarrow 0$  
\end{center}
\subsection{A $\theta$-characteristic on the curve of multi-jumping
  lines of an $n$-instanton}                               

Assume here that $E$ is an instanton vector bundle with second Chern
class $n$. The aim of this part is to study the scheme of multi-jumping
lines of $E$ when it satisfies the properties expected from its
determinantal structure.

For instance, according to [B-H], the general member of the
irreducible component of the 
moduli space containing the t'Hooft bundles satisfy those properties.
\begin{proposition}
\label{w2}If an $n$-instanton has only \prefix{k\leq 2}jumping lines,
and if its scheme of multi-jumping lines is a curve $\Gamma $ in $G$,
then we have:\\\centerline{$(R^1q_{*}p^{*}E)^{%
\otimes 4}=\OO_\Gamma (n);\omega _\Gamma
^o=(R^1q_{*}p^{*}E)^{\otimes 2}(n-4);(\omega _\Gamma ^o)^{\otimes 2}=%
\OO_\Gamma (3n-8)$}\\where $\omega _\Gamma ^o$ is a dualizing sheaf on  $\Gamma $.
\end{proposition}
In fact this result is just an analogous of the proposition 1.5 of
[E-S]. Let $K$ and $Q$ be the tautological bundles on $G$ such that we
have the sequence:
\begin{center}
$0\longrightarrow K\longrightarrow H^0(\OO_{\proj_3}(1))\dual\otimes
\OO_G\longrightarrow Q\longrightarrow 0$ 
\end{center}
The points/lines incidence variety of $\proj_3$ is
$F=\proj_G(Q\dual)$, so one has the exact sequence, where $q$ still be
the projection from $F$ to $G$:
\begin{center}
$0\longrightarrow \OO_F(-\tau)\longrightarrow q^*Q\dual
\longrightarrow (\omega_{F/G}(\tau ))\dual \longrightarrow 0$ 
\end{center}
We'd like here to make a relative (over $G$) Beilinson's construction.
So consider the resolution of the diagonal $\Delta$ of
$F\stackunder{G}{\times}F$:
\begin{center}
$0\longrightarrow \OO_F(-\tau )\stackunder{G}{\boxtimes }\omega
_{F/G}(\tau )\longrightarrow \OO_{F\stackunder{G}{\times}%
F}\longrightarrow \OO_\Delta \longrightarrow 0$
\end{center}
which gives when twisted by $p^{*}E(\tau )\stackunder{G}{ \boxtimes}
\OO_F$ the following spectral sequence which stops in $E_2^{p,q}$
and ends to $p^{*}E(\tau )$ in degree 0 and 0 in the other degrees.
\begin{center}
{\begin{picture}(2000,500)(-300,0)
\setsqparms[1`0`0`1;1500`300]
\putsquare(0,0)[{q^{*}(R^1q_{*}p^{*}E)\otimes \OO_F(\sigma -\tau 
)}`{q^{*}(R^1q_{*}p^{*}E(\tau ))}`{q^{*}q_{*}p^{*}E\otimes 
\OO_F(\sigma -\tau )}`{q^{*}q_{*}p^{*}E(\tau )};```d]
\putmorphism(1500,0)(1,0)[\phantom{q^{*}q_{*}p^{*}E(\tau)}``p]{500}1a
\putmorphism(1500,600)(0,1)[`\phantom{M}`q]{300}{-1}r
\put(1500,70){\line(0,1){150}}
\end{picture}}
\end{center}
By assumption $E$ has no \prefix{k\geq 3}jumping lines, so
$R^1q_{*}p^{*}E(\tau )=0$, and we have the following surjection where
$M$ denotes its kernel.
\begin{center}
$0\longrightarrow M\longrightarrow p^{*}E(\tau )\longrightarrow
q^{*}(R^1q_{*}p^{*}E)\otimes \OO_F(\sigma -\tau )\longrightarrow 0$
\end{center}
Denote by $h$ and $K$ the restrictions of $q$ and
$q^{*}(R^1q_{*}p^{*}E)\otimes \OO_F(\sigma -\tau )$  to
$q^{-1}(\Gamma)$. The sheaf $R^1q_{*}p^{*}E$ is locally free over
$\Gamma$, so $h_{*}K=R^1q_{*}p^{*}E\otimes h_{*}\OO_F(\sigma -\tau
)=0$, and we have: $h_{*}M_{| q^{-1}(\Gamma )}=h_{*}(p^{*}E(\tau
)_{|q^{-1} (\Gamma )})$. Furthermore, as $K$ is locally free over
$\Gamma$ and $\stackrel{2}{\Lambda }E=0$, one has $M_{| q^{-1}(\Gamma
)}=K\dual(2\tau )$, so $h_{*}(p^{*}E(\tau )_{| q^{-1}(\Gamma
)})=(R^1q_{*}p^{*}E)\dual \otimes h_{*}\OO_F(3\tau -\sigma )$. But it
means that  $h_{*}(p^{*}E(\tau )_{| q^{-1}(\Gamma )})\otimes
R^1q_{*}p^{*}E=\OO_G(-1)\otimes Sym_3Q_{| \Gamma }$, and using that
$h_{*}(p^{*}E(\tau )_{| q^{-1}(\Gamma )})$ is locally free of rank 4,
we found the relation: $(R^1q_{*}p^{*}E)^{\otimes
4}\otimes \det (h_{*}(p^{*}E(\tau )_{| q^{-1}(\Gamma )}))=\OO_\Gamma
(2)$. But in fact $q_*p^*E(\tau)$ is locally free over $G$ by base
change, so $\det(h_*(p^*E(\tau)_{|q^{-1}(\Gamma)}))= \OO_\Gamma
(2-n)$, because $c_1(q_{*}(p^{*}E(\tau )))=2-n$ from Riemann-Roch over
$F$. So we have:
\begin{center}
$(R^1q_{*}p^{*}E)^{\otimes 4}=\OO_\Gamma (n)$
\end{center}
On another hand, the following resolution of $\OO_F$ in $G\times
\proj_3$:
\begin{center}
$0\longrightarrow \OO_{G\times \proj_3}(-\sigma -2\tau
)\longrightarrow Q\dual \otimes \OO_{\proj_3}(-\tau
)\longrightarrow \OO_{G\times \proj_3}\longrightarrow
\OO_F\longrightarrow 0$ 
\end{center}
gives when twisted by $p^*E$ the following exact sequence deduced from
the Leray spectral sequence:
\begin{center}
$H^1(E(-1))\otimes Q\dual \stackrel{\phi }{\longrightarrow }
H^1(E)\otimes \OO_G\longrightarrow  R^1q_{*}p^{*}E\longrightarrow 0$
\end{center}
The Eagon-Northcott complexes $(E_i)$ associated to $\phi$ gives
resolutions of $M_i$, where we have because $\Gamma$ is a curve:
$M_0=\OO_\Gamma $, $M_1=R^1q_{*}p^{*}E$, $M_i=Sym_iM_1$, and furthermore:
\begin{center}
$\omega _\Gamma ^o=M_2\otimes \omega _G\otimes \det (H^1(E)\otimes
\OO_G) \otimes (H^1(E(-1))\otimes Q\dual)\dual$
\end{center}
But $\omega _G=\OO_G(-4)$, and $\det (H^1(E)\otimes \OO_G) \otimes
(H^1(E(-1)) \otimes Q\dual)\dual=\OO_G(n)$, so we have:
\begin{center}
$\omega _\Gamma ^o=(R^1q_{*}p^{*}E)^{\otimes 2}(n-4)$ ans $(\omega
_\Gamma ^o)^{\otimes 2}=\OO_\Gamma (3n-8)$
\end{center}
\begin{remark}
$\theta =(R^1q_{*}p^{*}E)\otimes \omega _\Gamma ^o(2-n)$ is a
$\theta$-characteristic of $\Gamma $.\\If $n=2n^{\prime }$, then so is
$\theta=(R^1q_{*}p^{*}E)(n^{\prime }-2)$.
\end{remark}
\begin{proposition}
If the scheme of multi-jumping lines $\Gamma $ of an $n$-instanton is a
curve in $G$ without \prefix{k\geq 3}jumping lines, then:\\
\centerline{$R^1q_{*}p^{*}E(-\tau )_{| \Gamma }\simeq
  (R^1q_{*}p^{*}E)\otimes Q_{| \Gamma }(-\sigma )$}\\Furthermore, if
$\Gamma $ is smooth, then its normal bundle in $G$ is:\\\centerline{$N_{\Gamma
,G}=Sym_2Q\otimes \omega _\Gamma (3-n)$}
\end{proposition}
This proposition is also due to a relative Beilinson's construction,
but this time the resolution of the diagonal is twisted by
$p^{*}E\stackunder{G}{\boxtimes }\OO_F$, so we have the spectral
sequence:
\begin{center}
{\begin{picture}(2000,500)(-400,0)
\setsqparms[1`0`0`1;1500`300]
\putsquare(0,0)[{q^{*}(R^1q_{*}p^{*}E(- \tau ))\otimes \OO_F(\sigma 
-\tau )}`{q^{*}(R^1q_{*}p^{*}E)}`{q^{*}(q_{*}p^{*}E(- \tau ))\otimes 
\OO_F(\sigma -\tau )}`{q^{*}q_{*}p^{*}E};d```d']
\putmorphism(1500,0)(1,0)[\phantom{q^{*}q_{*}p^{*}E(\tau)}``p]{500}1a
\putmorphism(1500,600)(0,1)[`\phantom{M}`q]{300}{-1}r
\put(1500,70){\line(0,1){150}}
\end{picture}}
\end{center}
which ends with $0\rightarrow K\rightarrow p^{*}E\rightarrow
N\rightarrow 0$ where $N$ is the kernel of $d$ and $K$ the cokernel of
$d^{\prime}$. It gives the exact sequences (*):
\begin{center}
$0\longrightarrow q^{*}(q_{*}p^{*}E(-\tau ))\otimes \OO_F(\sigma
-\tau )\longrightarrow q^{*}(q_{*}p^{*}E)\kercoker{K} p^{*}E
\longrightarrow N \rightarrow 0$

$0\longrightarrow N(\tau -\sigma )\longrightarrow q^{*}(R^1q_{*}p^{*}E(-\tau
))\longrightarrow q^{*}(R^1q_{*}p^{*}E)(\tau -\sigma )\longrightarrow 0$
\end{center}
Let's restrict this last sequence to $q^{-1}(\Gamma)$, then we can
  apply the projection formula to $q^{*}(R^1q_{*}p^{*}E(-\tau ))_{| 
  q^{-1}\Gamma}$  and to $q^{*}(R^1q_{*}p^{*}E)(\tau -\sigma )_{|
  q^{-1}\Gamma }$ because they are locally free, and we obtain the exact
sequence:
\begin{center}
$R^1q_{*}p^{*}E(-\tau )_{| \Gamma }\longrightarrow R^1q_{*}p^{*}E\otimes
q_{*}\OO_\Gamma (\tau -\sigma )\longrightarrow R^1q_{*}N_{|
  q^{-1}\Gamma }(\tau -\sigma )$ 
\end{center}
On another hand the restriction to $q^{-1}\Gamma$ of the first
sequence of (*) shows that $R^1q_{*}K_{| \Gamma }(\tau )=0$, so that
$R^1q_{*}N_{| q^{-1}\Gamma }(\tau )\simeq R^1q_{*}p^{*}E_\Gamma
(\tau )=0$ because $E$ has no \prefix{k \geq 3}jumping lines. So we
have a surjection from $R^1q_{*}p^{*}E(-\tau )_{| \Gamma }$ to 
$R^1q_{*}p^{*}E\otimes Q(-1)$ which is in fact an isomorphism because
those 2 sheaves are locally free of rank 2 on $\Gamma$. Furthermore,
the curve $\Gamma$ is given by the first Fitting ideal of the
following map $\phi$:
\begin{center}
$H^2(E(-3)\otimes \OO_G(-1)\stackrel{\phi }{\longrightarrow }%
H^1(E(-1))\otimes \OO_G \stackrel{j }{\longrightarrow} R^1q_{*}p^{*}E(-\tau
)\longrightarrow 0$
\end{center}
where $\phi$ is symetric with respect to Serre's duality. Assume now
that $\Gamma$ is smooth, and denote by $\LL_{\Gamma}$ the restriction
of $R^1q_{*}p^{*}E(-1)$ to $\Gamma$.

We just need now, to conclude the proof, to recall the results of
Tjurin (Cf [T1]) which give a description of the normal bundle.
Let's consider the hypernet of quadrics classically associated to an
instanton, and denote by $V=(H^0\OO_{\proj_3}(1))\dual$ and
$H=H^2E(-3)$. The hypernet gives an inclusion
$\proj(\stackrel{2}{\Lambda } V) \subset \proj(Sym_2H\dual)$, and
denote by $D_2$ the surface made of the quadrics of the hypernet of
rank at most $n-2$, assuming here that $G$ and
$\proj(\stackrel{2}{\Lambda } V)$ cut transversally the stratification
of $\proj(Sym_2H\dual)$ by the rank of the quadrics. The second
symetric power of the restriction
of $j$ to $\Gamma$ gives a surjection  $Sym_2H\dual
\otimes \OO_\Gamma \stackrel{s_2j_{| \Gamma }}{\rightarrow }Sym_2{\cal L}_\Gamma
\rightarrow 0$. Denote by $K$ the kernel of $s_2j_{| \Gamma }$, then
the fiber $\proj_(K_q)$ over some $q\in \Gamma$ is just made of the
quadrics of $\proj(H)$ which contain the singular locus of $q$. This
is also according to [T1] or [T3]\S2 lemma 1.1, the projective tangent
space in $q$ to the locus of rank at most n-2 quadrics. The above hypothesis
 of transversality implies that $\stackrel{2}{\Lambda }V\cap K_q$ is 3
 dimensional. Similarly, composing the inclusion
 $(\stackrel{2}{\Lambda }V)\otimes \OO_\Gamma \subset Sym_2H^{*}\otimes 
 \OO_\Gamma $ with $s_{2j_\Gamma}$ gives the sequence:
\begin{center}
$0\longrightarrow K^{\prime }\longrightarrow (\stackrel{2}{\Lambda }%
V)\otimes \OO_\Gamma \longrightarrow Sym_2{\cal L}_\Gamma
\longrightarrow 0$
\end{center}
where $\proj(K_q^{\prime })$ is the projective tangent space to $D_2$
at a point $q$ of $\Gamma$. So we have the commutative diagram of [T1]
restricted to $\Gamma$:
\begin{center}
{\begin{picture}(2200,1000)
\setsqparms[1`-1`-1`1;700`250]
\putsquare(400,250)[{K^{\prime }}`{\buildrel\hbox{\tiny2}\over{\Lambda} 
V\otimes \OO_\Gamma }`{\OO_\Gamma 
(-1)}`{\OO_\Gamma (-1)};```]
\putsquare(400,500)[{T_{D_2}(-1)_{\mid \Gamma }}`{T_{\proj\!_5}(-1)_{\mid 
\Gamma }}`\phantom{K^{\prime 
}}`\phantom{\buildrel\hbox{\tiny2}\over{\Lambda} V\otimes 
\OO_\Gamma };```]
\putsquare(1100,500)[\phantom{T_{\proj\!_5}(-1)_{\mid \Gamma 
}}`{N_{D_2,\proj\!_5}(-1)_{\mid \Gamma 
}}`\phantom{\buildrel\hbox{\tiny2}\over{\Lambda} V\otimes 
\OO_\Gamma }`{Sym_2{\cal L}_\Gamma};```]
\putmorphism(0,250)(1,0)[{0}`\phantom{\OO_\Gamma (-1)}`]{400}1l
\putmorphism(0,500)(1,0)[{0}`\phantom{K^{\prime }}`]{400}1l
\putmorphism(0,750)(1,0)[{0}`\phantom{T_{D_2}(-1)_{\mid \Gamma }}`]{400}1l
\putmorphism(1100,250)(1,0)[\phantom{\OO_\Gamma (-1)}`{0}`]{700}1l
\putmorphism(1800,500)(1,0)[\phantom{Sym_2{\cal L}_\Gamma }`{0}`]{500}1l
\putmorphism(1800,750)(1,0)[\phantom{N_{D_2,\proj\!_5}(-1)_{\mid \Gamma 
}}`{0}`]{500}1l
\putmorphism(1100,1000)(0,1)[{0}`\phantom{T_{\proj\!_5}(-1)_{\mid \Gamma 
}}`]{250}{-1}m
\putmorphism(1800,1000)(0,1)[{0}`\phantom{N_{D_2,\proj\!_5}(-1)_{\mid 
\Gamma }}`]{250}{-1}m
\putmorphism(400,250)(0,1)[\phantom{\OO_\Gamma (-1)}`{0}`]{250}{-1}m
\putmorphism(1100,250)(0,1)[\phantom{\OO_\Gamma (-1)}`{0}`]{250}{-1}m

\putmorphism(1800,500)(0,1)[\phantom{Sym_2{\cal L}_\Gamma}`\phantom{0}`]{250
}{-1}m
\putmorphism(400,1000)(0,1)[{0}`\phantom{T_{D_2}(-1)_{\mid \Gamma 
}}`]{250}{-1}m
\end{picture}
}
\end{center}
where the first two columns are the Euler relative exact sequences over
$\Gamma$ and \mbox{$\proj(\stackrel{2}{\Lambda }V)$}. But $\Gamma= D_2
\cap G$, so we have
$N_{\Gamma,G}(-1)=N_{D_2,\proj_5}(-1)_{|\Gamma}$, then \linebreak$N_{\Gamma
  ,G}\simeq Sym_2(R^1q_{*}p^{*}E(-\tau )_{| \Gamma })\otimes
\OO_\Gamma (1)$, which gives with the previous results:
\begin{center}
  $N_{\Gamma ,G}\simeq Sym_2Q\otimes \omega _\Gamma (3-n)$
\end{center}

\begin{center}
   \hrule
\end{center}
{\small F.Han: Universit\'{e} Denis Diderot (Paris 7), Institut de
  Math\'{e}matiques de Jussieu, 2 place Jussieu, 75251 Paris cedex 05, 
France.\\ email: han@math.jussieu.fr}}
\end{document}